\documentclass[12pt,a4paper,twoside]{article}
\usepackage[french] {babel}
\usepackage{amssymb,amsmath,amsfonts,graphics}
\usepackage{stmaryrd}
\usepackage[all]{xy}

\begin{document}

\newcommand{\hooklongrightarrow}{\lhook\joinrel\longrightarrow}
\begin{center}
\textbf{Fonctions $L$ en g\'eom\'etrie rigide I : $F$-modules convergents ou surconvergents et conjecture de Dwork\ }

\vskip20mm

Jean-Yves ETESSE
  \footnote{(CNRS - IRMAR, Universit\'e de Rennes 1, Campus de Beaulieu - 35042 RENNES Cedex France)\\
E-mail : Jean-Yves.Etesse@univ-rennes1.fr}
\end{center}
 
 \vskip60mm
 \noindent\textbf{Sommaire}\\

		\begin{enumerate}
		\item[0.] Introduction
		\item[1.] Rel\`evements de Teichm\" uller
		\item[2.] Fonctions $L$ des $F$-modules convergents
		\item[3.] Fonctions $L$ des $F$-modules surconvergents
		\item[4.] Formule des traces de Monsky g\'en\'eralis\'ee
		\item[5.] Conjecture de Dwork pour les $F$-modules surconvergents

		\end{enumerate}
\newpage

\noindent\textbf{R\'esum\'e}\\

Cet article est le premier d'une s\'erie de trois articles consacr\'es aux fonctions $L$. Dans celui-ci nous d\'efinissons les fonctions $L$ des $F$-modules convergents ou surconvergents \`a l'aide des rel\`evements de Teichm¬\"uller et nous \'etablissons la m\'eromorphie  des fonctions $L$ des $F$-modules convergents dans le disque unit\'e ferm\'e. Wan a \'etabli la conjecture de Dwork dans une s\'erie de trois articles [W 2, W 3, W 4] ; par un th\'eor\`eme d'isog\'enie de Katz, sa preuve se ram\`ene au cas ordinaire: nous prouvons ici, sur deux exemples explicites li\'es aux familles de courbes elliptiques, que la filtration par les pentes d'un $F$-module ordinaire surconvergent ne se remonte pas en une filtration surconvergente. Au passage nous montrons que le sous-$F$-isocristal unit\'e dans la cohomologie de de Rham de la famille de Legendre des courbes elliptiques ordinaires n'est pas surconvergent au sens de Berthelot. Dans le deuxi\`eme article nous donnerons une d\'efinition des fonctions $L$ des $F$-(iso)cristaux par voie cohomologique et nous montrerons comment elle rejoint celle donn\'ee ici: elle redonne celle utilis\'ee en cohomologie cristalline par Katz [K 1] ou [Et 1], ou celle en cohomologie rigide de [E-LS 1], ou celle utilis\'ee par Wan [W 2]; le but est alors de donner une preuve de la conjecture de Katz sur les z\'eros et p\^oles unit\'es $p$-adiques de ces fonctions $L$ en utilisant la cohomologie rigide. Dans le troisi\`eme article nous explicitons ces r\'esultats pour les sch\'emas ab\'eliens ordinaires .\\

\vskip5mm
\noindent\textbf{Abstract}\\

This article is the first one of a series of three articles devoted to $L$-functions. In this one we define the $L$-functions of convergent or overconvergent $F$-modules with the help of Teichm¬\"uller liftings and we establish the meromorphy of the $L$-functions of convergent $F$-modules in the closed unit disk. Wan has established Dwork conjecture in a series of three articles [W 2, W 3, W 4]; owing to an isogeny theorem of Katz, his proof reduces to the ordinary case: here we prove, on two explicit examples related to families of ellipitic curves, that the slope filtration on an ordinary overconvergent $F$-module doesn't lift to an overconvergent filtration. As a by-product we show that the unit-root sub-$F$-isocrystal of the de Rham cohomology of the Legendre family of ordinary elliptic curves is not overconvergent in Berthelot's sense. In the second article we'll give a definition of the $L$-functions of $F$-(iso)crystals by cohomological means and we'll show how it matches with the one given here: it gives back the one used in crystalline cohomology by Katz [K 1] or [Et 1], or the one used in rigid cohomology by [E-LS 1], or the one used by Wan [W 2]; the aim is then to give a proof of Katz conjecture on $p$-adic unit roots and poles of these $L$-functions using rigid cohomology. In the third article we give an explicit form of these results for ordinary abelian schemes.\\


Mots cl\'es: alg\`ebres de Monsky-Washnitzer, fonctions $L$, $F$- modules (sur)convergents, formule des traces, m\'eromorphie $p$-adique.


\newpage \vskip 10mm

\section*{0. Introduction}
\textbf{Notations:} Sauf mention du contraire, on suppose dans cet article que $k$ est un corps fini, $k = \mathbb{F}_{q}, q = p^a, \mathcal{V}$ est un anneau de valuation discr\`ete complet, d'id\'eal maximal $\mathfrak{m}$ et corps r\'esiduel $k = \mathbb{F}_{q}$. On suppose le corps des fractions $K$ de $\mathcal{V}$ de caract\'eristique 0,  on fixe une uniformisante $\pi$ et on note $e$ l'indice de ramification. On rel\`eve la puissance $q$ sur $k$ en un automorphisme $\sigma$ de $\mathcal{V}$ suivant la m\'ethode de [Et 4, I 1.1] en supposant que $\sigma(\pi) = \pi$: on note encore $\sigma$ son extension \`a $K$.\\

Cet article est le premier d'une s\'erie de trois articles consacr\'es \`a l'\'etude des fonctions $L$ qui apparaissent en cohomologie rigide, par exemple comme facteurs dans la fonction z\^eta d'une vari\'et\'e $X$ param\'etr\'ee par une autre vari\'et\'e $S$ au-dessus du corps fini $k$. Les \guillemotleft coefficients\guillemotright\  qui apparaissent alors sont les images directes par $f: X \rightarrow S$ du faisceau structural, et la question est de savoir si ces \textbf{images directes} sont adapt\'ees \`a la cohomologie rigide, i.e. si ce sont des $F$-isocristaux surconvergents, question qui a \'et\'e abord\'ee dans [Et 7] et [Et 9]: lorsque c'est le cas, ces fonctions $L$ sont m\'eromorphes [E-LS 1], et m\^eme rationnelles gr\^ace \`a Kedlaya [Ked 1]. L'un des buts de ces trois articles est de d\'ecomposer ces fonctions $L$, lorsqu'elles sont m\'eromorphes, en produits de facteurs $L_{\alpha}$ de \guillemotleft pentes\guillemotright\   $\alpha$ rationnelles diff\'erentes, li\'es aux \textbf{conjectures de Dwork et de Katz}, et de donner une \textbf{description explicite de ces facteurs $L_{\alpha}$} en termes de fonctions $L$ usuelles.\\

Plus pr\'ecis\'ement, soient  $S$ une $k$-vari\'et\'e (i.e. un $k$-sch\'ema s\'epar\'e de type fini) et $f: X \rightarrow S$ un morphisme de $k$-vari\'et\'es. Si l'on d\'esigne par $\vert X\vert $ l'ensemble des points ferm\'es de $X$, la fonction z\^eta de $X$ est d\'efinie par 

$$
Z(X, t)=  \displaystyle \mathop{\prod}_{x \in  \vert X\vert} det(1-t^{\textrm{deg}\ x})^{-1}, \mbox{o\`u} \ {\textrm{deg}\ x} = [k(x): k]\ ,
$$

\noindent et l'on sait par Grothendieck [G 2] que pour $\ell$ premier distinct de $p$ cette expression est donn\'ee par 

$$
Z(X,t)= \displaystyle \mathop{\prod}_{i} det(1 - t\ F | H^i_{\textrm{\'et}, c}(X_{\overline{k}}, \mathbb{Q}_{\ell}))^{(-1)^{i+1}} \qquad (1)\ ,
$$

\noindent o\`u $X_{\overline{k}}$ est l'image inverse de $X$ sur une cl\^oture alg\'ebrique $\overline{k}$ de $k$, et les $H^i_{\textrm{\'et}, c}$ sont des $\mathbb{Q}_{\ell}$-espaces vectoriels de dimension finie sur lesquels le Frobenius $F$ agit.\\

On peut aussi fibrer la situation au-dessus de $S$:

$$
Z(X, t) = \displaystyle \mathop{\prod}_{s \in  \vert S\vert } Z(X_{s}, t^{{\textrm{deg}\ s}}), \mbox{o\`u}\ X_{s}=f^{-1}(s)\ {\textrm {et deg}\ s} = [k(s): k]\ ;
$$
et l'on sait encore par Grothendieck [G 2] que l'on a les expressions
 $$
\begin{array}{rcl}
Z(X, t)&=&  \displaystyle \mathop{\prod}_{s \in  \vert S\vert } \displaystyle \mathop{\prod}_{i}det (1- t^{\textrm{deg}\ s}\ F_{s}^{\textrm{deg}\ s}\vert (R^{i}f_{\acute{e}t,c \ast}\mathbb{Q}_{\ell})_{\overline{s}})^{(-1)^{i+1}}\\
&=&\displaystyle \mathop{\prod}_{i}\ L(S, \mathcal{F}^{i},t)^{(-1)^{i}}\  \mbox{o\`u}\ \mathcal{F}^{i}=R^{i}f_{\acute{e}t,c \ast}\mathbb{Q}_{\ell}\\
&=&\displaystyle\mathop{\prod}_{i,j}det(1 - t\ F | H^j_{\textrm{\'et}, c}(S_{\overline{k}}, \mathcal{F}^{i}))^{(-1)^{i+j+1}} \qquad (2) \ ,
\end{array}
$$
\noindent o\`u $\overline{s}$ est un point g\'eom\'etrique au-dessus de $s$, $S_{\overline{k}}$ est l'image inverse de $S$ sur une cl\^oture alg\'ebrique $\overline{k}$ de $k$, et les $H^i_{\textrm{\'et}, c}$ sont des $\mathbb{Q}_{\ell}$-espaces vectoriels de dimension finie sur lesquels le Frobenius $F$ agit.\\

On obtient ainsi par voie $\ell$-adique la rationalit\'e de la fonction z\^eta (formule (1)) de m\^eme que celle de chacune des fonctions $L$ qui apparaissent dans la fibration de $X$ au-dessus de $S$ (formule (2)).\\

La premi\`ere preuve $p$-adique de rationalit\'e de la fonction z\^eta avait \'et\'e obtenue ant\'erieurement par Dwork [Dw 1], mais par une voie \guillemotleft pr\'ecohomologique\guillemotright. Pour avoir une preuve cohomologique de l'analogue $p$-adique de la formule (1) il faudra attendre les travaux de Monsky et Washnitzer [M-W] pour une vari\'et\'e $X$ affine et lisse sur $k$ (la cohomologie de Monsky-Washnitzer remplace la cohomologie \'etale $\ell$-adique dans la formule (1)) ou la cohomologie cristalline de Berthelot  pour une vari\'et\'e $X$ propre et lisse sur $k$ (la cohomologie cristalline \`a coefficients dans le faisceau structural remplace la cohomologie $\ell$-adique dans la formule (1)). Pour une vari\'et\'e $X$ quelconque la cohomologie rigide de Berthelot est utilis\'ee pour \'etablir la rationalit\'e de la fonction z\^eta dans le cas g\'en\'eral [E-LS 1].\\

Si l'on veut un \textbf{analogue $p$-adique de la formule (2)}, on peut dans un premier temps envisager le cas o\`u $S$ est propre et lisse sur $k$ et $\mathcal{F}^{i}$ correspond \`a un $F$-cristal localement libre de type fini: c'est la situation \'etudi\'ee dans [Et 1]. Mais si l'on souhaite consid\'erer une vari\'et\'e $S$ quelconque sur $k$ il faut remplacer la cohomologie $\ell$-adique \`a supports compacts de (2) par la cohomologie rigide \`a supports compacts de $S$, avec des coefficients adapt\'es \`a la cohomologie rigide, \`a savoir des $F$-isocristaux surconvergents. On est donc amen\'e \`a introduire l'analogue \guillemotleft rigide\guillemotright \ des $\mathcal{F}^{i}$, autrement dit les \textbf{images directes du faisceau structural}
$$
E^{i}:= R^{i}f_{rig \ast}(\mathcal{O}_{X/K}) \ ,
$$
 \noindent o\`u $K$ d\'esigne le corps des fractions de l'anneau $W(k)$ des vecteurs de Witt de $k$, et voir \textbf{dans quelles conditions on obtient leur surconvergence}: cette surconvergence est vraie lorsque $X/S$ est un sch\'ema ab\'elien [Et 4]; l'objet de [Et 7], [Et 9] est l'obtention d'autres conditions de surconvergence pour ces images directes. Le cas d'un $f$ propre et lisse quelconque reste ouvert malgr\'e des r\'esultats de Shiho [Shi 2], [Shi 3], et constitue une conjecture de Berthelot. Lorsque ces conditions de surconvergence sont satisfaites on obtient alors gr\^ace \`a [E-LS 1] la m\'eromorphie de $L(S,E^{i}, t)$ et m\^eme sa rationalit\'e en prenant en compte la finitude de la cohomologie rigide \`a coefficients prouv\'ee par Kedlaya [Ked], ce qui fournit l'analogue $p$-adique de (2).\\
 
 Une fois connue la rationalit\'e (ou la m\'eromorphie) de la fonction z\^eta ou des fonctions $L$ pr\'ec\'edentes, on se pose la question de conna\^itre la partie constitu\'ee des z\'eros et p\^oles unit\'es $p$-adiques de ces fonctions: c'est l'objet des \textbf{conjectures de Dwork et de Katz.}\\
 
 Reprenons l'exemple de la situation relative $f :\ X \rightarrow S$ au-dessus de $k$. La fonction z\^eta-unit\'e de $X/S$ est d\'efinie par Dwork et Sperber dans [Dw-S] par
 $$
 Z_{u}(X/S, t)=  \displaystyle \mathop{\prod}_{s \in \vert S\vert }Z_{u}(X_{s},t^{\textrm{deg}\ s})\ ,
 $$
 \noindent o\`u $Z_{u}(X_{s}, T)$ est obtenue \`a partir de la fonction z\^eta de $X_{s},\ Z(X_{s},T)$, en ne conservant dans cette fraction rationnelle que les facteurs $(1-a T)^{+/-}$ avec $\vert a \vert_{p}=1$.\\
 
 Initialement \textbf{la conjecture de Dwork} consistait \`a se demander si $ Z_{u}(X/S, t)$ \'etait $p$-adiquement m\'eromorphe [Dw 3] [Dw 5]. Or il a \'et\'e \'etabli dans [E-LS 2] que 
 $$
  Z_{u}(X/S, t) = L(S, \mathbb{R}f_{!}\mathbb{Q}_{p}):= \displaystyle \mathop{\prod}_{i}L(S, R^{i}f_{\acute{e}t,c \ast}\mathbb{Q}_{p})^{(-1)^{i}}
 $$
 \noindent o\`u les $\mathcal{G}^{i}:=R^{i}f_{\acute{e}t,c \ast}\mathbb{Q}_{p}$ sont des faisceaux $p$-adiques constructibles. Les $\mathcal{G}^{i}$ peuvent \^etre associ\'es \`a des $F$-isocristaux convergents, mais pas surconvergents en g\'en\'eral comme le prouve un exemple de Crew [C]: par cons\'equent on ne peut utiliser la cohomologie rigide et esp\'erer obtenir la rationalit\'e de ces fonctions $L$ via [E-LS 1] et donc de $Z_{u}(X/S, t)$, sauf si $S$ \'etait propre sur $k$ [E-LS 1] ou si $S$ \'etait propre et lisse sur $k$ [Et 1]. Dans le cas g\'en\'eral d'une vari\'et\'e $S$ quelconque Wan a m\^eme d\'emontr\'e [W 1] que de telles fonctions $L$ (associ\'ees \`a des $F$-isocristaux convergents) ne sont pas n\'ecessairement $p$-adiquement m\'eromorphes. Interm\'ediaire en un certain sens entre $F$-isocristaux convergents et surconvergents se trouve la notion de $F$-modules surconvergents [cf \S 3]: on oublie la connexion et sa surconvergence pour ne garder que la surconvergence du Frobenius; pour de tels $F$-modules surconvergents $\mathcal{H}$ sur $S$ affine et lisse sur $k$, Wan a alors d\'efini les fonctions $L(S, \mathcal{H},t)$ et prouv\'e leur m\'eromorphie (cf [W 1], [W 2], [Dw-S]). De plus Wan d\'efinit dans ce nouveau contexte [W 2] la partie $L_{\alpha}(S, \mathcal{H},t)$ de pente $\alpha \in \mathbb{Q}$ de $L(S, \mathcal{H},t)$ et la nouvelle formulation de la conjecture de Dwork consiste \`a dire que $L_{\alpha}(S, \mathcal{H},t)$ est $p$-adiquement m\'eromorphe: gr\^ace \`a la surconvergence du Frobenius, Wan prouve alors la m\'eromorphie de $L_{\alpha}$ ([W 3], [W 4]).\\
 
  Ce \textbf{premier article} est consacr\'e au cas d'un $S$ affine et lisse. Nous commen\c cons au \S1 par une \'etude d\'etaill\'ee des rel\`evements de Teichm\" uller: celle-ci nous servira d'une part \`a relier les fonctions $L$ d'un $F$-module sur $S$ au cas de l'espace affine et d'autre part \`a \'etablir ult\'erieurement le lien avec les fonctions $L$ de $F$-isocristaux. Apr\`es avoir pos\'e au \S2 la d\'efinition de la fonction $L(S, \mathcal{H},t)$ d'un $F$-module sur $S$,  nous la relions au \S3 au cas de l'espace affine. Au \S4, apr\`es un rappel sur la formule des traces de Monsky, nous prouvons que la fonction $L$ d'un $F$-module convergent est m\'eromorphe dans le disque unit\'e ferm\'e. Au \S5 nous montrons comment le th\'eor\`eme d'isog\'enie de Katz ram\`ene la preuve de Wan au cas ordinaire: c'est dans ce contexte que nous prouvons, sur deux exemples explicites li\'es aux familles de courbes elliptiques ordinaires, que la filtration par les pentes d'un $F$-module ordinaire surconvergent ne se remonte pas en une filtration surconvergente. Au passage nous montrons que le sous-$F$-isocristal unit\'e dans la cohomologie de de Rham de la famille de Legendre des courbes elliptiques ordinaires n'est pas surconvergent au sens de Berthelot.\\
 
Le \textbf{deuxi\`eme article} [Et 10] est consacr\'e aux fonctions $L$ des $F$-(iso)cristaux: nous globalisons les d\'efinitions et r\'esultats pr\'ec\'edents au cas de vari\'et\'es non n\'ecessairement affines avec pour coefficients des $F$-isocristaux Dwork-surconvergents et nous abordons la \textbf{conjecture de Dwork} dans ce contexte, ainsi que la \textbf{conjecture de Katz} relative aux z\'eros et p\^oles unit\'es $p$-adiques de la fonction $L$ d'un $F$-cristal unit\'e.\\

 En faisant une hypoth\`ese suppl\'ementaire sur $f$ dans le \textbf{troisi\`eme article} [Et 11] , \`a savoir si $X/S$ constitue un sch\'ema ab\'elien ordinaire, \textbf{nous explicitons dans la fraction rationnelle $L(S, R^{i}f_{rig \ast}(\mathcal{O}_{X/K}))$} (c'est une fraction rationnelle gr\^ace \`a la surconvergence des images directes $R^{i}f_{rig \ast}(\mathcal{O}_{X/K}$ [Et 4]) \textbf{la partie $L_{\alpha}$ de pente $\alpha \in \mathbb{Q}$ comme une fonction $L$ usuelle} d'un certain $F$-cristal localement libre de type fini associ\'e au cristal de Dieudonn\'e du sch\'ema ab\'elien: ceci n\'ecessite au pr\'ealable une caract\'erisation d\'etaill\'ee des sch\'emas ab\'eliens ordinaires en termes du cristal de Dieudonn\'e du groupe $p$-divisible associ\'e au sch\'ema ab\'elien. Si $S$ est propre sur $k$ il r\'esulte alors de [E-LS 1] (ou de [Et 1] si $S$ est propre et lisse sur $k$) que ces fonctions $L_{\alpha}$ sont en fait rationnelles.

\section*{1. Rel\`evements de Teichm\" uller}

Dans ce \S1 on supposera simplement que $k$ est un corps parfait de caract\'eristique $p>0$. Soit $X = \textrm{Spec}\ A_{0}$ un $k$-sch\'ema lisse. Pour $x \in\  \vert X \vert\  =  \{ \textrm{points ferm\'es de}\  X \}$, soit $i_{x} = Spec\ k(x) \hookrightarrow X$ l'immersion ferm\'ee canonique: $k(x) = A_{0}/\mathfrak{m}_{x}$ est une extension finie \'etale de $k$ de degr\'e deg $x = [k(x) : k]$. Notons $W = W(k)$ (resp. $W(x) = W(k(x))$ l'anneau des vecteurs de Witt \`a coefficients dans $k$ (resp. $k(x)$),

$$
\mathcal{V}(x) = W(x) \otimes_{W} \mathcal{V} \simeq W(x) [\pi]\ ,
$$

\noindent $K_{0} = \textrm{Frac}\ W$, $K_{0}(x) = \textrm{Frac}(W(x)), K(x) = \textrm{Frac}(\mathcal{V}(x))$, $\sigma_{x}$ la puissance $p^a$ sur $k(x)$, $\sigma_{W(x)} = W(\sigma_{x})$ le rel\`evement canonique de $\sigma_{x}$ \`a $W(x)$, $\sigma_{\mathcal{V}(x)} = \sigma_{W(x)} \otimes_{W} \mathcal{V}$ et $\sigma_{K(x)}$ (resp. $\sigma_{K_{0}(x)})$ son extension naturelle \`a $K(x)$ (resp. $K_{0}(x))$ d\'efinie par $\sigma_{K(x)}(u/v) = \sigma_{\mathcal{V}(x)}(u)/\sigma_{\mathcal{V}(x)}(v)$ (resp. $\sigma_{K_{0}(x)}(u/v) = \sigma_{W(x)}(u)/\sigma_{W(x)}v))$. Le morphisme $\sigma_{K(x)}$ co¬\"{\i}ncide, d'apr\`es [Et 4, I.1.1] et [B-M, (1.2.7) (ii)], avec le morphisme $\sigma' : K' \rightarrow K'$ (au-dessus de $\sigma : K \rightarrow K)$ de [Et 4, I.1.1] pour $k' = k(x)$.\\

Soit $A$ une $\mathcal{V}$-alg\`ebre lisse relevant $A_{0}$ et fixons une pr\'esentation $A = \mathcal{V}[t_{1},...,t_{d}]/(f_{1},...,f_{\epsilon})$. On d\'esigne par $\hat{A}$ ((resp. $A^{\dag}$) le s\'epar\'e compl\'et\'e (resp. le compl\'et\'e faible) $\mathfrak{m}$-adique de $A$: on a des isomorphismes

$$\qquad \hat{A}\ \simeq\ \mathcal{V} \{t_{1},...,t_{d} \} / (f_{1},...,f_{\epsilon}) \  ,$$
$$\qquad A^{\dag}\ \simeq\ \mathcal{V}[t_{1},...,t_{d}]^{\dag} / (f_{1},...,f_{\epsilon}).$$

Soient $P$ le compl\'et\'e formel de la fermeture projective de $\mathcal{X} = Spec\ A$ dans $\mathbb{P}^d_{\mathcal{V}}$, $P_{1}
= Spf(\mathcal{V}(x))$, $X_{1} = Spec(k(x))$. D'apr\`es [Et 4, (1.2.1)] il existe un carr\'e commutatif

$$
\xymatrix{
A^{\dag} \otimes_{\mathcal{V}} \mathcal{V}(x) =: A^{\dag}(x)  \ar[rr]^{\qquad \quad F_{A^{\dag}(x)}} & & A^{\dag}(x) \\
A^{\dag} \ar[u] \ar[rr]^{F_{A^{\dag}}}   & & A^{\dag}  \ar[u]
}
$$

\noindent au-dessus du carr\'e commutatif

$$
\xymatrix{
\mathcal{V}(x)  \ar[r]^{\sigma_{\mathcal{V}(x)}}  & \mathcal{V}(x) &\\
\mathcal{V} \ar[u] \ar[r]^{\sigma_{\mathcal{V}}}  & \mathcal{V} \ar[u] &,
}
$$

\noindent o\`u $F_{A^{\dag}}$ est un rel\`evement \`a $A^{\dag}$ du Frobenius (puissance $q$) de $A_{0}$ ; d'o\`u un diagramme commutatif

$$
\xymatrix{
\mathcal{V}(x)  \ar[r] \ar[d]_{\sigma_{\mathcal{V}(x)}} & A^{\dag}(x) \ar[r] \ar[d]_{F_{A^{\dag}(x)}}  & \hat{A}(x) := \hat{A} \otimes_{\mathcal{V}} \mathcal{V}(x) \ar[d]^{F_{\hat{A}(x)}} &\\
\mathcal{V}(x) \ar[r]   & A^{\dag}(x) \ar[r] & \hat{A}(x) & .
}
$$

Par cons\'equent le morphisme $s : A_{0} \rightarrow k(x)$ se rel\`eve de mani\`ere unique d'apr\`es Katz [K 1] en un morphisme

$$
\tau(x) : \hat{A}(x) \rightarrow \mathcal{V}(x)
$$

\noindent tel que le diagramme

$$
\xymatrix{
\hat{A}(x)  \ar[r]^{\tau (x)} \ar[d]_{F_{\hat{A}(x)}} & \mathcal{V}(x) \ar[d]^{\sigma_{\mathcal{V}(x)}}\\
\hat{A}(x)  \ar[r]^{\tau (x)}  & \mathcal{V}(x)
}
$$

\noindent commute : $\tau (x)$ est appel\'e le \textbf{rel\`evement de Teichm\"{u}ller de $s$ (ou de $x$) relativement \`a $F_{\hat{A}(x)}$}.
 Les morphismes compos\'es

$$
\hat{\tau}(x) : \hat{A} \hookrightarrow \hat{A}(x) \displaystyle \mathop{\longrightarrow}^{\tau (x)} \mathcal{V}(x)\ ,
$$

$$
\tau^{\dag}(x) : A^{\dag} \hookrightarrow \hat{A} \displaystyle \mathop{\longrightarrow}^{\hat{\tau} (x)} \mathcal{V}(x)
$$

\noindent sont appel\'es respectivement \textbf{
rel\`evement de Teichm\"{u}ller de $s$ (ou de $x$) 
relativement \`a $F_{\hat{A}}$ ou $F_{A^{\dag}}: \hat{\tau}(x) (\mbox{resp.} \ \tau^{\dag}(x))$ 
 est un point de Teichm\"{u}ller de $\hat{A}$ (resp. de $A^{\dag}$) relativement \`a $F_{\hat{A}}$ (resp. $F_{A^{\dag}}$).
									}\\
									
Pour $n\in \mathbb{N}$ les applications $\hat{\tau}(x)$ et $\tau^{\dag}(x)$ donnent, en quotientant par $\pi^{n+1}$, une application
$$
\tau_{n}(x): A_{n}=\hat{A}/\pi^{n+1}\hat{A}=A^{\dag}/\pi^{n+1}/A^{\dag}\rightarrow \mathcal{V}_{n}(x)=\mathcal{V}(x)/\pi^{n+1}\mathcal{V}(x) \  .
$$

Au lieu de consid\'erer $x$ comme point ferm\'e de $X$ on peut aussi le consid\'erer comme point ferm\'e de $\mathbb{A}_{k}^{d}= Spec\ k[\underline{t}]$ o\`u $\underline{t}= (t_{1},...,t_{d})$: on va voir qu'alors  $\tau^{\dag}(x)$ s'\'etend en un rel\`evement de Teichm\"uller de $x$ \`a $\mathcal{V}[\underline{t}]^{\dag}$. Notons $R=\mathcal{V}[\underline{t}],\ R^{\dag}$ (resp. $\hat{R}$) son compl\'et\'e faible (resp. son s\'epar\'e compl\'et\'e) et $I$ le noyau de la surjection canonique
$$
\mu: R^{\dag}\twoheadrightarrow A^{\dag}\ .
$$
\noindent On peut relever (de mani\`ere non unique) le Frobenius $F_{A^{\dag}}$ de $A^{\dag}$ en un endomorphisme $F_{R^{\dag}}$ de $R^{\dag}$ de la fa\c con suivante: il suffit de choisir des \'el\'ements $F_{R^{\dag}}(t_{i})$ de $R^{\dag}$ tels que

$$
F_{R^{\dag}}(t_{i})\equiv t_{i}^{q} (\mbox{mod}\ \pi),\ F_{R^{\dag}}(t_{i})\in \mu^{-1}(F_{A^{\dag}}(\mu(t_{i})))\ ,
$$
\noindent choix qu'il est possible de faire car $F_{A^{\dag}}(\mu(t_{i}))\equiv t^{q}_{i}$ ( mod ($\pi, I$)). On \'etend ce choix d'\'el\'ements $F_{A^{\dag}}(\mu(t_{i}))$ en un endomorphisme $F_{R^{\dag}}$ de la $\mathcal{V}$-alg\`ebre $R^{\dag}$ tel que le diagramme
$$
\xymatrix{
R^{\dag} \ar@{->>}[r]^{\mu} \ar[d]_{F_{R^{\dag}}} & A^{\dag} \ar[d]^{F_{A^{\dag}}}\\
R^{\dag} \ar@{->>}[r]^{\mu}  &A^{\dag}
}
$$
\noindent commute; en particulier on a $F_{R^{\dag}}(I)\subset I$.\\
Les morphismes $F_{R^{\dag}}$ et $F_{A^{\dag}}$ sont finis et fid\`elements plats puisque leur r\'eduction mod $ \pi$ le sont [Et 3, th\'eo 17]; par changement de base de $R^{\dag}$ \`a $\hat{R}$ appliqu\'e au diagramme pr\'ec\'edent on en d\'eduit le diagramme commutatif

$$
\xymatrix{
\hat{R} \ar@{->>}[r]^{\hat{\mu}} \ar[d]_{F_{\hat{R}}} &\hat{A} \ar[d]^{F_{\hat{A}}}\\
\hat{R} \ar@{->>}[r]^{\hat{\mu}}  &\hat{A}}
$$

\noindent Par cons\'equent le morphisme compos\'e
$$
k[\underline{t}]\twoheadrightarrow A_{0}  \displaystyle \mathop{\longrightarrow}^{s} k(x)
$$  
\noindent correspondant au point $x$ de $\mathbb{A}^{d}_{k}$ se rel\`eve de mani\`ere unique [K 1] en un morphisme

$$
\hat{\tau}_{R}(x) : \hat{R} \hookrightarrow \hat{A} \displaystyle \mathop{\longrightarrow}^{\hat{\tau} (x)} \mathcal{V}(x)
$$
\noindent tel que le diagramme
$$
\xymatrix{
\hat{R}  \ar[r]^{\hat{\tau}_{R} (x)} \ar[d]_{F_{\hat{R}}} & \mathcal{V}(x) \ar[d]^{\sigma_{\mathcal{V}(x)}}\\
\hat{R}  \ar[r]^{\hat{\tau}_{R} (x)}  & \mathcal{V}(x)
}
$$

\noindent commute : ${\hat{\tau}_{R} (x)}$ est appel\'e le \textbf{rel\`evement de Teichm\"{u}ller de  $x$ relativement \`a $F_{\hat{R}}$}. Le morphisme compos\'e
$$
{\tau}^{\dag}_{R}(x) : R^{\dag} \hookrightarrow \hat{R} \displaystyle \mathop{\longrightarrow}^{\hat{\tau}_{R} (x)} \mathcal{V}(x)
$$
\noindent est appel\'e le \textbf{rel\`evement de Teichm\"{u}ller de  $x$ relativement \`a $F_{R^{\dag}}$}.
Par l'unicit\'e des rel\`evements de Teichm\"uller prouv\'ee par Katz [K 1], il y a ainsi \textbf{bijection entre les points de Teichm\"uller de $A^{\dag}$ relativement \`a $F_{A^{\dag}}$ et les points de Teichm\"uller de $R^{\dag}$ relativement \`a $F_{R^{\dag}}$ qui se r\'eduisent mod $ \pi$ en des points de $X$}.\\

De m\^eme, pour toute extension finie $k \hookrightarrow k'$, il y a une bijection entre les trois ensembles suivants:
\begin{enumerate}
\item[(i)] l'ensemble des points $x\in X(k')$= \{points de X \`a valeur dans k'\},
\item[(ii)] l'ensemble des points de Teichm\"uller de $A^{\dag}$ relativement \`a $F_{A^{\dag}}$ \`a valeur dans $\mathcal{V}(k'):=\ W(k')\otimes_{W} \mathcal{V}$,
\item[(iii)] l'ensemble des points de Teichm\"uller de $R^{\dag}$ relativement \`a $F_{R^{\dag}}$ \`a valeur dans $\mathcal{V}(k'):=\ W(k')\otimes_{W} \mathcal{V}$ qui se r\'eduisent mod $ \pi$ en des points de $X$ \`a valeur dans $k'$.
\end{enumerate}

 Les morphismes $\hat{\tau}(x)$  et $\tau^{\dag}(x)$ (resp. $\hat{\tau}_{R}(x)$  et $\tau^{\dag}_{R}(x)$) sont surjectifs car la r\'eduction de $\tau^{\dag}(x) \ \textrm{mod}\  \pi$ (resp. $\tau^{\dag}_{R}(x) \ \textrm{mod}\  \pi$) est le morphisme surjectif $s : A_{0} \twoheadrightarrow k(x)$ de d\'epart (resp. le morphisme surjectif $k[\underline{t}]\twoheadrightarrow k(x)$) [M-W, theo 3.2]. Donc $\mathcal{V}(x)$ est un quotient de $\hat{A}$ et $\mathcal{V}(x) \simeq W(x) [\pi]$, qui est un anneau de valuation discr\`ete, est une extension finie \'etale de $\mathcal{V}$ de rang deg $x$. Le noyau du morphisme

$$
\xymatrix{
\hat{\tau}_{K}(x) := \hat{\tau}(x) \otimes_{\mathcal{V}} K : \hat{A}_{K} \ar@{->>}[r] & K(x) = \textrm{Frac} (\mathcal{V}(x))
}
$$

\noindent est ainsi un id\'eal maximal $\mathfrak{q}_{x}$ de $\hat{A}_{K}$ et le diagramme

$$
\begin{array}{c}
\xymatrix{
A^{\dag}_{K} \ar@{^{(}->}[r] \ar[d]_{F_{A^{\dag}_{K}}} & \hat{A}_{K} \ar[d]_{F_{\hat{A}_{K}}} \ar@{->>}[r]^{\hat{\tau}_{K }(x)} & K(x) \ar[d]^{\sigma_{K (x)}}\\
A^{\dag}_{K} \ar@{^{(}->}[r] & \hat{A}_{K} \ar@{->>}[r]^{\hat{\tau}_{K} (x)} & K(x)
}
\end{array}
\leqno{(1.1)}
$$

\noindent commute. On notera $\tau^{\dag}_{K}(x)$ la fl\`eche compos\'ee
$$
\tau^{\dag}_{K}(x) = \tau^{\dag}(x) \otimes_{\mathcal{V}} K :  A^{\dag}_{K} \displaystyle \mathop{\longrightarrow}^{\varphi}  \hat{A}_{K} \displaystyle -\hspace{-10pt}-\hspace{-10pt}-\hspace{-10pt}-\hspace{-10pt}-\hspace{-10pt} \mathop{\twoheadrightarrow}^{\hspace{-10pt}\hat{\tau}_{K }(x)}  K(x)\ ; \leqno{(1.2)}
$$\\
remarquons que $\tau^{\dag}_{K}(x)$ est aussi surjectif  car $\varphi$ induit une bijection entre les id\'eaux maximaux de $\hat{A}_{K}$ et ceux de $A^{\dag}_{K}$ [G-K 2, theo 1.7]. Par le morphisme de sp\'ecialisation [B 2, (0.2.2.1)]

$$
sp : Spm\  \hat{A}_{K} \rightarrow Spec\  A_{0}
$$

\noindent l'image de $\{ \mathfrak{q}_{x} \}$ n'est autre que $\{ \mathfrak{m}_{x} \}$ ; de plus $\hat{\tau}(x) : \hat{A} \twoheadrightarrow \mathcal{V}(x)$ est localis\'e en $\{ \mathfrak{p}_{x} \} \in Spec\ \hat{A}$, o\`u
$\mathfrak{p}_{x}$ est l'unique id\'eal maximal de $\hat{A}$ au-dessus de $\mathfrak{m}_{x}$. En d\'efinissant l'application (encore appel\'ee \textbf{rel\`evement de Teichm¬\"uller})\\

\noindent (1.3) $\qquad \qquad \qquad \qquad \hat{T}_{K} : Spm\ A_{0} \rightarrow Spm\ \hat{A}_{K}$\\

\noindent par $\hat{T}_{K}(x) = \{ \textrm{Ker}\ \hat{\tau}_{K}(x) \} = \{ \mathfrak{q}_{x} \}$, on vient de prouver que $\hat{T}_{K}$ est \textbf{une section du  morphisme de sp\'ecialisation}, consid\'er\'e comme une application

$$
sp : \vert Spm\ \hat{A}_{K} \vert\ \rightarrow \vert Spec\ A_{0} \vert .
$$

La proposition suivante nous sera utile au \S5:
\vskip 3mm
\noindent \textbf{Proposition (1.4)}. \textit{Soit $\phi: \mathcal{M}\rightarrow\mathcal{N}$ un morphisme de $\hat{A}$-modules tels que $\mathcal{N}$ soit de type fini. Les assertions suivantes sont \'equivalentes:}
\begin{itemize}
\item[(i)] $\phi$ \textit{est surjectif.}
\item[(ii)] \textit{Pour tout point ferm\'e $x$ de $Spec\ A_{0}$ le morphisme}
$$
\phi_{x}:=\hat{\tau}(x)^{\ast}(\phi)\ :\ \mathcal{M}_{x}=\hat{\tau}(x)^{\ast}(\mathcal{M})\rightarrow\mathcal{N}_{x}=\hat{\tau}(x)^{\ast}(\mathcal{N})
$$
\noindent \textit{est surjectif.}
\end{itemize}
\vskip 3mm
\noindent \textit{D\'emonstration}. L'implication $(i)\Rightarrow(ii)$ est claire par exactitude \`a droite du produit tensoriel.\\

Pour l'implication r\'eciproque, supposons $(ii)$. Soit $x$ un point ferm\'e de $Spec\ A_{0}$ correspondant \`a l'id\'eal maximal $ \mathfrak{m}_{x} $ de $A_{0}$ et $\mathfrak{p}_{x}$ l'unique id\'eal maximal de $\hat{A}$ au-dessus de $\mathfrak{m}_{x}$. Le rel\`evement de Teichm\" uller
$\hat{\tau}(x) : \hat{A} \twoheadrightarrow \mathcal{V}(x)$
se factorise via
$\hat{\tau}'(x) : \hat{A}_{\mathfrak{p}_{x}} \twoheadrightarrow \mathcal{V}(x)$
et $\hat{\tau}'(x)$ induit un isomorphisme
$$
 \hat{A}_{\mathfrak{p}_{x}}/ \mathfrak{p}_{x}\hat{A}_{\mathfrak{p}_{x}}\simeq \mathcal{V}(x)/\pi\ \mathcal{V}(x)= k(x)\ .
 \leqno{(1.4.1)}
$$
\noindent Notons $\phi_{\mathfrak{p}_{x}}=\phi\otimes_{\hat{A}}\hat{A}_{\mathfrak{p}_{x}}$. La surjectivit\'e de $\phi_{x}$ fournit celle de $\phi_{x}\ \mbox{mod}\ \pi$; or l'isomorphisme (1.4.1) identifie $\phi_{\mathfrak{p}_{x}}\ \mbox{mod}\ \mathfrak{p}_{x}$ \`a $\phi_{x}\ \mbox{mod}\ \pi$. D'o\`u la surjectivit\'e 
de $\phi_{\mathfrak{p}_{x}}$ par Nakayama $[EGA\ 0_{I}; (7.1.14)]$
et ceci pour tout id\'eal maximal $\mathfrak{p}_{x}$ de $\hat{A}$. La proposition en r\'esulte. $\square$\\

\section*{2. Fonctions $L$ des $F$-modules convergents}

Avec les notations du \S1 et $\mathcal{A}=\hat{A}$ ou $\hat{A}_{K}$ ou $A_{n}$ on d\'esigne par $\mathbf{F^a\mbox{-}\textrm{\bf Mod}(\mathcal{A})}$ (resp. $\mathbf{F^a\mbox{-}\textrm{\bf Modlib}(\mathcal{A})}$) la cat\'egorie des $\mathcal{A}$-modules projectifs de type fini (resp. des $\mathcal{A}$ -modules libres de type fini) $\mathcal{M}$ munis d'un morphisme de Frobenius (non n\'ecessairement un isomorphisme) 

$$\phi_{\mathcal{M}} : \mathcal{M}^{\sigma} :=  F^{\ast}_{\mathcal{A}}(\mathcal{M}) \rightarrow \mathcal{M}
$$

\noindent  Un tel $\mathcal{M}$ est appel\'e un $F$-\textbf{module convergent} (resp. un $F$-\textbf{module  libre convergent}). On note  $\mathbf{F^a\mbox{-}\textrm{\bf Mod}(\hat{A})^{0}}$ (resp. $\mathbf{F^a\mbox{-}\textrm{\bf Modlib}(\hat{A})^{0}}$) la sous-cat\'egorie de $\mathbf{F^a\mbox{-}\textrm{\bf Mod}(\hat{A})}$ (resp. $\mathbf{F^a\mbox{-}\textrm{\bf Modlib}(\hat{A})}$) form\'ee des objets $\mathcal{M}$ unit\'es, i.e. tel que le Frobenius soit un isomorphisme: un tel $\mathcal{M}$ est appel\'e $F$-\textbf{module convergent unit\'e} (resp. $F$-\textbf{module libre convergent unit\'e}).\\

Soit $\mathcal{M} \in F^a\mbox{-}\textrm{Mod}(\hat{A}_{K})$.
La fibre $\mathcal{M}_{x}$ de $\mathcal{M}$ en $x \in \vert X \vert$ est par d\'efinition\\

\noindent (2.1) $\qquad \qquad \qquad \mathcal{M}_{x} := \hat{\tau}_{K}(x)^{\ast} (\mathcal{M})$ , \\

\noindent et $\phi_{\mathcal{M}}$ induit\\

\noindent (2.2) $\qquad \qquad \phi_{x} = \phi_{\mathcal{M}} \otimes_{\hat{A}_{K}} K(x) : \sigma^{\ast}_{K(x)}(\mathcal{M}_{x}) \rightarrow \mathcal{M}_{x} $\ , \\

\noindent d'apr\`es la commutativit\'e du diagramme (1.1).\\

L'it\'er\'e deg $x$ fois de $\phi_{x}$ est un endomorphisme $K(x)$-lin\'eaire du $K(x)$-espace vectoriel de dimension finie $\mathcal{M}_{x}$

$$
\phi^{\textrm{deg}\  x}_{x} : \mathcal{M}_{x} \rightarrow \mathcal{M}_{x}\ .
$$

\noindent Notons\\

\noindent (2.3) $\qquad \qquad \textrm{det} (\mathcal{M}_{x},T) = \textrm{det} (1-T\ \phi^{\textrm{deg}\  x}_{x}, \mathcal{M}_{x})$\\

\noindent le \guillemotleft polyn\^ome caract\'eristique\guillemotright \ de $\phi^{\textrm{deg}\  x}_{x}$.\\

Pour $\mathcal{M} \in F^a\mbox{-}\textrm{Mod} (\hat{A})$ (resp. $\mathcal{M} \in F^a\mbox{-}\textrm{Mod} (A_{n})$) 

\noindent on d\'efinit de m\^eme\\
$$
\mathcal{M}_{x} := \hat{\tau}(x)^{\ast} (\mathcal{M})\  (\mbox{resp.}\  \mathcal{M}_{x} := \hat{\tau}_{n}(x)^{\ast} (\mathcal{M}))
$$

\noindent (2.4)$ \qquad \qquad \textrm{det}(\mathcal{M}_{x}, T) = \textrm{det}(1-T\ \phi^{\textrm{deg}\  x}_{x}, \mathcal{M}_{x}).$\\

Soit $(\mathcal{M}, \phi)\in F^a\mbox{-}\textrm{Mod} (\hat{A})$; posons $\mathcal{M}_{n}= \mathcal{M}/\pi^{n+1}\mathcal{M}$ et $\phi_{n}=\phi\  \mbox{mod}\  \pi^{n+1}$; et pour $x \in \vert X \vert$ notons 
$$
(\mathcal{M}_{n})_{x}= \hat{\tau}_{n}(x)^{\ast} (\mathcal{M}_{n})), (\phi_{n})_{x}=\hat{\tau}_{n}(x)^{\ast}(\phi_{n}),
$$
$$
(\mathcal{M}_{x})_{n} =\mathcal{M}_{x}/\pi^{n+1}\mathcal{M}_{x} , (\phi_{x})_{n}=\phi_{x}\  \mbox{mod}\ \pi^{n+1} \ .
$$
Puisque ($\mathcal{M}_{n})_{x}=(\mathcal{M}_{x})_{n}=:\mathcal{M}_{x,n}$ et $(\phi_{n})_{x}= (\phi_{x})_{n}=: \phi_{x,n}$
on a la relation\\

\noindent (2.4 bis) $ \qquad \textrm{det}(1-T\ (\phi_{x,n})^{\textrm{deg}\  x}, \mathcal{M}_{x,n})\equiv\textrm{det}(1-T\ \phi^{\textrm{deg}\  x}_{x}, \mathcal{M}_{x})\  \mbox{mod}\ \pi^{n+1}\ .$ \\

\noindent \textbf{Lemme (2.5)}. \textit{Avec les notations pr\'ec\'edentes, on a :}

\begin{itemize}
\item[(i)] \textit{Si  $\mathcal{M} \in F^a\mbox{-}Mod(\hat{A}_K)$ , alors det $(\mathcal{M}_{x},T) \in K[T]  $ .}
\item[(ii)] \textit{Si $\mathcal{M} \in F^a\mbox{-}Mod(\hat{A})$, alors det  $(\mathcal{M}_{x},T)  \in \mathcal{V}[T] $.}
\item[(iii)] \textit{Si $\mathcal{M} \in F^a\mbox{-}Mod(A_{n})$, alors det  $(\mathcal{M}_{x},T)  \in \mathcal{V}_{n}[T] $.}
\end{itemize}

\vskip 3mm
\noindent \textit{D\'emonstration}. Pour (i), soient $(e_{i})_{i=1,...,r}$ (resp. $(e_{i} \otimes 1)_{i=1,...,r})$ une base locale de $\mathcal{M}$ (resp. de $\mathcal{M}^{\sigma})$, et $C(\underline{X})$ la matrice de $\phi_{\mathcal{M}}$ dans ces bases respectives. Alors 
$$
C(\underline{X}) = \displaystyle \mathop{\Sigma}_{\underline{u} \in \mathbb{N}^d } a_{\underline{u}}\ \underline{X}^{\underline{u}},
$$

\noindent avec $a_{\underline{u}}\ \in \pi^{\alpha}\ M_{r}(\mathcal{V})$ pour un $\alpha \in \mathbb{Z}$ et 

$$
det(\mathcal{M}_{x}, T) = det \{ 1-T\ C(\underline{t}(x)^{q^{(\textrm{deg}\  x)-1}}) \times ... \times C(\underline{t}(x)^{q}) \times C(\underline{t}(x)) \}
$$

\noindent o\`u $\underline{t}(x)^{\beta} = t_{1}(x)^{\beta} \times ... \times t_{d}(x)^{\beta}$ et les $t_{j}(x)$ sont les coordonn\'ees de $\hat{\tau}_{K}(x)$. Il est clair que det $(\mathcal{M}_{x},T)$ a des coefficients invariants par l'action du Frobenius $\sigma_{K(x)}$, car $\sigma_{K(x)}$ envoie $\underline{t}(x)$ sur $\underline{t}(x)^q$ ; d'o\`u le (i).\\

Les cas (ii) et (iii) sont analogues. $\square$\\

Si l'on note $\tilde{\mathcal{M}}_{x}$ l'espace vectoriel $\mathcal{M}_{x}$ vu comme $K$-espace vectoriel et $\tilde{\phi}_{x}$ son endomorphisme de Frobenius on a

$$
det(\mathcal{M}_{x}, T) = det(1-T\ \phi^{\textrm{deg}\ x}_{x}) = det_{K}(1-T\ \tilde{\phi}^{\textrm{deg}\ x}_{x})^{-1/\textrm{deg}\ x}\ .
$$

\vskip 3mm
\noindent \textbf{D\'efinition-proposition (2.6)}. \textit{La fonction $L$ de $(\mathcal{M}, \phi_{\mathcal{M}}) \in F^a\mbox{-}Mod(\hat{A}_{K})$ est d\'efinie par}

$$
L(Spec\ A_{0}, \mathcal{M}, t) = \displaystyle \mathop{\prod}_{x \in \vert Spec\ A_{0} \vert} det(1-t^{\textrm{deg}\ x}\  \phi^{\textrm{deg}\ x}_{x} \mid \mathcal{M}_{x})^{-1} \in K[[t]]
$$
$\qquad \qquad \qquad \qquad \qquad = \displaystyle \mathop{\prod}_{x \in \vert Spec\ A_{0} \vert} 
det(1-t^{\textrm{deg}\ x}\  \tilde{\phi}^{\textrm{deg}\ x}_{x} \mid \tilde{\mathcal{M}}_{x})^{-1/\textrm{deg}\ x}\ .$\\

\noindent \textit{Si $(\mathcal{M}, \phi_{\mathcal{M}}) \in F^a\mbox{-}\textrm{Mod}(\hat{A})$ on d\'efinit de m\^eme $L(Spec\ A_{0}, \mathcal{M}, t)$ et alors}

$$
L(Spec\ A_{0}, \mathcal{M}, t) = L(Spec\ A_{0}, \mathcal{M}_{K},t) \in \mathcal{V}[[t])\ ,
$$

\noindent \textit{o\`u l'on a pos\'e}
$$
(\mathcal{M}_{K}, \phi_{\mathcal{M}_{K}}) := (\mathcal{M}, \phi_{\mathcal{M}}) \otimes_{\hat{A}} \hat{A}_{K}\ .
$$
 \textit{La fonction $L$ de $(\mathcal{M}, \phi_{\mathcal{M}}) \in F^a\mbox{-}Mod(A_{n})$ est d\'efinie par}

$$
L(Spec\ A_{0}, \mathcal{M}, t) = \displaystyle \mathop{\prod}_{x \in \vert Spec\ A_{0} \vert} det(1-t^{\textrm{deg}\ x}\  \phi^{\textrm{deg}\ x}_{x} \mid \mathcal{M}_{x})^{-1} \in \mathcal{V}_{n}[[t]] \ .
$$
\noindent \textit{Si $(\mathcal{M}, \phi_{\mathcal{M}}) \in F^a\mbox{-}\textrm{Mod}(\hat{A})$, il r\'esulte de (2.4 bis) que}
$$
L(Spec\ A_{0}, \mathcal{M}, t)\equiv L(Spec\ A_{0}, \mathcal{M}_{n}, t) \ \mbox{mod}\ \pi^{n+1} \ .
$$

\vskip 3mm
\noindent \textbf{Lemme (2.7)}. \textit{Pour \'etablir la m\'eromorphie $p$-adique de $L(Spec\ A_{0}, \mathcal{M}, t)$ pour $(\mathcal{M}, \phi_{\mathcal{M}}) \in F^a\mbox{-}Mod(\hat{A} _{K})$ on peut supposer qu'il existe $(\mathcal{M'}, \phi_{\mathcal{M'}}) \in F^a\mbox{-}Mod(\hat{A})$ tel que $\mathcal{M'}$ est libre et $(\mathcal{M}, \phi_{\mathcal{M}}) = (\mathcal{M'}, \pi^{\alpha} \phi_{\mathcal{M'}})  \otimes_{\hat{A}} \hat{A}_{K}$ pour un $\alpha \in \mathbb{Z}$. On a alors}

$$
L(Spec\ A_{0}, \mathcal{M}, t) = L(Spec\ A_{0}, \mathcal{M'}, \pi^{\alpha} t).
$$

\vskip 3mm
\noindent \textit{D\'emonstration}. Puisque $\mathcal{M}$ est projectif de type fini sur $\hat{A}_{K}$, et qu'un ouvert de $Spec\ \hat{A}_{K}$ est intersection d'un ouvert de $Spec\ \hat{A}$ avec $Spec\ \hat{A}_{K}$, il existe un recouvrement fini de $Spec\ \hat{A}_{K}$ par des ouverts $U_{K} = Spec\ B_{K}$, o\`u $B = \hat{A}\ [1/g], g \in \hat{A}$, tels que

$$
\mathcal{M} \otimes_{\hat{A}_{K}} B_{K}  \simeq \displaystyle \mathop{\oplus}_{i=1}^r B_{K}\ e_{i}.
$$

\noindent Relevons $g\  \textrm{mod}\  \pi =: g_{0}$ en $f \in A$ ; comme on a un isomorphisme $\hat{B} \simeq \widehat{A[1/f]}$ [Et 3, cor 1 du th\'eo 4] on en d\'eduit que

$$
\mathcal{N}  := \mathcal{M} \otimes_{\hat{A}_{K}} \hat{B}_{K} \simeq \displaystyle \mathop{\oplus}_{i=1}^r \widehat{A[1/f]_{K}}\  e_{i}
 $$

\noindent et que la matrice $C(\underline{X})$ du Frobenius $\phi_{\mathcal{N}}$ est \`a coefficients dans $\pi^{\alpha} \widehat{A[1/f]}$, avec $\alpha \in \mathbb{Z}$.\\

La fonction $L(Spec\ A_{0},\mathcal{M}, t)$ est d\'efinie par un produit eul\'erien sur les points ferm\'es de $Spec\ A_{0}$ et ceux-ci sont en bijection avec les points ferm\'es de $Spf\ \hat{A}$ : un recouvrement ouvert fini de $Spf\ \hat{A}$ est fourni par des $Spf\ \widehat{A[1/f]}$, $f \in A$ comme ci-dessus ; on peut donc prendre

$$
\mathcal{M'} = \displaystyle \mathop{\oplus}_{i=1}^r \widehat{A[1/f]}\  e_{i}
$$

\noindent avec pour matrice du Frobenius $\phi_{\mathcal{M'}}$ la matrice $\pi^{- \alpha}\  C(\underline{X})$, \`a coefficients dans $\widehat{A[1/f]}$ :

$$
\xymatrix{
\mathcal{M'}^{\sigma} \ar[rr]^{\phi_{\mathcal{M'}}}  \ar@{^{(}->}[d] && \mathcal{M'} \ar@{^{(}->}[d] \\
\mathcal{N}^{\sigma} = \mathcal{M}^{\prime\sigma}_{K} \ar[rr]_{\pi^{- \alpha}{\phi_{\mathcal{N}}}}&& \mathcal{N} = \mathcal{M}'_{K} .
}
$$

Le couple $(\mathcal{M'}, \phi_{\mathcal{M'}})$ est ce que Wan appelle une $\sigma$-module convergent [W 2], [W 3]. $\square$

\newpage
\section*{3.  Fonctions $L$ des $F$-modules surconvergents} 

Avec les notations du \S1 et $\mathcal{A} = A^{\dag}$ ou $A^{\dag}_{K}$ et par analogie avec le \S2 on d\'esigne par $\mathbf{F^a\mbox{-}\textrm{\bf Mod}(\mathcal{A})}$ ( resp. $\mathbf{F^a\mbox{-}\textrm{\bf Modlib}(\mathcal{A})}$) la cat\'egorie des $\textbf{F-\textrm{modules}}$\\
$\textbf{\textrm {surconvergents}}$, i.e. la cat\'egorie des $\mathcal{A}$-modules projectifs  (resp. $\mathcal{A}$-modules libres) de type fini $M$ muni d'un morphisme de Frobenius (non n\'ecessairement un isomorphisme)
$$
\phi_{M} : M^{\sigma} = F^{\ast}_{\mathcal{A}}(M) \rightarrow M
$$

\noindent On note  $\mathbf{F^a\mbox{-}\textrm{\bf Mod}(A^{\dag})^{0}}$ (resp. $\mathbf{F^a\mbox{-}\textrm{\bf Modlib}(A^{\dag})^{0}}$) la sous-cat\'egorie de $\mathbf{F^a\mbox{-}\textrm{\bf Mod}(A^{\dag})}$ (resp. $\mathbf{F^a\mbox{-}\textrm{\bf Modlib}(A^{\dag})}$) form\'ee des objets $M$ unit\'es, i.e. tel que le Frobenius soit un isomorphisme: un tel $M$ est appel\'e $F$-\textbf{module surconvergent unit\'e} (resp. $F$-\textbf{module libre surconvergent unit\'e}).\\

Soit $M \in F^a\mbox{-}\textrm{Mod}(A^{\dag}_{K})$. La fibre $M_{x}$ de $M$ en $x \in \vert X \vert $ est par d\'efinition\\

\noindent (3.1) $\qquad \qquad \qquad \qquad M_{x} := \tau^{\dag}_{K}(x)^{\ast}(M)$,\\

\noindent et $\phi_{M}$ induit\\

\noindent (3.2) $\qquad \qquad \qquad \phi_{x} = \phi_{M} \otimes_{A^{\dag}_{K}} K(x) : \sigma^{\ast}_{K(x)}(M_{x}) \rightarrow M_{x}\ ,$\\

\noindent d'apr\`es la commutativit\'e du diagramme (1.1).\\

L'it\'er\'e deg $x$ fois de $\phi_{x}$ est un endomorphisme $K(x)$-lin\'eaire du $K(x)$-espace vectoriel de dimension finie $M_{x}$

$$
\phi^{\textrm{deg}\ x}_{x} : M_{x} \rightarrow M_{x}\ ;
$$

\noindent on notera $\tilde{M}_{x}$ l'espace vectoriel $M_{x}$ vu comme $K$-espace vectoriel et $\tilde{\phi}_{x}$ son morphisme de Frobenius. \\

Notons\\

\noindent {(3.3)}  $\qquad \qquad \qquad  \textrm{det}(M_{x}, T) =  \textrm{det}(1 - T \phi^{\textrm{deg}\ x}_{x}, M_{x})\ ;$\\

\noindent pour $M \in F^a\mbox{-}\textrm{Mod}(A^{\dag})$ on pose de m\^eme \\

\noindent (3.4)$ \qquad M_{x}=\tau^{\dag}(x)^{\ast}(M)\ ,  \qquad  \textrm{det}(M_{x},T)\   =  \textrm{det}(1 - T\ \phi^{\textrm{deg}\ x}_{x}, M_{x}).$\\

On d\'emontre le lemme suivant comme (2.5).

\vskip 3mm
\noindent \textbf{Lemme (3.5)}. \textit{Avec les notations pr\'ec\'edentes, on a : }

\begin{itemize}
\item[(i)] \textit{Si $M \in F^a\mbox{-}\textrm{Mod} (A^{\dag}_{K})$, alors $det(M_{x}, T) \in K[T]\  $}.
\item[(ii)] \textit{Si $M \in F^a\mbox{-}\textrm{Mod} (A^{\dag})$, alors $det(M_{x}, T) \in \mathcal{V}[T]\ $.}
\end{itemize}

\vskip 3mm
De m\^eme on a :

\vskip 3mm
\noindent \textbf{D\'efinition et proposition (3.6)}. \textit{Soit $(M, \phi_{M}) \in F^a \mbox{-}\textrm{Mod}(A^{\dag}_{K})$ et \\
$(\mathcal{M}, \phi_{\mathcal{M}}) \in F^a\mbox{-}\textrm{Mod}(\hat{A}_{K})$ son image canonique par l'extension des scalaires de $A^{\dag}_{K}$ \`a  $\hat{A}_{K}$. La fonction $L$ de $(M, \phi_{M})$ est d\'efinie par }

$$
L(Spec\ A_{0}, M, t) = \displaystyle \mathop{\prod}_{x \in \vert Spec\ A_{0} \vert}\ det(1-t^{\textrm{deg}\  x}\  \phi^{\textrm{deg}\ x}_{x}, M_{x})^{-1} \in K[[t]]
$$
$\qquad \qquad \qquad \qquad \qquad = \displaystyle \mathop{\prod}_{x \in \vert Spec\ A_{0} \vert}\ det(1-t^{\textrm{deg}\  x}\ \tilde{\phi}^{\textrm{deg}\ x}_{x}, \tilde{M}_{x})^{-1/\textrm{deg}\ x}\ ,
$

\noindent \textit{et on a}

$$L(Spec\ A_{0}, M, t) = L(Spec\ A_{0}, \mathcal{M}, t).$$

\vskip 3mm

\textit{Si $(M, \phi_{M}) \in F^a\mbox{-}\textrm{Mod} (A^{\dag})$ on d\'efinit de m\^eme $L(Spec\ A_{0}, M, t)$ et alors}

$$
L(Spec\ A_{0}, M, t) = L(Spec\ A_{0}, M_{K}, t)=L(Spec\ A_{0}, \mathcal{M}, t) \in \mathcal{V}[[t]] ,
$$
$\qquad\qquad\qquad\qquad\qquad  \equiv L(Spec\ A_{0}, M_{n}, t) \in \mathcal{V}_{n}[[t]] ,$ \\

\noindent \textit{o\`u l'on a pos\'e}

$$
M_{K} := M \otimes_{A^{\dag}} A^{\dag}_{K}\  \ ,\ \mathcal{M} := M  \otimes_{A^{\dag}} \hat{A} \   \  ,\ M_{n}:= M/\pi^{n+1}M \ .
$$

\vskip 3mm
Comme (2.7) on montre :

\vskip 3mm
\noindent \textbf{Lemme (3.7)}. \textit{Pour \'etablir la m\'eromorphie $p$-adique de $L(Spec\ A_{0}, M, t)$ pour $(M, \phi_{M}) \in F^a\mbox{-}\textrm{Mod}(A^{\dag}_{K})$ on peut supposer qu'il existe $(M', \phi_{M'}) \in F^a\mbox{-}\textrm{Mod}(A^{\dag})$ tel que $M'$ est libre et $(M, \phi_{M}) = (M', \pi^{\alpha} \phi_{M'}) \otimes_{A^{\dag}} A^{\dag}_{K}$ pour un $\alpha \in \mathbb{Z}$. On a alors}

$$
L(Spec\ A_{0}, M, t) = L(Spec\ A_{0}, M', \pi^{\alpha} t)\ .
$$

\vskip 3mm

\subsection*{3.8. R\'eduction  au cas de l'espace affine}
Soient $A= \mathcal{V}[x_{1},...,x_{d}]/(f_{1},...,f_{\epsilon})$ une $\mathcal{V}$-alg\`ebre lisse relevant $A_{0}, X= \mbox{Spec}A_{0},  \underline{x}=(x_{1},...,x_{d}), R=\mathcal{V}[\underline{x}], (M,\phi_{M})\in F^{a}\mbox{-Modlib}(A^{\dag})$: un tel $M$ est ce que Wan appelle un $\sigma$-module surconvergent. Dans ce paragraphe nous allons montrer que la fonction $L(X, M, t)$ est une fonction $L$ sur un certain espace affine.\\

Relevons d'abord $M$ de $X$ \`a $\mathbb{A}^{d}_{\mathbb{F}_{q}}$. On choisit une base $(e_{1},...,e_{m}) $ de $M$ sur $A^{\dag}$ et on note $C\in \mathcal{M}_{m}(A^{\dag})$ la matrice $m\times m$ \`a coefficients dans $A^{\dag}$ de l'application $F_{A^{\dag}}$-lin\'eaire $\phi_{M}$. En notant $\mu$ la surjection canonique
$$
\mu:\mathcal{V}[x_{1},...,x_{d}]^{\dag} \twoheadrightarrow A^{\dag} \ ,
$$
\noindent on choisit ensuite $C_{0}$, une matrice $m\times m$ \`a coefficents dans $\mathcal{V}[\underline{x}]^{\dag}$ telle que $\mu(C_{0})=C$. Alors les relations
$$
M_{0}=\mathop{\oplus}_{i=1}^{m}\mathcal{V}[\underline{x}]^{\dag}e_{i} \ \ ,\ \phi_{0}(e_{i})=C_{0}e_{i}
$$
\noindent d\'efinissent un \'el\'ement $(M_{0},\phi_{0})\in F^{a}\mbox{-Modlib}(\mathcal{V}[\underline{x}]^{\dag})$, i.e. un $F$-module libre surconvergent sur $\mathbb{A}^{d}_{\mathbb{F}_{q}}$.\\

Pour se r\'eduire au cas de l'espace affine nous allons faire une petite digression par le $F$-cristal de Dwork [B 1]. Quitte \`a passer \`a une extension totalement ramifi\'ee de $\mathcal{V}$, on peut supposer que $\mathcal{V}$ contient une racine primitive $p$-i\`eme de l'unit\'e: \`a toute racine $\pi^{\ast}$ de l'\'equation
$$
(\pi^{\ast})^{p-1}=-p
$$
\noindent correspond une unique racine primitive $p$-i\`eme de l'unit\'e $\lambda \in K_{0}:=\mathbb{Q}_{p}(\pi^{\ast})$ telle que l'on ait la congruence [Mo 1, theo 4.3]
$$
\lambda \equiv 1+\pi^{\ast} \ (\mbox{mod} \ \pi^{\ast^{2}}) \  \mbox{dans} \  \mathcal{O}_{K_{0}} \ .
$$
\noindent La s\'erie
$$
E(t)= exp(\pi^{\ast}(t-t^{q}))
$$
\noindent qui converge pour $ord_{p}(t)>-\frac{p-1}{qp}$ [Mo 1, theo 4.1] est donc surconvergente et $E(1)$ est une racine primitive $p$-i\`eme de l'unit\'e [Mo 1, theo 4.3]. Soit $\psi_{0}$ le caract\`ere additif non trivial de $\mathbb{F}_{p}$ tel que $\psi_{0}(1)=E(1)$: pour tout $r\in \mathbb{N}^{\ast}$ posons 
$$
\psi_{r}=\psi_{0}\circ Tr_{\mathbb{F}_{q^{r}}/\mathbb{F}_{p}} \ ;
$$
\noindent $\psi_{r}$ est un caract\`ere additif non trivial de $\mathbb{F}_{q^{r}}$. Au caract\`ere $\psi:=\psi_{1}$ est associ\'e le $F$-cristal de Dwork $\mathcal{L}_{\psi}$ sur la droite $\mathbb{A}^{1}_{\mathbb{F}_{q}}$ [B 6]: $\mathcal{L}_{\psi}$ est un $F$-cristal surconvergent (au sens de Berthelot) libre de rang 1 et si l'on note $\theta_{\psi}$ une base de $\mathcal{L}_{\psi}$, le Frobenius
$$
\phi_{\mathcal{L}_{\psi}}: F^{\ast}\mathcal{L}_{\psi}\rightarrow \mathcal{L}_{\psi}
$$
\noindent est donn\'e par [B 1, (1.5.2)]
$$
\phi_{\mathcal{L}_{\psi}}(\theta_{\psi}\otimes 1)=exp(\pi^{\ast}(t-t^{q}))\theta_{\psi}=E(t)\theta_{\psi} \ .
$$
\noindent Posons $\underline{y}=(y_{1}, ..., y_{s})$; pour $i\in \llbracket1,s\rrbracket$, notons $\overline{f_{i}}$ la r\'eduction de $f_{i}$ mod $\pi$ et 
$$
f(\underline{x},\underline{y})=y_{1}\overline{f_{1}}(\underline{x})+ ...+y_{s}\overline{f_{s}}(\underline{x}) \in \mathbb{F}_{q}[\underline{x}, \underline{y}] \ .
$$
\noindent On \'etend le Frobenius $F_{R^{\dag}}$ (d\'efini au \S1) de $R^{\dag}=\mathcal{V}[\underline{x}]^{\dag}$ \`a  $R[\underline{y}]^{\dag}=\mathcal{V}[\underline{x},\underline{y}]^{\dag}$ en $\sigma$, d\'efini par exemple par $y_{i}^{\sigma}=y_{i}^{q}$. Le polyn\^ome $f$ d\'efinit un morphisme encore not\'e $f$\\

$$
\begin{array}{ccc}
 f:\mathbb{A}^{d+s}_{\mathbb{F}_{q}}= Spec\ (k[\underline{x},\underline{y}])& \longrightarrow &Spec \ k[t]  = \mathbb{A}^{1}_{\mathbb{F}_{q}}  \\
 f(\underline{x},\underline{y})&\longleftarrow & t  \ .
 \end{array}
 $$
 
\noindent Par image inverse par $f$, le $F$-cristal de Dwork $\mathcal{L}_{\psi}$ fournit le $F$-cristal $\mathcal{L}_{\psi,f}:=f^{\ast}(\mathcal{L}_{\psi})$ surconvergent sur $\mathbb{A}^{d+s}_{\mathbb{F}_{q}}$, libre de rang 1, de base $e$, et dont le Frobenius  
$$
\phi_{\mathcal{L}_{\psi,f}}: F^{\ast}\mathcal{L}_{\psi,f}\rightarrow \mathcal{L}_{\psi,f}
$$
\noindent est donn\'e par
$$
\phi_{\mathcal{L}_{\psi,f}}(e\otimes 1)=exp(\pi^{\ast}\sum_{i=1}^{s}y_{i}f_{i}(\underline{x})-\pi^{\ast}\sum_{i=1}^{s}y_{i}^{\sigma}f_{i}(\underline{x}^{\sigma}))e=: \phi_{f}(\underline{x},\underline{y})e \ .
$$
\noindent Pour $(\underline{x},\underline{y})\in \mathbb{A}^{d+s}_{\mathbb{F}_{q}}(\mathbb{F}_{q^{r}})$ notons $\tau^{\dag}_{R[\underline{y}]}(\underline{x},\underline{y})=(\tilde{\underline{x}},\tilde{\underline{y}})$ son rel\`evement de Teichm\" uller [\S1] \`a $R[\underline{y}]^{\dag}=\mathcal{V}[\underline{x},\underline{y}]^{\dag}$; alors l'action du Frobenius it\'er\'ee $r=deg\ x$ fois sur $\mathcal{L}_{\psi,f}$ au point $(\underline{x},\underline{y})$ est donn\'ee par (cf la preuve de (2.5)):\\

$$
\begin{array}{rcl}
Tr_{\mathcal{V}}(\phi^{deg\ x}_{{(\mathcal{L}_{\psi,f})_{(\underline{x},\underline{y})})}})&=&\phi_{f}(\tilde{\underline{x}},\tilde{\underline{y}})\times\phi_{f}(\tilde{\underline{x}}^{q},\tilde{\underline{y}}^{q})\times...\times \phi_{f}(\tilde{\underline{x}}^{q^{r-1}},\tilde{\underline{y}}^{q^{r-1}})\\
{}\\
&=&\psi_{r}(f(\underline{x},\underline{y})) \ ,
\end{array} 
\leqno{(3.8.1)}
$$
\noindent la derni\`ere \'egalit\'e r\'esultant de [B 1, (1.4)(ii)]. Puisque $\overline{f_{1}},...,\overline{f_{\epsilon}}$ d\'efinissent $X$ dans $\mathbb{A}^{d}_{\mathbb{F}_{q}}$, un argument standard sur les sommes de caract\`eres permet de remarquer les \'egalit\'es:
$$
\begin{array}{cccc}
\displaystyle\mathop{\sum}_{\underline{y}\in\mathbb{A}^{s}_{\mathbb{F}_{q}}(\mathbb{F}_{q^{r}})}\psi_{r}(f(\underline{x},\underline{y}))&=&q^{rs} &\mbox{si} \ \underline{x}\in X(\mathbb{F}_{q^{r}}) \\
&=&0 &\mbox{si} \ \underline{x}\in (\mathbb{A}^{d}_{\mathbb{F}_{q}}\setminus X)(\mathbb{F}_{q^{r}}) \ ,
\end{array}
\leqno{(3.8.2)}
$$
\noindent o\`u, pour un sch\'ema $Y$ sur $\mathbb{F}_{q}$,  l'on a pos\'e $Y(\mathbb{F}_{q^{r}})=\lbrace \mbox{points de}\ Y \mbox{\`a valeurs dans}\  \mathbb{F}_{q^{r}}\rbrace$.\\

Revenons \`a pr\'esent \`a notre fonction $L$ de d\'epart; on a:
$$
\begin{array}{rcl}
L(X,M,t)&=&\displaystyle\mathop{\prod}_{x\in \vert X\vert}\mbox{det} (1-t^{deg\ x}\phi_{x}^{deg\ x}\vert M_{x})^{-1}\\
&=& \mbox{exp}(\displaystyle\mathop{\sum}_{r=1}^{\infty}\frac{t^{r}}{r}S_{r}(X,M))\\
\end{array}
$$
\noindent avec 
$$
S_{r}(X,M)= \displaystyle \mathop{\sum}_{x\in X(\mathbb{F}_{q^{r}})}Tr_{\mathcal{V}}(\phi_{x}^{r}\vert M_{x}) \ ;
$$
\noindent en effet, pour v\'erifier la concordance des deux \'ecritures de la fonction $L$, il suffit de passer au logarithme dans les deux expressions et regrouper les facteurs du produit eul\'erien par degr\'es (cf [Et 1, II 2]). Dans la d\'efinition de $S_{r}(X,M)$ il est \`a noter que l'on peut remplacer la somme sur les $x\in X(\mathbb{F}_{q^{r}})$ par une somme sur les points de Teichm\" uller $\tau^{\dag}(x)$ de $A^{\dag}$ \`a valeurs dans $\mathcal{V}(\mathbb{F}_{q^{r}}):= W(\mathbb{F}_{q^{r}})\otimes _{W}\mathcal{V}$ relativement \`a $F_{A^{\dag}}$, ceci gr\^ace \`a la bijection entre ces deux ensembles \'etablie au \S1.\\

Pour passer de $X$ \`a l'espace affine $\mathbb{A}^{d+s}_{\mathbb{F}_{q}}$ il va nous suffire de r\'e\'ecrire l'expression de $S_{r}(X,M)$ \`a l'aide de $\mathcal{L}_{\psi,f}$, en utilisant (3.8.1) et (3.8.2) de la fa\c con suivante:
$$
\begin{array}{rcl}
S_{r}(X,M)&=&\frac{1}{q^{rs}}\displaystyle\mathop{\sum}_{(\underline{x},\underline{y})\in (X \times\mathbb{A}^{s})(\mathbb{F}_{q^{r}})}\ Tr_{\mathcal{V}}(\phi_{\underline{x}}^{r})\times \psi_{r}(f(\underline{x},\underline{y}))\\
{}\\
&=&\frac{1}{q^{rs}}\displaystyle\mathop{\sum}_{(\underline{x},\underline{y})\in \mathbb{A}^{d+s}(\mathbb{F}_{q^{r}})}\ Tr_{\mathcal{V}}(\phi_{0,\underline{x}}^{r})\times Tr_{\mathcal{V}}(\phi^{r}_{{(\mathcal{L}_{\psi,f})_{(\underline{x},\underline{y})}}})\\
{}\\
&=&\frac{1}{q^{rs}}\displaystyle\mathop{\sum}_{(\underline{x},\underline{y})\in \mathbb{A}^{d+s}(\mathbb{F}_{q^{r}})}\ Tr_{\mathcal{V}}(\phi^{r}_{N_{(\underline{x},\underline{y})}})\times Tr_{\mathcal{V}}(\phi^{r}_{{(\mathcal{L}_{\psi,f})_{(\underline{x},\underline{y})}}})
\end{array}
$$
\noindent o\`u $N:= pr^{\ast}(M_{0})$ est l'image inverse de $M_{0}$ par la projection canonique $pr:\mathbb{A}^{d+s}_{\mathbb{F}_{q}}\rightarrow \mathbb{A}^{d}_{\mathbb{F}_{q}}$ et $\phi_{N}=pr^{\ast}\ (M_{0})$. D'o\`u
$$
\begin{array}{rcl}
S_{r}(X,M)&=&\frac{1}{q^{rs}}\displaystyle\mathop{\sum}_{(\underline{x},\underline{y})\in \mathbb{A}^{d+s}(\mathbb{F}_{q^{r}})}\ Tr_{\mathcal{V}}((\phi_{N}\otimes\phi_{\mathcal{L}_{\psi,f}})^{r}_{(\underline{x},\underline{y})})\\
{}\\
&=&\frac{1}{q^{rs} }\ S_{r}(\mathbb{A}^{d+s}_{\mathbb{F}_{q}},N\otimes \mathcal{L}_{\psi,f})
\end{array}
$$
\noindent Par suite 
$$
\begin{array}{rcl}
L(X,M,t)&=&\mbox{exp}(\displaystyle\mathop{\sum}_{r=1}^{\infty}\frac{1}{r}(\frac{t}{qs})^{r}S_{r}(\mathbb{A}^{d+s}_{\mathbb{F}_{q}},N\otimes \mathcal{L}_{\psi,f}))\\
{}\\
&=&L(\mathbb{A}^{d+s}_{\mathbb{F}_{q}},N\otimes \mathcal{L}_{\psi,f},\frac{1}{q^{s}}t)\\
{}\\
L(X,M,t)&=&L(\mathbb{A}^{d+s}_{\mathbb{F}_{q}},pr^{\ast}(M_{0})\otimes f^{\ast}(\mathcal{L}_{\psi}),\frac{1}{q^{s}}t) \ .
\end{array}
\leqno{(3.8.3)}
$$
\noindent Ceci est le lien que nous recherchions entre les fonctions $L$ sur $X$ et les fonctions $L$ sur l'espace affine. C'est ce passage par le cas plus simple de l'espace affine que Wan utilise dans [W 4] pour la derni\`ere phase de sa preuve de la conjecture de Dwork consacr\'ee aux $F$-modules libres surconvergents de rang 1.

\section*{4. Formule des traces de Monsky g\'en\'eralis\'ee}
Soit $X=Spec\ A_{0}$ un $k$-sch\'ema lisse; on consid\`ere une $\mathcal{V}$-alg\`ebre lisse $A$ relevant $A_{0}$ comme dans (3.8) telle que $Spec\ A$ soit de dimension $n$ sur $Spec\ \mathcal{V}$. Pour tout entier $i\in \mathbb{N}$, on note $\Omega^{i}_{A^{\dag}}:= \Omega^{i}_{A^{\dag}/{\mathcal{V}}}$: puisque $A$ est lisse sur $\mathcal{V}$,  $\Omega^{i}_{A^{\dag}}$ est un $A^{\dag}$-module projectif de type fini; en particulier  $\Omega^{n}_{A^{\dag}}$ est projectif de rang 1 sur $A^{\dag}$. Le morphisme de Frobenius $F_{A^{\dag}}$ s'\'etend \`a $\Omega^{i}_{A^{\dag}}$ en tant que $\mathcal{V}$-endomorphisme $\sigma$-lin\'eaire, injectif et fini localement libre, not\'e
$$
\sigma_{i}:\Omega^{i}_{A^{\dag}} \rightarrow\Omega^{i}_{A^{\dag}}\ \ .
$$
\noindent Puisque $\sigma_{i}=\wedge^{i}\sigma_{1}$ et que $\sigma_{1}(\omega)\equiv 0$ mod $\pi$, pour $\omega \in \Omega^{1}_{A^{\dag}}$, on en d\'eduit \\

(4.1) $\qquad \sigma_{i}(\omega)\equiv 0$ mod $\pi^{i}$ , pour $\omega \in \Omega^{i}_{A^{\dag}}$ et $i\geqslant 1$.\\

\noindent Il en r\'esulte une application trace [M-W, theo 8.3]
$$
Tr_{i}:\Omega^{i}_{A^{\dag}} \rightarrow \sigma_{i}(\Omega^{i}_{A^{\dag}})
$$
\noindent telle que, pour $a\in A^{\dag}, \omega\in \Omega^{i}_{A^{\dag}}$, on ait
$$
Tr_{i}(F_{A^{\dag}}(a)\omega)=F_{A^{\dag}}(a).Tr_{i}(\omega) \ \ .
$$
\noindent De plus, en tant qu'endomorphisme de $\Omega^{i}_{A^{\dag}}$, on a d'apr\`es [M-W, theo 8.5] la relation suivante
$$
\sigma_{i}^{-1}\circ Tr_{i}\circ\sigma_{i}=[A_{0}:F_{A_{0}}(A_{0})]= q^{n} \ \ .
$$
\noindent Ainsi l'application $\sigma_{i}$ agissant sur $\Omega^{i}_{A^{\dag}}$ a un inverse \`a gauche quand on tensorise par $\mathbb{Q}$: cet \guillemotleft inverse \`a gauche\guillemotright\  est l'application $\psi_{i}$ d\'efinie par
$$
\psi_{i}=\sigma_{i}^{-1}\circ Tr_{i}: \Omega^{i}_{A^{\dag}} \rightarrow \Omega^{i}_{A^{\dag}} \ \ .
$$
L'application $\psi_{i}$ est un exemple d'op\'erateur de Dwork, i.e. d'une application $\sigma^{-1}$-lin\'eaire dans le sens suivant:
$$
\psi_{i}(F_{A^{\dag}}(a)\omega)= a\psi_{i}(\omega)\ ,\ a\in A^{\dag}\ , \ \omega\in  \Omega^{i}_{A^{\dag}} \ .
$$
\noindent Les op\'erateurs de Dwork sont des op\'erateurs \`a trace, d'o\`u leur importance ([M-W], [E-LS 1, \S\ 5.1], [W 3]).\\

Consid\'erons maintenant $(M, \phi_{M})\in F^{a}\mbox{-Mod}(A^{\dag})$ et le $A^{\dag}$-module suivant
$$
M^{\ast}:= Hom_{A^{\dag}}(M, \Omega^{n}_{A^{\dag}}) \ \ .
$$
\noindent D\'efinissons un op\'erateur de Dwork $\phi^{\ast}$ sur $M^{\ast}$ comme suit: si $f\in M^{\ast}$ et $m\in M$, on pose
$$
\phi^{\ast}(f)(m)=(\sigma_{n}^{-1}\circ Tr_{n})(f(\phi_{M}(m)))=\psi_{n}(f(\phi_{M}(m))) \ \ ;
$$
\noindent on v\'erifie que $\phi^{\ast}$ est bien un op\'erateur de Dwork, car $\psi_{n}$ en est un:
$$
\phi^{\ast}(F_{A^{\dag}}(a)f)= a\phi^{\ast}(f)\ ,\ a\in A^{\dag}\ ,  \ f\in M^{\ast} \  .
$$
\noindent Pour $i\in \llbracket 0, n\rrbracket$, soit
$$
\Omega^{i}M= M\otimes_{A^{\dag}}\Omega^{i}_{A^{\dag}} \ \ ;
$$
\noindent ce module $\Omega^{i}M$ est projectif de type fini sur $A^{\dag}$ et l'on pose
$$
\phi_{i}:= \phi_{M}\otimes \sigma_{i} \ \ ;
$$
d'apr\`es (4.1) on a \\

(4.2) $\qquad \phi_{i}\equiv 0 \  \mbox{mod}\ \pi^{i} \  $ pour $i\geqslant 1$.\\

\noindent Alors la paire $(\Omega^{i}M,\phi_{i})$ est un $F$-module surconvergent. \\
Soit
$$
M^{\ast}_{i}:= Hom_{A^{\dag}}(\Omega^{i}M, \Omega^{n}_{A^{\dag}}) \ \ .
$$
\noindent On d\'efinit de mani\`ere analogue un op\'erateur de Dwork $\phi_{i}^{\ast}$ sur $M_{i}^{\ast}$ par 
$$
\phi^{\ast}_{i}(f)(m)=(\sigma_{n}^{-1}\circ Tr_{n})(f(\phi_{i}(m)))=\psi_{n}(f(\phi_{i}(m))), m\in \Omega^{i}M, f\in M^{\ast}_{i}\ ;
$$
\noindent par (4.2) on a la congruence\\

(4.3) $\qquad \phi^{\ast}_{i}\equiv 0 \  \mbox{mod}\ \pi^{i}$ pour $i\geqslant 1$ \ .\\

\noindent Pour $i=0$, on a 
$$
(\Omega^{0}_{A^{\dag}}, \sigma_{0})=(A^{\dag}, F_{A^{\dag}}) \ , (M^{\ast}_{0}, \phi_{0}^{\ast})= (M^{\ast}, \phi^{\ast}) \ .
$$
Le $A_{K}^{\dag}$-module $M_{i}^{\ast}\otimes_{\mathcal{V}}K$ est un $K$-espace vectoriel de dimension infinie, qui est limite inductive d'une suite d'espaces de Banach $p$-adiques avec des bases orthonormales. Puisque $(M, \phi_{M})$ est surconvergent l'op\'erateur de Dwork $\phi_{i}^{\ast}$ donne, apr\`es tensorisation par $K$ sur $\mathcal{V}$, un op\'erateur nucl\'eaire sur l'espace $p$-adique $M_{i}^{\ast}\otimes_{\mathcal{V}}K$. Il r\'esulte de la th\'eorie spectrale $p$-adique de Serre [S 2] que le d\'eterminant de Fredholm $det (1-t\phi_{i}^{\ast}\vert M_{i}^{\ast}\otimes_{\mathcal{V}}K)$ est bien d\'efini et est une fonction enti\`ere de $t$ sur $K$: les coefficients de cette fonction enti\`ere sont en fait dans $\mathcal{V}$ car la construction ci-dessus est partie de $\phi_{i}$ \`a coefficients dans $\mathcal{V}$.\\

La formule des traces de Monsky g\'en\'eralis\'ee \'etablie par Wan [W 3, appendix] s'\'enonce alors comme suit:
\vskip 3mm
\noindent \textbf{Th\'eor\`eme (4.4) (Wan)  [W 3]}. \textit{Avec les notations pr\'ec\'edentes soit $(M, \phi_{M})\in F^{a}\mbox{-Mod}(A^{\dag})$ un $F$-module surconvergent sur $X/\mathbb{F}_{q}$. Alors }
$$
L(X,M,t)=\prod_{i=0}^{n}det (1-t\phi_{n-i}^{\ast}\vert M_{n-i}^{\ast}\otimes_{\mathcal{V}}K)^{(-1)^{i+1}} \ \ ,
$$
\noindent\textit{o\`u les d\'eterminants de Fredholm $det (1-t\phi_{n-i}^{\ast}\vert M_{n-i}^{\ast}\otimes_{\mathcal{V}}K)$ sont bien d\'efinis et sont des fonctions enti\`eres de $t$. En particulier $L(X,M,t)$ est $p$-adiquement  m\'eromorphe}\\

\vskip 3mm
\noindent \textit{D\'emonstration}. Si $M$ est de rang 1, ce th\'eor\`eme est cons\'equence du th\'eor\`eme 5.3 de [Mo 2]. Le cas g\'en\'eral est fourni par le th\'eor\`eme (10.10) de [W 3], dont la d\'emonstration utilise la m\'ethode des groupes de Grothendieck  de [Mo 2].$\square$\\

\noindent \textbf{Corollaire (4.5)}. \textit{Soient $X = Spec\ A_{0}$ un $k$-sch\'ema lisse, $A$ une $\mathcal{V}$-alg\`ebre lisse relevant $A_{0}$ et $(M, \phi_{M}) \in F^a\mbox{-}\textrm{Mod}(A^{\dag}_{K})$.\\
Alors $L(X, M, t)$ est $p$-adiquement m\'eromorphe.}

\vskip 3mm
\noindent \textit{D\'emonstration}. Par le lemme (3.7), on peut remplacer $M$ par un $A^{\dag}$-module libre $M' $: $(M', \phi_{M'})$ est alors ce que Wan appelle un $\sigma$-module surconvergent [W 2] [W 3]. L'extension au cas affine et lisse de la formule des traces de Monsky-Washnitzer \'etablie par Wan (cf Th\'eor\`eme (4.4)) prouve alors la m\'eromorphie de $L(M', t)$, donc celle de $L(M, t)$. On peut aussi se ramener au cas de l'espace affine [ \S\ 3.8] et utiliser la formule des traces de Monsky [Dw 5, 7(a)]. $\square$\\

\textit{Comme cons\'equence du th\'eor\`eme (4.4) nous allons, dans la suite de ce \S4, montrer que la fonction $L$ d'un $F$-module convergent \`a coefficients dans $\hat{A}$ est m\'eromorphe dans le disque unit\'e ferm\'e $\vert t \vert_{p}\leqslant 1$ [Th\'eor\`eme (4.9)].}\\

Auparavant nous allons devoir \'etablir quelques propositions.
\vskip 3mm
\noindent \textbf{Proposition (4.6)}. \textit{Avec les notations pr\'ec\'edentes soit $\mathcal{L}$ un $\hat{A}_{K}$-module projectif de type fini. Alors}\\
\begin{enumerate}

		\item[(i)] \textit{Il existe un $A^{\dag}_{K}$-module projectif de type fini $L$  et un  $\hat{A}$-module de type fini $\mathcal{M}$ tels que}\\
		
		$\mathcal{L}\simeq L\otimes_{A^{\dag}_{K}}\hat{A}_{K}$  \textit{et} $\mathcal{L}\simeq\mathcal{M}\otimes_{\hat{A}}\hat{A}_{K}$. \\
		
		\item[(ii)] \textit{Soit $M:=\mathcal{M}\cap L\subset\mathcal{L}$; alors $M$ est un $A^{\dag}$-module de type fini et on a des isomorphismes}\\
		
		 $\mathcal{M}\simeq M\otimes_{A^{\dag}}\hat{A}$ \textit{et} $L\simeq M\otimes_{A^{\dag}}A^{\dag}_{K}$ .
	 
\end{enumerate}
\vskip 3mm
\noindent \textit{D\'emonstration.}\\
\textit{ Pour(i)}. Puisque $(A^{\dag}, A^{\dag}/\pi A^{\dag})$ est un couple hens\'elien [Et 3, th\'eor\`eme 3], l'existence de $L$ r\'esulte d'un th\'eor\`eme de Elkik [E$\ell$, cor 1 p.573]. Pour $\mathcal{M}$ il suffit de prendre l'image de l'application compos\'ee $\hat{A}^{r}\displaystyle \mathop{\longrightarrow}^{can} \hat{A}^{r}_{K}\displaystyle \mathop{\longrightarrow}^s \mathcal{L}$  o\`u $s$ est une surjection et $can$ l'injection canonique.\\
\noindent\textit{Pour (ii)}. Comme on a des injections $L\hookrightarrow \mathcal{L}, \mathcal{M}\hookrightarrow \mathcal{L}$ le $A^{\dag}$-module $M:=\mathcal{M}\cap L$ est le produit fibr\'e de $\mathcal{M}$ et $L$ au-dessus de $\mathcal{L}$, donc par [F-R, prop. 4.2] on a les isomorphismes cherch\'es. $\square$

\vskip 3mm
\noindent \textbf{Corollaire (4.7)}. \textit{Avec les notations pr\'ec\'edentes soit $\mathcal{M}$ un $\hat{A}$-module de type fini tel que $\mathcal{M}_{K}:=\mathcal{M}\otimes_{\hat{A}}\hat{A}_{K}$ soit projectif. Alors}
\begin{enumerate}
\item[(i)] \textit{Il existe un $A^{\dag}$-module de type fini $M$ tel que}  $\mathcal{M}\simeq M\otimes_{A^{\dag}}\hat{A}$.
\item[(ii)] \textit{Si de plus $\mathcal{M}$ est projectif alors le $M$ du (i) est projectif de type fini.}
\end{enumerate}
\vskip 3mm
\noindent \textit{D\'emonstration.} Il suffit d'appliquer la proposition (4.6) \`a $\mathcal{L}=\mathcal{M}_{K}$; le cas projectif r\'esultant de la pleine fid\'elit\'e de $\hat{A}$ sur $A^{\dag}$ [Et 3, prop 2: (2) (ii)]. $\square$

\vskip 3mm
\noindent \textbf{Proposition (4.8)}. \textit{Avec les notations du \S4 soient $\mathcal{M}$ et $\mathcal{N}$ deux $\hat{A}$-modules de type fini tels que $\mathcal{M}_{K}$ et $\mathcal{N}_{K}$ soient projectifs sur $\hat{A}_{K}$. Soit $\psi : \mathcal{M}\rightarrow \mathcal{N}$ une application $\hat {A}$-lin\'eaire et pour tout entier $n\in \mathbb{N}$, $\psi_{n}: \mathcal{M}_{n}:=\mathcal{M}/ \pi^{n+1}\mathcal{M}\rightarrow\mathcal{N}_{n}$ sa r\'eduction modulo $\pi^{n+1}$. Consid\'erons comme dans (4.7) un $A^{\dag}$-module de type fini $M$ (resp. $N$) relevant $\mathcal {M}$ (resp. $\mathcal{N}$). Alors} 
\begin{enumerate}
\item[(i)] \textit{Il existe une application $A^{\dag}$-lin\'eaire $\varphi(n): M\rightarrow N$ relevant $\psi_{n}$.}
\item[(ii)] \textit{Si $\psi$ est surjective alors $\varphi(n)$ est surjective.}
\item[(iii)] \textit{Supposons de plus que $\mathcal{M}=\mathcal{N}$ ou que $\mathcal{N}$ est projectif sur $\hat{A}$; alors, si $\psi$ est un isomorphisme, $\varphi(n)$ est aussi un isomorphisme.}
\end{enumerate}
\vskip 3mm
\noindent \textit{D\'emonstration.}\\
\textit{ Pour(i)}. Puisque $M$ est de type fini, on rel\`eve une famille g\'en\'eratrice $(\overline{e}_{i})_{i=1,...,r}$ de $\mathcal{M}_{n}\simeq M_{n}:= M/\pi^{n+1}M$ en une famille $(e_{i})_{i=1,...,r}$ de $M$: d'apr\`es le lemme de Nakayama [Bour, AC II, \S3, \no2, cor 2 de prop 4], cette famille est g\'en\'eratrice pour $M$ car $\pi A^{\dag}\subset Rad A^{\dag}$. On note $f_{1},...,f_{r}$ des rel\`evements de $\psi_{n}(\overline{e}_{1}),..., \psi_{n}(\overline{e}_{r})$ dans $N$, et on d\'efinit $\varphi(n)$ par $A^{\dag}$-lin\'earit\'e en posant $\varphi(n)(e_{i})=f_{i}$.\\
\textit{ Pour(ii)}. On applique [Bour, AC II, \S3, \no2, cor 1 de prop 4].\\
\textit{ Pour(iii)}. Le cas $\mathcal{M}=\mathcal{N}$ r\'esulte de [Ma, theo 2.4 p. 9]. Supposons $\mathcal{N}$ projectif: si $\psi$ est un isomorphisme $\psi_{n}$ en est un aussi par $[EGA \ 0_{I} (6.7.2)]$,
 et donc $\varphi(n)$ aussi [loc. cit.] car $N$ est un $A^{\dag}$-module projectif et $M$ est de type fini sur $A^{\dag}$. $\square$\\
 
 Nous sommes \`a pr\'esent en mesure d'\'enoncer le r\'esultat  de m\'eromorphie promis avant la proposition (4.6):\\
 
 \vskip 3mm
\noindent \textbf{Th\'eor\`eme (4.9)}. \textit{Avec les notations du \S4, soit $(\mathcal{M}, \phi)\in F^{a}\mbox{-Mod}(\hat{A})$. Alors}
$$
L(X,\mathcal{M}, t)^{(-1)^{dim X -1}}\in 1+t\mathcal{V}\{t\} \ ;
$$
\noindent \textit{en particulier $L(X,\mathcal{M}, t)$ est m\'eromorphe dans le disque unit\'e ferm\'e $\vert t\vert_{p}\leqslant 1$.}\\

\vskip 3mm
\noindent \textit{D\'emonstration.} Pour tout entier $n\in \mathbb{N}$ on note
$$\psi_{n}: F^{\ast}_{A_{n}}(\mathcal{M}/\pi^{n+1}\mathcal{M})\rightarrow \mathcal{M}/\pi^{n+1}\mathcal{M}$$
 l'application obtenue en r\'eduisant $\phi$ modulo $\pi^{n+1}$. Par la proposition (4.8) on rel\`eve $\psi_{n}$ en $\phi(n): F_{A^{\dag}}^{\ast}(M)\rightarrow M$; par le d\'ebut du \S4 on en d\'eduit, pour tout entier $i, 0\leqslant i\leqslant dim\ X$, des endomorphismes $\phi(n)_{i}^{\ast}$ de $M^{\ast}_{i}\otimes_{\mathcal{V}}K$ tels que\\
 
\noindent(4.9.1)\qquad\qquad$det(1-t \phi(n)_{i}^{\ast})\in 1+\pi^{i}\mathcal{V}\llbracket\pi^{i}t\rrbracket \  [(4.3)]$\\

 \noindent et que $det(1-t \phi(n)_{i}^{\ast})$ soit une fonction enti\`ere de la variable $p$-adique $t$ [(4.4)], en particulier\\
  
 \noindent(4.9.2)\qquad\qquad$det(1-t \phi(n)_{i}^{\ast})\in 1+t\mathcal{V}\{t\}$ .\\
 
 \noindent Par le th\'eor\`eme (4.4) on en d\'eduit que\\
 
  \noindent(4.9.3)\qquad\qquad$L(X,(M, \phi(n)), t)^{(-1)^{dim X -1}}\in 1+t\mathcal{V}\{t\} $.\\
  
  \noindent Notons $L(\phi(n)):=L(X,(M, \phi(n)), t)$ et $L(\phi):=L(X,(\mathcal{M}, \phi), t)$;  par construction on a\\
  
   \noindent(4.9.4) $L( \phi(n))^{(-1)^{dim X -1}}  \mbox{mod}\  \pi^{n+1}\equiv  L( \phi)^{(-1)^{dim X -1}} \mbox{mod}\ \pi^{n+1}\in 1+t\mathcal{V}_{n}\llbracket t\rrbracket $,\\
   
 \noindent donc par (4.9.3) on en d\'eduit\\
 
  \noindent(4.9.5)\qquad\qquad $L( \phi)^{(-1)^{dim X -1}} \ \mbox{mod}\  \pi^{n+1}\in 1+t\mathcal{V}_{n}[ t] $ .\\
  
  \noindent En passant \`a la limite sur $n$ ceci ach\`eve la preuve du th\'eor\`eme (4.9). $\square$

\newpage
\section*{5. La conjecture de Dwork pour les $F$-modules surconvergents}

\subsection*{5.1. Enonc\'e de la conjecture et du th\'eor\`eme}	
	Avec les notations des \S1, 2, 3 d\'ecomposons le \guillemotleft{polyn\^{o}me caract\'{e}ristique}\guillemotright \ de $(\mathcal{M, \phi_{\mathcal{M}}}) \in F^a\mbox{-}\textrm{Mod}(\hat{A}_{K})$ au point $x \in \vert X \vert $ en
	
$$
det(\mathcal{M}_{x},t) := det(1-t\ \phi^{\textrm{deg}\ x}_{x} \vert \mathcal{M}_{x}) \in K[t]
$$
$
\qquad \qquad \qquad \qquad \qquad \qquad \quad = \displaystyle \mathop{\pi}_{j}\  (1 - a_{j,x}\  t) ,
$

\noindent o\`u les $a_{j,x}$ sont dans une cl\^oture alg\'ebrique $K^{\textrm{alg}}$ de $K$ : plus exactement les $a_{j,x}$ sont dans une extension finie (\'eventuellement ramifi\'ee) $K'(x) \subset K ^{\textrm{alg}}$ de $K(x)$ ; soient $\pi'(x)$ une uniformisante de $K'(x)$ et $\sigma_{K'(x)}$ un rel\`evement \`a $K'(x)$ de la puissance $p^a$ de $k(x)$ tel que $\sigma_{K'(x)}(\pi'(x)) = \pi'(x)$  [Et 4, 1.1]. Notons $\pi_{x} = \pi^{\textrm{deg}\ x}$ et $\textrm{ord}_{\pi_{x}}$ la valuation de $K'(x)$ normalis\'ee par 
$\textrm{ord}_{\pi_{x}}(\pi_{x})  = 1$.\\

Pour tout nombre rationnel $\alpha \in \mathbb{Q}$ on d\'efinit la partie de pente $\alpha$ du \guillemotleft polyn\^ome caract\'eristique\guillemotright\   $det(\mathcal{M}_{x}, t)$ par le produit\\

\noindent (5.1.1) $\qquad \qquad \qquad det_{\alpha}(\mathcal{M}_{x}, t) := \displaystyle \mathop{\prod}_{\textrm{ord}_{\pi_{x}}(a_{j,x}) = \alpha}\ (1 - a_{j,x}\ t)$ .\\

\vskip 2mm
\noindent Si $a_{j,x}$ est l'inverse d'une racine de $det(\mathcal{M}_{x}, t)$, alors $\sigma_{K'(x)} (a_{j,x})$ en est une aussi par le m\^eme argument que pour la d\'emonstration du lemme (2.5), et puisque $\sigma_{K'(x)} : K'(x) \rightarrow K'(x)$ est une extension isom\'etrique on a $\textrm{ord}_{\pi_{x}}(a_{j,x}) = \textrm{ord}_{\pi_{x}}(\sigma_{K'(x)}(a_{jx}))$ ; par cons\'equent $det_{\alpha}(\mathcal{M}_{x}, t)$ est \`a coefficients dans $K$. Comme $\mathcal{M}$ est un $\hat{A}_{K}$-module projectif de type fini, il existe $\alpha_{0} \in \mathbb{Q}$ tel que pour tout $x \in \vert Spec\ A_{0} \vert$ et tout $\alpha \in \mathbb{Q}$, $\alpha < \alpha_{0}$, on ait $det_{\alpha}(\mathcal{M}_{x}, t) = 1$. D'autre part $det(\mathcal{M}_{x}, t)$ s'exprime par un produit fini\\

\noindent (5.1.2) $\qquad \qquad \qquad det(1 - t\ \phi^{\textrm{deg}\ x}_{x}, \mathcal{M}_{x}) = \displaystyle \mathop{\prod}_{\alpha \in \mathbb{Q}} det_{\alpha}(\mathcal{M}_{x}, t)$ . \\

La partie de pente $\alpha$ de la fonction $L(X, \mathcal{M}, t)$ est d\'efinie par l'expression\\

\noindent (5.1.3) $\qquad \qquad \qquad L_{\alpha}(X, \mathcal{M}, t) := \displaystyle \mathop{\prod}_{x \in \vert X \vert}  det_{\alpha}(\mathcal{M}_{x}, t^{\textrm{deg}\ x})^{-1} \in K[[t]]$ .\\

\noindent Comme $det_{\alpha}(\mathcal{M}_{x}, t) = 1$ pour tout $\alpha < \alpha_{0}$ et tout $x \in \vert X \vert $, le th\'eor\`eme de sp\'ecialisation de Grothendieck [K 2, (2.3)] montre qu'il n'existe qu'un nombre fini de possibilit\'es pour les pentes de Newton de $\mathcal{M}_{x}$ (donc pour les $\alpha$ ci-dessus) pour tous les $x \in \vert X \vert $.\\
\noindent D'apr\`es (5.1.2) on a ainsi la relation\\

\noindent (5.1.4) $\qquad \qquad \qquad L(X, \mathcal{M}, t) = \displaystyle \mathop{\prod}_{\alpha \in \mathbb{Q}} L_{\alpha}(X, \mathcal{M}, t)$ ;\\

\noindent o\`u le produit (5.4) est fini d'apr\`es la remarque pr\'ec\'edente.\\

Pour $r \in \mathbb{N}^{\ast}$ et $\alpha \in \mathbb{Q}$, on d\'efinit plus g\'en\'eralement\\

\noindent (5.1.5) $\qquad \qquad \qquad det^{(r)}(\mathcal{M}_{x}, t) := \displaystyle \mathop{\prod}_{j}\  (1 - (a_{j,x})^r\ t),$\\

\vskip 3mm
\noindent (5.1.6) $\qquad \qquad \qquad det^{(r)}_{\alpha}(\mathcal{M}_{x}, t) := \displaystyle \mathop{\prod}_{\textrm{ord}_{\pi_{x}}(a_{j,x}) = \alpha}\ (1 - (a_{j,x})^r\ t) $ ,\\

\vskip 3mm
\noindent (5.1.7) $\qquad \quad   L^{(r)}(X, \mathcal{M}, t) := \displaystyle \mathop{\prod}_{x \in \vert X \vert}\  det^{(r)}  (\mathcal{M}_{x}, t^{\textrm{deg}\ x})^{-1} \in K[[t]]$ , \\

\vskip 3mm
\noindent (5.1.8) $\qquad   \quad L^{(r)}_{\alpha}(X, \mathcal{M}, t) := \displaystyle \mathop{\prod}_{x \in \vert X \vert}\  det^{(r)}_{\alpha}  (\mathcal{M}_{x}, t^{\textrm{deg}\ x})^{-1} \in K[[t]]$ .\\

\noindent On a encore :\\

\noindent (5.1.9) $\qquad   \quad L^{(r)} (X, \mathcal{M}, t) = \displaystyle \mathop{\prod}_{\alpha \in \mathbb{Q}}\ L^{(r)}_{\alpha}  (X, \mathcal{M}, t)$ .\\

\noindent Plus pr\'ecis\'ement, gr\^ace \`a l'expression \'etablie par Wan [W 2, lemma 4.4], on a en fait\\

\noindent (5.1.10) $\qquad   \quad L^{(r)} (X, \mathcal{M}, t) = \displaystyle \mathop{\prod}_{i \geqslant 1}\ L(X, \textrm{Sym}^{r-i}\ \mathcal{M} \ \otimes \displaystyle \mathop{\Lambda}^i\ \mathcal{M}, t)^{i \times (-1)^{i-1}}$ ,\\

\vskip 3mm
\noindent (5.1.11) $\qquad   \quad L^{(r)}_{\alpha} (X, \mathcal{M}, t) = \displaystyle \mathop{\prod}_{i \geqslant 1}\ L_{\alpha} (X, \textrm{Sym}^{r-i}\ \mathcal{M} \ \otimes \displaystyle \mathop{\Lambda}^i\ \mathcal{M}, t)^{i \times (-1)^{i-1}}$ .\\

\vskip 3 mm
Les d\'efinitions pr\'ec\'edentes s'\'etendent \`a $M \in F^a\mbox{-}\textrm{Mod}(A^{\dag}_{K})$.\\

Nous pouvons \`a pr\'esent \'enoncer la conjecture de Dwork g\'en\'eralis\'ee pour les $F^a$-modules surconvergents.

\vskip 3mm
\noindent (5.1.12) \textit{\textbf{Conjecture (Dwork)}. Soit $M \in F^a\mbox{-}\textrm{Mod}(A^{\dag}_{K})$. Alors, pour tout nombre rationnel $\alpha \in \mathbb{Q}$ et tout entier $r \in \mathbb{N}^{\ast}$, la fonction $L^{(r)}_{\alpha}(X, M, t)$ est $p$-adiquement m\'eromorphe.}\\

Wan a prouv\'e cette conjecture dans une s\'erie de trois articles [W 2][W 3][W 4] (cf (5.1.13) plus bas).\\

\textbf{Notre objectif} dans ce \S\ 5 est d'apporter la pr\'ecision suivante \`a la preuve de Wan (pour une formulation plus pr\'ecise cf \S\  5.4). A isog\'enie pr\`es la fonction $L$-unit\'e $L_{0}(X, M, t)$ est la fonction $L(X,\mathcal{M}_{0}, t)$ du sous-$F$-module unit\'e $\mathcal{M}_{0}$ du compl\'et\'e $\mathcal{M}$ de $M$; cependant, bien que cette fonction $L(X,\mathcal{M}_{0}, t)$ soit $p$-adiquement m\'eromorphe par Wan, nous montrons qu'\textbf{il n'existe pas en g\'en\'eral de sous-module $M_{0}$ de $M$ dont $\mathcal{M}_{0}$ serait le compl\'et\'e  (m\^eme si le Frobenius de $\mathcal{M}_{0}$ est surconvergent)} et auquel on pourrait appliquer la formule des traces de Monsky du \S\ 4 pour obtenir la m\'eromorphie: nous en donnerons deux contre-exemples issus de familles de courbes elliptiques ordinaires, dont l'un est la famille de Legendre.\\
 Apr\`es un bref rappel au \S\ 5.2 du contexte historique o\`u Dwork a pos\'e sa conjecture et l'a r\'esolue pour la famille de Legendre des courbes elliptiques ordinaires, nous exposerons au \S\ 5.3 un des aspects de la preuve de Wan ("d\'ecomposition de Hodge-Newton,
th\'eor\`eme d'isog\'enie de Katz") qui nous servira \`a expliciter nos contre-exemples au \S\ 5.4.\\

 Cette conjecture (5.1.12) est \'etablie ci-dessous en (5.1.14) comme cons\'equence du th\'eor\`eme suivant de Wan:\\

\noindent \textbf{Th\'eor\`eme (5.1.13) (Wan)[W 3, theo 1.1]}. \textit{Soient $X$ une vari\'et\'e affine et lisse d\'efinie sur un corps fini $\mathbb{F}_{q}$ de caract\'eristique $p>0$, et $(M, \phi)$ un $\sigma $-module surconvergent sur $X/ \mathbb{F}_{q}$. Alors, pour tout nombre rationnel $s$, la fonction $L$ de pente pure $s, L_{s}(\phi, t)$, attach\'ee \`a $(M, \phi)$ est $p$-adiquement m\'eromorphe.}\\

Rappelons la terminologie de Wan. Si $X= Spec\ A_{0}$, on rel\`eve $A_{0}$ en une $\mathcal{V}$-alg\`ebre lisse $A$, on note $A^{\dag}$ (resp.$\hat{A}$) sa compl\'et\'ee faible (resp. son s\'epar\'e compl\'et\'e) et $F_{A^{\dag}}$ (resp. $F_{\hat{A}}$) un rel\`evement du Frobenius de $A_{0}$ \`a $A^{\dag}$ (resp. \`a $\hat{A}$). Un $\sigma$-module surconvergent est ce que nous avons appel\'e un $F$-module surconvergent libre, i.e. la donn\'ee d'un $A^{\dag}$-module libre de type fini $M$ muni d'un morphisme de Frobenius
$$
\phi: F_{A^{\dag}}^{\ast} M=: M^{\sigma} \rightarrow M \ \ .
$$
\noindent L'expression $\ L_{s}(\phi,t):= L_{s}(X, (M,\phi),t)$ a \'et\'e d\'efinie en (5.1.3).\\

La d\'emonstration par Wan de ce th\'eor\`eme (5.1.13) s'effectue en deux phases:\\
\textbf{La premi\`ere phase}, regroup\'ee dans le premier article [W 3], est une phase de r\'eduction:
\begin{itemize}
\item de la pente $s$ \`a la pente z\'ero,
\item d'une pente z\'ero de rang $r$ \`a une pente z\'ero de rang 1,
\item de $X$ affine et lisse sur $\mathbb{F}_{q}$ \`a l'espace affine $\mathbb{A}^{n}_{\mathbb{F}_{q}}$.
\end{itemize}

Cette partie est de nature alg\'ebrique via:
\begin{itemize}
\item la formule des traces de Monsky,
\item le th\'eor\`eme de sp\'ecialisation de Grothendieck,
\item la d\'ecomposition de Hodge-Newton,
\item le th\'eor\`eme d'isog\'enie de Katz.
\end{itemize}

\noindent\textbf{La deuxi\`eme phase}, regroup\'ee au sein du deuxi\`eme article [W 4], est de nature plus analytique:
\begin{itemize}
\item par le travail avec des modules de rang infini,
\item par des processus de passage \`a la limite dans des familles de fonctions $L$ m\'eromorphes. $\square$
\end{itemize}
 \vskip 3mm
 
Avant d'examiner plus en d\'etail au \S 5.3 les parties "d\'ecomposition de Hodge-Newton,
th\'eor\`eme d'isog\'enie de Katz" de la preuve de Wan, le th\'eor\`eme (5.1.13) nous permet d'\'enoncer:\\

\noindent \textbf{Corollaire (5.1.14)}. \textit{Soient $X = Spec(A_{0})$ un $k$-sch\'ema lisse, $A$ une $\mathcal{V}$-alg\`ebre  lisse relevant $A_{0}$ et $M \in F^a\mbox{-}\textrm{Mod}(A^{\dag}_{K})$. Alors, pour tout  $\alpha \in \mathbb{Q}$ et tout $r \in \mathbb{N}^{\ast}$, les fonctions $L^{(r)}_{\alpha}(X, M, t)$, et $L^{(r)}(X, M, t)$ sont $p$-adiquement m\'eromorphes.}
 \vskip 3mm
\noindent \textit{D\'emonstration.}
\noindent  Par stabilit\'e de la cat\'egorie $F^a\mbox{-}\textrm{Mod}(A^{\dag}_{K})$ par puissances sym\'etriques et ext\'erieures [W 2, \S\ 3] la conjecture (5.1.12) se ram\`ene au cas $r = 1$. Par la m\^eme d\'emonstration que celle du lemme (3.7) on peut supposer qu'il existe un $F^a$-module libre surconvergent $(M', \phi_{M'}) \in F^a\mbox{-}\textrm{Mod}(A^{\dag})$ tel que

$$
(M, \phi_{M}) = (M', \pi^{\beta} \phi_{M'}) \otimes_{A^{\dag}} A^{\dag}_{K}\  \textrm{pour un}\  \beta \in \mathbb{Z}\ ,
$$

\noindent et $ \qquad \qquad \quad L(M,t) = L(M', \pi^{\beta} t)$ .\\

\noindent On est ramen\'e \`a montrer que $L_{\alpha}(Spec\ A_{0}, M', t)$ est $p$-adiquement m\'eromorphe, ce qui r\'esulte du th\'eor\`eme de Wan (5.1.13). $\square$\\

\subsection*{5.2. Petit rappel historique}

La conjecture de Dwork, qui date du d\'ebut des ann\'ees 70 [Dw 5], est venue des tentatives de Dwork pour comprendre les variations analytiques $p$-adiques des parties \guillemotleft pures\guillemotright\  d'une vari\'et\'e lorsque celle-ci \'evolue au sein d'une famille sur un corps fini de caract\'eristique $p>0$, i.e de comprendre les variations des z\'eros et des p\^oles de la fonction z\^eta d'une famille de vari\'et\'es lorsque le param\`etre varie.\\

Dans le cas de la famille de Legendre des courbes elliptiques ordinaires que Dwork \'etudie en d\'etail dans [Dw 2], Dwork avait prouv\'e sa conjecture sur la m\'eromorphie de la fonction $L$ associ\'ee [Dw 3]: dans ce cas le point cl\'e de sa preuve est l'existence d'un \guillemotleft rel\`evement excellent\guillemotright\  (excellent lifting) du Frobenius ( d\'efini dans [Dw 4, \S5], [Dw 5, \S2] comme laissant stable le premier cran de la filtration de Hodge par le Frobenius agissant sur la cohomologie de de Rham relative d'une famille de vari\'et\'es) qui permet d'appliquer la formule des traces de Monsky et d'en d\'eduire la m\'eromorphie [Dw 3]; voir \`a ce propos l'article de Katz [K 3, \S A3] o\`u cette stabilit\'e du premier cran de la filtration de Hodge par le Frobenius est automatique lorsque l'on choisit le \guillemotleft rel\`evement canonique\guillemotright\ du Frobenius ( le \guillemotleft rel\`evement canonique\guillemotright\ est donc un \guillemotleft rel\`evement excellent\guillemotright\ ). Dwork a \'egalement prouv\'e sa conjecture pour certaines familles de surfaces $K3$ ordinaires [Dw 5]. L'ennui de cette m\'ethode c'est qu'elle repose sur l'existence des \guillemotleft rel\`evements excellents\guillemotright\  du Frobenius, existence qui n'est pas assur\'ee dans le cas g\'en\'eral comme l'a prouv\'e Sperber [Sp, \S3].\\

On peut se demander \'egalement l'int\'er\^et qu'il y a de conna\^itre la m\'eromorphie de la fonction $L$ dans la conjecture de Dwork: c'est qu'elle donne des estimations $p$-adiques de sommes de caract\`eres, estimations qui g\'en\'eralisent celles obtenues par les conjectures de Weil. Illustrons cette remarque sur un exemple. Consid\'erons un morphisme propre et lisse $f: X \rightarrow Y$ de vari\'et\'es d\'efinies sur $\mathbb{F}_{q}$: alors la fonction z\^eta unit\'e de cette famille est donn\'ee par [E-LS 2]:

$$
Z_{u}(Y/X, t)=  \displaystyle \mathop{\prod}_{m} L(X, R^{m}f_{\ast}\mathbb{Q}_{p},t)^{(-1)^{m}}
$$                       
\noindent et Dwork conjecture, pour tout entier $m$, la m\'eromorphie de la fonction $L(X, R^{m}f_{\ast}\mathbb{Q}_{p},t)$. Or $ L(X, R^{m}f_{\ast}\mathbb{Z}_{p},t)$ est donn\'ee par la fonction g\'en\'eratrice
$$
L(X, R^{m}f_{\ast}\mathbb{Z}_{p},t)= exp\ (\sum_{k=1}^{\infty} \frac{S_{k}(R^{m}f_{\ast}\mathbb{Z}_{p})}{k} t^{k}) \ \  ,
\leqno{(5.2.1)}
$$
\noindent o\`u $S_{k}(R^{m}f_{\ast}\mathbb{Z}_{p})$ est une somme de caract\`eres sur les points ferm\'es de $X$. La m\'eromorphie (suppos\'ee) de la fonction $L(X, R^{m}f_{\ast}\mathbb{Z}_{p},t)$ et le th\'eor\`eme de factorisation $p$-adique de Weierstra\ss  \ montrent qu'il existe des entiers $p$-adiques $\alpha_{i} (1\leqslant i < + \infty)$ et des entiers $p$-adiques $\beta_{j} (1\leqslant j < + \infty)$, alg\'ebriques sur $K$ [Ro, chap 6, 2.2, theo 2 (a),p 312], tels que
$$
L(X, R^{m}f_{\ast}\mathbb{Z}_{p},t)= \frac{ \displaystyle \mathop{\prod}_{i \geqslant 1}(1- \alpha_{i}t)}{\displaystyle \mathop{\prod}_{j \geqslant 1}(1- \beta_{j}t)}
\leqno{(5.2.2)}
$$
\noindent avec
$$
\lim_{i \to + \infty} \alpha_{i} = \lim_{j \to + \infty} \beta_{j} = 0 \ \ .
$$
\noindent En prenant les logarithmes dans (5.2.1) et (5.2.2), il en r\'esulte, pour tout entier $k$, l'\'egalit\'e
$$
S_{k}= \sum_{j \geqslant 1} \beta_{j}^{k} - \sum_{i \geqslant 1} \alpha_{i}^{k} \ \ .
\leqno{(5.2.3)}
$$
\noindent Cette formule g\'en\'eralise la formule classique du nombre de points rationnels de $X$ dans $\mathbb{F}_{q^{k}}$ qui r\'esulte des conjectures de Weil, et dans ce cas la somme dans (5.2.3) ne comporte qu'un nombre fini de termes non nuls. Dans le cas g\'en\'eral, si l'on ne consid\`ere qu'un nombre fini de termes $\beta_{j}$ et $\alpha _{i}$ dans (5.2.3), on obtient une formule asymptotique $p$-adique pour des sommes de caract\`eres $p$-adiques: plus on prend de termes et plus la pr\'ecision augmente.\\

\newpage
\subsection*{5.3. $F$-modules convergents ordinaires}
Dans ce \S5.3. on supposera simplement que $k$ est un corps parfait de caract\'eristique $p>0$, $X= Spec\ A_{0}$ est un $k$-sch\'ema lisse et $A$ une $\mathcal{V}$-alg\`ebre lisse relevant $A_{0}$.
\subsubsection*{5.3.1. D\'efinitions}
 Soit $(\mathcal{M}, \phi_{\mathcal{M}})\in F^a\mbox{-}\textrm{Mod}(\hat{A})$ un $F$-module convergent . \\

A la suite de Katz et Deligne [K 2, II, 2.4, Rks p 148], [De$\ell$ 2], on dit que \textbf{$\mathcal{M}$ est un $F$-module convergent ordinaire de niveau $m$} s'il existe une filtration de $\mathcal{M}$ par des sous $F$-modules convergents ( donc localement libres de type fini)\\

\noindent (5.3.1.1) $ \qquad \qquad 0 \subset \mathcal{M}_{0}Ê\subset \mathcal{M}_{1} \subset ... \subset \mathcal{M}_{i-1} \subset \mathcal{M}_{i} \subset ... \subset \mathcal{M}_{m} = \mathcal{M}$\\

\noindent tels que $(\mathcal{M}_{i}/\mathcal{M}_{i-1}, \phi_{\mathcal{M}_{i}/\mathcal{M}_{i-1}})$, o\`u $ \phi_{\mathcal{M}_{i}/\mathcal{M}_{i-1}}$est induit par $ \phi_{\mathcal{M}}$, soit de la forme $ (\mathcal{U}_{i}, \pi^{i}\  \phi_{i})$, o\`u $ (\mathcal{U}_{i},  \phi_{i})$ est un $F$-module convergent unit\'e (donc localement libre de type fini). On dit que $\mathcal{M} \in F^a\mbox{-}\textrm{Mod}(\hat{A})$ est $\textbf{ordinaire}$ s'il existe $m \in \mathbb{N}$ tel que $\mathcal{M}$ soit ordinaire de niveau $m$.\\

Pour les d\'efinitions et propri\'et\'es des polygones de Hodge et de Newton de $\mathcal{M}$ nous renvoyons le lecteur \`a [K 2, I, 1.2, 1.3 et 2.3 p 142].\\

Puisque $k$ est parfait, si $\mathcal{M} \in F^a\mbox{-}\textrm{Mod}(\hat{A})$ est ordinaire, il est clair qu'en tout point $x \in \vert X\vert$ les polygones de Hodge et de Newton de $\mathcal{M}$ co¬\"{\i}ncident et qu'ils sont constants, i.e. ind\'ependants de $x \in \vert X\vert$ [loc. cit.]. La r\'eciproque est vraie puisque $X$ est un sch\'ema sur un corps parfait $k$ de caract\'eristique $p > 0$ et $X$ est lisse sur $k$ [K 2, II, 2.4 Rks p 148], [W 2, lemma 3.6], [W 3, theo 7.2, cor 7.3].
\subsubsection*{5.3.2. Rappels des r\'esultats de Wan sur la d\'ecomposition de Hodge-Newton et le th\'eor\`eme d'isog\'enie}
Nous utiliserons en 5.4 les r\'esultats suivants de Wan.
\vskip 3mm
\noindent \textbf{Lemme (5.3.2.1) (Wan) [W3, 4.5]}. \textit {Avec les notations de 5.3, soit $(M, \phi_{M}) \in F^a\mbox{-}\textrm{Modlib}\ (A^{\dag})$ un $F$-module libre surconvergent. Alors, quitte \`a r\'etr\'ecir $X$ si n\'ecessaire, le $F$-module surconvergent $(M, \phi_{M})$  admet une filtration $\phi_{M}$-stable par des $A^{\dag}$-modules libres de type fini \`a quotients libres de type fini}
$$  0\subset N_{0}\subset M$$
\noindent \textit {tels que la restriction $\phi_{{N}_{0}}$ de  $\phi_{M}$ \`a $N_{0}$ est nilpotente et le quotient  $M/N_{0}$ est un $F$-module libre surconvergent dont le Frobenius est injectif.}\\

Dans le cas o\`u l'on part d'un $F$-module libre convergent ordinaire on a le r\'esultat suivant de Wan qui am\'eliore (5.3.1.1):
\vskip 3mm
\noindent \textbf{Proposition (5.3.2.2)(Wan)[W 3, 4.9, 4.1]}. \textit {Avec les notations de 5.3 soit $(\mathcal{M},\phi_{\mathcal{M}}) \in F^a\mbox{-}\textrm{Modlib}\ (\hat{A})$ un $F$-module libre convergent ordinaire de niveau $m$. Alors  $(\mathcal{M},\phi_{\mathcal{M}})$ admet une filtration finie $\phi_{\mathcal{M}}$-stable par des $\hat{A}$-modules libres de type fini \`a quotients libres de type fini}
$$ 0\subset \mathcal{M}_{0}Ê\subset \mathcal{M}_{1} \subset ... \subset \mathcal{M}_{i-1} \subset \mathcal{M}_{i} \subset ... \subset \mathcal{M}_{m} = \mathcal{M}$$
\noindent \textit {tels que
\begin{itemize}
\item [(a)] Le quotient  $\mathcal{M}_{i}/\mathcal{M}_{i-1}$ est un $F$-module libre convergent de la forme $(\mathcal{U}_{i},\pi^{i} \phi_{i})$ o\`u $(\mathcal{U}_{i}, \phi_{i}) \in F^a\mbox{-}\textrm{Modlib}\ (A^{\dag})^{0}$ est un $F$-module libre convergent unit\'e pour chaque entier $i\in\llbracket 0, m\rrbracket$ .
\item[(b)]  Pour chaque entier $i\in\llbracket 0, m\rrbracket$ on a une d\'ecomposition
$$ \mathcal{M}=\mathcal{M}_{i}\oplus \mathcal{M}_{(i+1)}$$
dans laquelle $\mathcal{M}_{(i+1)}$ est un sous $\hat{A}$-module libre de type fini de $\mathcal{M}$ tel que $\phi_{\mathcal{M}}(F^{\ast}_{\hat{A}}\mathcal{M}_{(i+1)})\subset \pi^{i+1}\mathcal{M}$; en particulier $\mathcal{M}/\mathcal{M}_{i}$ est muni d'un Frobenius divisible par $\pi^{i+1}$.
\item[(c)] $\phi_{\mathcal{M}}$ est injectif.
 \end{itemize}
 }
 \vskip 3mm
\noindent \textit{D\'emonstration.} Seul le point (c) reste \`a prouver. Or l'injectivit\'e de $\phi_{\mathcal{M}}$ est \'equivalente \`a la non nullit\'e de det($\phi_{\mathcal{M}}$); celle-ci est claire car 
$$det(\phi_{\mathcal{M}})=\prod_{i=0}^{i=m}det (\pi^{i}\phi_{i})\not=0.\qquad  \square$$

\vskip 6mm
Voici la version de Wan du th\'eor\`eme d'isog\'enie de Katz:
\vskip 3mm
\noindent \textbf{Th\'eor\`eme (5.3.2.3)(Wan)[W3, 7.2]}. \textit {Avec les notations de 5.3, soient $(M, \phi_{M}) \in F^a\mbox{-}\textrm{Modlib}\ (A^{\dag})$ un $F$-module libre surconvergent \`a \guillemotleft pentes de Newton\guillemotright \  des entiers naturels et $\mathcal{M}= M\otimes_{A^{\dag}} \hat{A} \in F^a\mbox{-}\textrm{Modlib}\ (\hat{A})$ le $F$-module convergent associ\'e; on suppose que $ \phi_{M}$ est injectif. Alors, quitte \`a r\'etr\'ecir $X$ si n\'ecessaire, le $F$-module surconvergent $(M, \phi_{M})$ est isog\`ene \`a un $F$-module surconvergent $(M', \phi_{M'})$ dont le $F$-module convergent associ\'e $(\mathcal{M}',\phi_{\mathcal{M}'})$ est ordinaire. C'est-\`a-dire que le $F$-module convergent $(\mathcal{M}',\phi_{\mathcal{M}'})$ admet une filtration finie $\phi_{\mathcal{M}'}$-stable par des $\hat{A}$-modules libres de type fini \`a quotients libres de type fini}\\

$ \qquad\qquad 0\subset \mathcal{M}'_{0}Ê\subset \mathcal{M}'_{1} \subset ... \subset \mathcal{M}'_{i-1} \subset \mathcal{M}'_{i} \subset ... \subset \mathcal{M}'_{m} = \mathcal{M}'$\\

\noindent \textit {tels que\\
\begin{itemize}
\item [(a)] Le quotient  $\mathcal{M}'_{i}/\mathcal{M}'_{i-1}$ est un $F$-module libre convergent de la forme $(\mathcal{U}_{i},\pi^{i} \phi_{i})$ o\`u $(\mathcal{U}_{i}, \phi_{i}) \in F^a\mbox{-}\textrm{Modlib}\ (A^{\dag})^{0}$ est un $F$-module libre convergent unit\'e pour chaque entier $i\in\llbracket 0, m\rrbracket$ .\\
\item[(b)]  Pour chaque entier $i\in\llbracket 0, m\rrbracket$ on a une d\'ecomposition
$$\mathcal{M}'=\mathcal{M}'_{i}\oplus \mathcal{M}'_{(i+1)}$$
 dans laquelle $\mathcal{M}'_{(i+1)}$ est un sous $\hat{A}$-module libre de type fini de $\mathcal{M}'$ tel que $\phi_{\mathcal{M}'}(F^{\ast}_{\hat{A}}\mathcal{M}'_{(i+1)})\subset \pi^{i+1}\mathcal{M}'$.
 \end{itemize}
 }
 \vskip 3mm
\noindent \textit{D\'emonstration.} C'est [W 3, 7.2] hormis le $F$-module surconvergent $M'$ dont l'existence est mentionn\'ee dans la preuve de [W 2, 7.2] par Wan, et l'assertion $(b)$ qui r\'esulte de (5.3.2.2)(b).\ $\square$

\subsection*{5.4 Non surconvergence de la partie unit\'e: deux contre-exemples}
\subsubsection*{5.4.1 Pentes des fonctions $L$ }

Soit $X= Spec A_{0}$ un $k$-sch\'ema lisse, $A$ une $\mathcal{V}$-alg\`ebre lisse relevant $A_{0}$, $(M, \phi_{M}) \in F^a\mbox{-}\textrm{Mod}\ (A^{\dag})$ un $F$-module surconvergent et $\mathcal{M}= M\otimes_{A^{\dag}} \hat{A} \in F^a\mbox{-}\textrm{Mod}\ (\hat{A})$ le $F$-module convergent associ\'e. Pour la m\'eromorphie de $L^{(r)}_{\alpha}(X, M, t)$, et $L^{(r)}(X, M, t)$, il suffit de traiter le cas $r=1$ et $M$ libre [cf la preuve de (5.1.14)]. Par le th\'eor\`eme de sp\'ecialisation de Grothendieck on a vu en (5.1) que $M$ (et donc $\mathcal{M}$) n'a qu'un nombre fini de \guillemotleft pentes de Newton\guillemotright \  $\alpha\in\mathbb{Q}^{+}$. Quitte \`a passer \`a une extension finie totalement ramifi\'ee de $\mathcal{V}$ on peut supposer que toutes les  \guillemotleft pentes de Newton\guillemotright \  de $M$ sont des entiers $\alpha\in\mathbb{N}$.\\

Soit $N_{0}\subset M$ comme dans le lemmme (5.3.2.1). Puisque
$$
L(X, M, t)=:L(M,t)=L(N_{0},t)\times L(M/N_{0},t)=L(M/N_{0},t)\ ,
$$
\noindent car $\phi_{N_{0}}$ est nilpotent, on est ramen\'e au cas o\`u $\phi_{M}$ est injectif.

On peut alors appliquer le th\'eor\`eme d'isog\'enie de Katz (5.3.2.3) \`a notre $(M, \phi_{M}) \in F^a\mbox{-}\textrm{Modlib}\ (A^{\dag})$ \`a \guillemotleft pentes de Newton\guillemotright \  des entiers $\alpha\in\mathbb{N}$; soit $M'$ isog\`ene \`a $M$ comme dans ce th\'eor\`eme (5.3.2.3). Alors $(M_{K}, \phi_{M_{K}})$ est isomorphe \`a $(M'_{K}, \phi_{M'_{K}})$, donc ils ont m\^eme fonction $L$
$$L(X,M,t)=L(X,M_{K},t)=L(X,M'_{K},t)=L(X,M',t)$$
\noindent et m\^emes fonctions $L_{\alpha}$
$$L_{\alpha}(M,t)=L_{\alpha}(M_{K},t)=L_{\alpha}(M'_{K},t)=L_{\alpha}(M',t)=L_{\alpha}(\mathcal{M}',t)= L( \mathcal{M}'_{\alpha}/\mathcal{M}'_{\alpha-1}, t)\ . $$

Pour la pente z\'ero, notons $\mathcal{M}'_{0}\subset \mathcal{M}'$ le sous-$F$-module convergent unit\'e de $\mathcal{M}'$ ; on a donc
$$
L_{0}(M,t)= L_{0}(\mathcal{M}',t)=L( \mathcal{M}'_{0}, t)\ .
$$
\noindent S'il existait un sous-$F$-module surconvergent unit\'e $M'_{0}$ de $M'$ tel que 
$$M'\otimes_{A^{\dag}}\hat{A}=\mathcal{M}_{0}'$$
 on aurait
$$
L_{0}(M,t)= L( M'_{0}, t)\ ;$$
la formule de traces de Monsky g\'en\'eralis\'ee [Th\'eor\`eme (4.4)] prouverait alors que cette derni\`ere fonction est $p$-adiquement m\'eromorphe.\textbf{Toute la suite de cet article est consacr\'ee \`a  prouver qu'un tel $M'_{0}$ n'existe pas en g\'en\'eral}; pour ce faire on peut supposer que \textbf{ $(M, \phi_{M})$ est un $F$-module surconvergent ordinaire, i.e. dont le $F$-module convergent associ\'e     $(\mathcal{M},\phi_{\mathcal{M}})$ est ordinaire}.

\subsubsection*{5.4.2 La pente z\'ero : partie unit\'e des $F$-modules convergents ordinaires}

Au passage notons la caract\'erisation suivante de la partie unit\'e d'un $F$-module convergent ordinaire:
\vskip 3mm
\noindent \textbf{Lemme (5.4.2.1)}. \textit {Avec les notations de 5.3, soit $(\mathcal{M},\phi_{\mathcal{M}}) \in F^a\mbox{-}\textrm{Modlib}\ (\hat{A})$ un $F$-module libre convergent ordinaire, donc muni d'une filtration satisfaisant \`a (5.3.2.2). Pour chaque entier $n\in \mathbb{N}$ on note 
 $ \phi^{n}_{\mathcal{M}} :  F^{n\ast}_{\hat{A}}(\mathcal{M}) \rightarrow \mathcal{M}$
 l'it\'er\'e de $\phi_{\mathcal{M}}$. Alors, on a}
$$
\mathcal{M}_{0}= \bigcap_{n}\mbox{Im} \lbrace\phi^{n}_{\mathcal{M}} :  F^{n\ast}_{\hat{A}}(\mathcal{M}) \rightarrow \mathcal{M}\rbrace\ .
$$
 \vskip 3mm
\noindent \textit{Preuve du lemme.} Puisque $\phi_{\mathcal{M}_{0}}$ est bijectif on a une inclusion \'evidente

$$
\mathcal{M}_{0}\subset \bigcap_{n}\mbox{Im} \lbrace\phi^{n}_{\mathcal{M}} :  F^{n\ast}_{\hat{A}}(\mathcal{M}) \rightarrow \mathcal{M}\rbrace\ .
$$
\noindent Pour prouver l'\'egalit\'e il suffit d'apr\`es [Prop. (1.4)] d'\'etablir la surjectivit\'e en tout point ferm\'e $x$ de $X$, i.e., en utilisant les notations du \S1, que l'on a un isomorphisme
$$
\hat{\tau}(x)^{\ast}(\mathcal{M}_{0})=:(\mathcal{M}_{0})_{x}\simeq (\bigcap_{n}\mbox{Im} \lbrace\phi^{n}_{\mathcal{M}} :  F^{n\ast}_{\hat{A}}(\mathcal{M}) \rightarrow \mathcal{M}\rbrace)_{x}\ .
$$
\noindent Or la compos\'ee des applications canoniques
$$
(\mathcal{M}_{0})_{x}\rightarrow (\bigcap_{n}\mbox{Im}\  \phi^{n}_{\mathcal{M}})_{x}\rightarrow \bigcap_{n}(\mbox{Im}\  \phi^{n}_{\mathcal{M}_{x}})
$$
\noindent et des \'egalit\'es
$$
\bigcap_{n}(\mbox{Im}\  \phi^{n}_{\mathcal{M}_{x}})=(\mathcal{M}_{x})_{0}\ [De\ell\ 2; 1.3.2, (1.3.3.3)\ et\ Rq\ 1.2.5]\ ,
$$
$$
(\mathcal{M}_{x})_{0}=(\mathcal{M}_{0})_{x}\ [K 2]\ [W 3; 4.12]\ ,
$$
\noindent est l'identit\'e de $(\mathcal{M}_{0})_{x}$. De m\^eme la compos\'ee de l'application canonique
$$
 (\bigcap_{n}\mbox{Im}\  \phi^{n}_{\mathcal{M}})_{x}\rightarrow \bigcap_{n}(\mbox{Im}\  \phi^{n}_{\mathcal{M}_{x}})\ ,
$$
\noindent  des \'egalit\'es
$$
\bigcap_{n}(\mbox{Im}\  \phi^{n}_{\mathcal{M}_{x}})=(\mathcal{M}_{x})_{0}=(\mathcal{M}_{0})_{x}\ ,
$$
\noindent  et de l'application canonique
$$
(\mathcal{M}_{0})_{x}\rightarrow (\bigcap_{n}\mbox{Im}\  \phi^{n}_{\mathcal{M}})_{x}
$$
est l'identit\'e de $ (\bigcap_{n}\mbox{Im}\  \phi^{n}_{\mathcal{M}})_{x}$. D'o\`u le lemme. \ $\square$\\

\vskip 3mm
\noindent \textbf{Lemme (5.4.2.2)}.\textit {Avec les notations de 5.3, soit $(M,\phi_{M}) \in F^a\mbox{-}\textrm{Modlib}\ (A^{\dag})$ un $F$-module libre surconvergent ordinaire , i.e.  dont le $F$-module convergent associ\'e $(\mathcal{M},\phi_{\mathcal{M}})$ est ordinaire, donc muni d'une filtration satisfaisant \`a (5.3.2.2). Pour chaque entier $n\in \mathbb{N}$ on note}
 $$ 
 \phi^{n}_{\mathcal{M}} :  F^{n\ast}_{\hat{A}}(\mathcal{M}) \rightarrow \mathcal{M} \qquad [resp.\  \phi^{n}_{M} :  F^{n\ast}_{A^{\dag}}(M) \rightarrow M]
 $$
  \noindent\textit{l'it\'er\'e de $\phi_{\mathcal{M}}$ [resp.  $\phi_{M}$]. On suppose de plus qu'il existe un sous-$A^{\dag}$-module $M_{0}$ de $M$ tel que $\mathcal{M}_{0}=M_{0}\otimes_{A^{\dag}}\hat{A}$. Alors, on a des isomorphismes canoniques
  $$
  M_{0}=\mathcal{M}_{0}\cap M\ ,
  $$
$$
M_{0}= \bigcap_{n}\mbox{Im} \lbrace\phi^{n}_{M} :  F^{n\ast}_{A^{\dag}}(M) \rightarrow M\rbrace\ ,
$$
et $M_{0}$ est $\phi_{M}$-stable, donc muni d'un Frobenius
$$
 \phi_{M_{0}} :  F^{\ast}_{A^{\dag}}(M_{0}):=M_{0}^{\sigma} \rightarrow M_{0}
$$
d\'efini par $\phi_{M_{0}}=\phi_{M\vert M_{0}^{\sigma}}=\phi_{\mathcal{M}\vert M_{0}^{\sigma}}$}\ .
\vskip 3mm
\noindent \textit{Preuve du lemme.}  Notons 
$$M'_{0}=\mathcal{M}_{0}\cap M \ ;$$
 comme on a des inclusions \'evidentes
$$
M_{0}\subset M'_{0}\subset \mathcal{M}_{0} \ ,
$$
que le compl\'et\'e de $M_{0}$ est $\mathcal{M}_{0}$ et que $M'_{0}$ est un $A^{\dag}$-module de type fini, il en r\'esulte une \'egalit\'e
$$
M_{0}\otimes_{A^{\dag}}\hat{A}=M'_{0}\otimes_{A^{\dag}}\hat{A}\ ;
$$
d'o\`u, par fid\`ele platitude de $\hat{A}$ sur $A^{\dag}$, la premi\`ere \'egalit\'e du lemme
$$
M_{0}=M'_{0}=\mathcal{M}_{0}\cap M \\ .
$$
Notons 
$$
M^{\sigma}= F^{\ast}_{A^{\dag}}(M)\ , \qquad M^{\sigma}_{0}=F^{\ast}_{A^{\dag}}(M_{0})\ ,
$$

$$\mathcal{M}^{\sigma}= F^{\ast}_{\hat{A}}(\mathcal{M})\ , \qquad \mathcal{M}_{0}^{\sigma}= F^{\ast}_{\hat{A}}(\mathcal{M}_{0})\ ;$$

\noindent on a alors des isomorphismes
$$
\mathcal{M}^{\sigma}\simeq M^{\sigma}\otimes_{A^{\dag}}\hat{A}\simeq F^{\ast}_{A^{\dag}}(\mathcal{M})\ , \qquad \mathcal{M}_{0}^{\sigma}\simeq M_{0}^{\sigma}\otimes_{A^{\dag}}\hat{A}\simeq F^{\ast}_{A^{\dag}}(\mathcal{M}_{0})\ .
 $$
 Puisque $F_{A^{\dag}}$ est plat, on en d\'eduit des isomorphismes canoniques [Bour, AC I, \S2, \no6, prop. 6 et Rq. 1]
 $$
  \mathcal{M}_{0}^{\sigma}\cap M^{\sigma}\simeq F^{\ast}_{A^{\dag}}(\mathcal{M}_{0})\cap F^{\ast}_{A^{\dag}}(M)\simeq F^{\ast}_{A^{\dag}} (\mathcal{M}_{0}\cap M)= M^{\sigma}_{0}\subset \mathcal{M}^{\sigma} \ ;
 $$
 d'o\`u un morphisme de Frobenius 
$$
\phi_{M_{0}}: M_{0}^{\sigma}\rightarrow M_{0}
$$
d\'efini par $\phi_{M_{0}}=\phi_{\mathcal{M}_{0} \vert M_{0}^{\sigma}}= \phi_{M\vert M_{0}^{\sigma}}= \phi_{\mathcal{M}\vert M_{0}^{\sigma}}$, et  $\phi_{M_{0}}$ est un isomorphisme car $\phi_{\mathcal{M}_{0}}=\phi_{M_{0}}\otimes_{A^{\dag}}\hat{A}$ en est un.\\

Comme on a des inclusions \'evidentes
$$
M_{0}\subset \bigcap_{n}\mbox{Im} \lbrace\phi^{n}_{M} :  F^{n\ast}_{A^{\dag}}(M) \rightarrow M\rbrace\subset \bigcap_{n}\mbox{Im} \lbrace\phi^{n}_{\mathcal{M}} :  F^{n\ast}_{\hat{A}}(\mathcal{M}) \rightarrow \mathcal{M}\rbrace= \mathcal{M}_{0} \ ,
$$
que le compl\'et\'e de $M_{0}$ est $\mathcal{M}_{0}$ et que $\bigcap_{n}\mbox{Im} \lbrace\phi^{n}_{M} :  F^{n\ast}_{A^{\dag}}(M) \rightarrow M\rbrace$ est un $A^{\dag}$-module de type fini, il en r\'esulte une \'egalit\'e
$$
M_{0}\otimes_{A^{\dag}}\hat{A}=\left(\bigcap_{n}\mbox{Im} \lbrace\phi^{n}_{M} :  F^{n\ast}_{A^{\dag}}(M) \rightarrow M\rbrace\right)\otimes_{A^{\dag}}\hat{A}\ ;
$$
d'o\`u, par fid\`ele platitude de $\hat{A}$ sur $A^{\dag}$, la deuxi\`eme \'egalit\'e du lemme
$$
M_{0}=\bigcap_{n}\mbox{Im} \lbrace\phi^{n}_{M} :  F^{n\ast}_{A^{\dag}}(M) \rightarrow M\rbrace \\ .\quad \square
$$

Plus g\'en\'eralement que pour la seule partie de pente z\'ero, on a:

\vskip 3mm
\noindent \textbf{Proposition (5.4.2.3)}.\textit {
\begin{enumerate}
\item[(a)] Avec les notations de 5.3, soient $M$ un $A^{\dag}$-module de type fini, $\mathcal{M}=M\otimes_{A^{\dag}}\hat{A}$ et $\mathcal{N}$ un sous-$\hat{A}$-module de $\mathcal{M}$. Alors les propri\'et\'es suivantes sont \'equivalentes:
	\begin{enumerate}
	\item[(i)] il existe un sous-$A^{\dag}$-module $N$ de $M$ tel que $\mathcal{N}=N\otimes_{A^{\dag}}\hat{A}$\ .
	\item[(ii)] $\left(\mathcal{N}\cap M\right)\otimes_{A^{\dag}}\hat{A}\simeq \mathcal{N}$\ .
	\end{enumerate}
	Si ces propri\'et\'es (i) et (ii) sont satisfaites, alors $N=\mathcal{N}\cap M$\ .
\item[(b)] Soit $(M,\phi_{M}) \in F^a\mbox{-}\textrm{Mod}\ (A^{\dag})$ un $F$-module  surconvergent, $(\mathcal{M},\phi_{\mathcal{M}})$ le $F$-module convergent associ\'e $(\mathcal{M},\phi_{\mathcal{M}})$ et 
$$
(\mathcal{N},\phi_{\mathcal{N}})\hookrightarrow(\mathcal{M},\phi_{\mathcal{M}})
 $$ 
 un sous-$F$-module convergent.\\
  On suppose de plus qu'il existe un sous-$A^{\dag}$-module $N$ de $M$ tel que $\mathcal{N}=N\otimes_{A^{\dag}}\hat{A}$. Alors, on a un isomorphisme canonique
  $$
  N=\mathcal{N}\cap M\ ,
  $$
et $N$ est $\phi_{M}$-stable, donc muni d'un Frobenius
$$
 \phi_{N} :  F^{\ast}_{A^{\dag}}(N):=N^{\sigma} \rightarrow N
$$
d\'efini par $\phi_{N}=\phi_{\mathcal{N}\vert N^{\sigma}}=\phi_{M\vert N^{\sigma}} = \phi_{\mathcal{M}\vert N^{\sigma}}$\ . De plus $\phi_{\mathcal{N}}$ est un isomorphisme si et seulement si $\phi_{N}$ en est un.
\item[(c)] Situation comme en (a) en supposant (i) et (ii) v\'erifi\'ees: on suppose de plus que $M$, $\mathcal{M}$, $\mathcal{N}$ sont munis de connexions $\nabla_{M}, \nabla_{\mathcal{M}}, \nabla_{\mathcal{N}}$, telles que $\nabla_{M}=\left(\nabla_{\mathcal{M}}\right)_{\vert {M}}, \nabla_{\mathcal{N}}=\left(\nabla_{\mathcal{M}}\right)_{\vert \mathcal{N}}$. Alors $N$ est muni d'une connexion d\'efinie par
$$
\nabla_{N}=\left(\nabla_{\mathcal{N}}\right)_{\vert N}=\left(\nabla_{M}\right)_{\vert N}=\left(\nabla_{\mathcal{M}}\right)_{\vert N}\ .
$$
\item[(d)] Supposons (b) et (c) v\'erifi\'ees. Si $\phi_{M}$ et $\phi_{\mathcal{N}}$ sont horizontaux alors $$
 \phi_{N} : (N^{\sigma},\nabla^{\sigma}) \rightarrow (N,\nabla)
$$
est aussi un morphisme horizontal.
\item[(e)] Les propri\'et\'es (a),(b),(c),(d) subsistent en rempla\c cant $A^{\dag}$ par $A^{\dag}_{K}$ et $\hat{A}$ par $\hat{A}_{K}$.
\end{enumerate}
}
\vskip 3mm

\noindent \textit{Preuve de la proposition.} \\
(a) Supposons (i) et notons 
$$N'=\mathcal{N}\cap M \ ;$$
 comme on a des inclusions \'evidentes
$$
N\subset N'\subset \mathcal{N} \ ,
$$
que le compl\'et\'e de $N$ est $\mathcal{N}$ et que $N'$ est un $A^{\dag}$-module de type fini, il en r\'esulte une \'egalit\'e
$$
N\otimes_{A^{\dag}}\hat{A}=N'\otimes_{A^{\dag}}\hat{A}\ ,
$$
donc, par fid\`ele platitude de $\hat{A}$ sur $A^{\dag}$, l'\'egalit\'e $N=\mathcal{N}\cap M$. D'o\`u (ii).\\
L'implication r\'eciproque est claire.\\

(b) D'apr\`es (a) on a 
 $$
  N=\mathcal{N}\cap M\ .
  $$
Notons 
$$
M^{\sigma}= F^{\ast}_{A^{\dag}}(M)\ , \qquad N^{\sigma}=F^{\ast}_{A^{\dag}}(N)\ ,
$$

$$\mathcal{M}^{\sigma}= F^{\ast}_{\hat{A}}(\mathcal{M})\ , \qquad \mathcal{N}^{\sigma}= F^{\ast}_{\hat{A}}(\mathcal{N})\ ;$$

\noindent on a alors des isomorphismes
$$
\mathcal{M}^{\sigma}\simeq M^{\sigma}\otimes_{A^{\dag}}\hat{A}\simeq F^{\ast}_{A^{\dag}}(\mathcal{M})\ , \qquad \mathcal{N}^{\sigma}\simeq N^{\sigma}\otimes_{A^{\dag}}\hat{A}\simeq F^{\ast}_{A^{\dag}}(\mathcal{N})\ .
 $$
 Puisque $F_{A^{\dag}}$ est plat, on d\'eduit de [Bour, AC I, \S2, \no6, prop. 6 et Rq. 1] les  isomorphismes du lemme suivant:
 \vskip 3mm
\noindent \textbf{Lemme (5.4.2.4)}.
 $$
  \mathcal{N}^{\sigma}\cap M^{\sigma}\simeq F^{\ast}_{A^{\dag}}(\mathcal{N})\cap F^{\ast}_{A^{\dag}}(M)\simeq F^{\ast}_{A^{\dag}} (\mathcal{N}\cap M)= N^{\sigma}\subset \mathcal{M}^{\sigma} \ .
 $$
  \vskip 3mm
 D'o\`u, par intersection \`a partir de $\phi_{\mathcal{N}}$ et $\phi_{M}$, un morphisme de Frobenius 
$$
\phi_{N}: N^{\sigma}\rightarrow N
$$
d\'efini par $\phi_{N}=\phi_{\mathcal{N} \vert N^{\sigma}}= \phi_{M\vert N^{\sigma}}= \phi_{\mathcal{M}\vert N^{\sigma}}$.\\
 Par fid\`ele platitude de $\hat{A}$  sur $A^{\dag}$, $\phi_{N}$ est un isomorphisme si et seulement si $\phi_{\mathcal{N}}=\phi_{N}\otimes_{A^{\dag}}\hat{A}$ en est un.\\

(c) On d\'eduit de [Bour, AC I, \S2, \no6, prop. 6 et Rq. 1] le lemme suivant:
\vskip 3mm
\noindent \textbf{Lemme (5.4.2.5)}. \textit {Avec les notations pr\'ec\'edentes on a un isomorphisme}
$$
N\otimes_{A^{\dag}}\Omega^{1}_{A^{\dag}}\simeq (\mathcal{N}\otimes_{\hat{A}}\Omega^{1}_{\hat{A}})\cap(M\otimes_{A^{\dag}}\Omega^{1}_{A^{\dag}})\ .
$$
Gr\^ace \`a ce lemme (5.4.2.5) les connexions 
$$\nabla_{\mathcal{N}}: \mathcal{N}\rightarrow \mathcal{N}\otimes_{\hat{A}}\Omega^{1}_{\hat{A}}$$
 et
 $$
 \nabla_{M}: M\rightarrow M\otimes_{A^{\dag}}\Omega^{1}_{A^{\dag}}
 $$
 permettent de d\'efinir encore par intersection, \`a la mani\`ere de [Et 4, \S4], une connexion 
 $$
 \nabla_{N}: N\rightarrow N\otimes_{A^{\dag}}\Omega^{1}_{A^{\dag}}
 $$
  telle que 
  $$
  \nabla_{N}=\nabla_{\mathcal{N}\vert N}= \nabla_{M\vert N}=\nabla_{\mathcal{M}\vert N} \ .
  $$
  
  (d) Par l'injectivit\'e des fl\`eches
  
  $$
(N,\phi_{N})\hookrightarrow(M,\phi_{M})
 $$ 
  et
  $$
(N,\phi_{N})\hookrightarrow(\mathcal{N},\phi_{\mathcal{N}})
 $$ 
  la commutation de $\phi_{N}$ et $\nabla_{N}$ r\'esulte de celle de $\phi_{M}$, $\nabla_{M}$ d'une part et de celle de $\phi_{\mathcal{N}}, \nabla_{\mathcal{N}}$ d'autre part.\\
  
  (e) Moyennant le lemme suivant, les preuves sont analogues.\\
   \vskip 3mm
\noindent \textbf{Lemme (5.4.2.6)}. \textit{En munissant $A^{\dag}_{K}=A^{\dag}\otimes_{\mathcal{V}}K$ (resp. $\hat{A}_{K}=\hat{A}\otimes_{\mathcal{V}}K$) de la norme extension naturelle de la norme $p$-adique de $A^{\dag}$ (resp. $\hat{A}$), on note $\widehat{\left( A^{\dag}_{K} \right)}$ le s\'epar\'e compl\'et\'e de $A^{\dag}_{K}$. Soit $M$ un $A^{\dag}_{K}$-module de type fini: on le munit de la norme $p$-adique et on note $\hat{M}$ son s\'epar\'e compl\'et\'e pour cette topologie. Alors: 
  \begin{enumerate}
  \item[(i)] $\hat{A}_{K}$ est s\'epar\'e et complet.
  \item[(ii)]  On a un isomorphisme isom\'etrique $\widehat{\left( A^{\dag}_{K} \right)}\simeq\hat{A}_{K}$.
  \item[(iii)]  On a un isomorphisme isom\'etrique  $\hat{M}\simeq\hat{A}_{K}\otimes_{A^{\dag}_{K}}M$.
  \end{enumerate}
    }
       \vskip 3mm

\noindent \textit{Preuve du lemme.} \\

(i) Que $\hat{A}_{K}$ soit s\'epar\'e est clair. Soit $(a_{\nu})_{\nu}$ une suite de Cauchy dans $\hat{A}_{K}$: pour tout $\varepsilon>0$ il existe un entier $N_{0}$ tel que, pour tous entiers naturels $\nu, \mu, \nu\geqslant N_{0}$, on ait $\vert a_{\nu+\mu}-a_{\nu}\vert<\varepsilon$; en particulier pour $0<\varepsilon<1$ on a $\vert a_{\nu+\mu}-a_{\nu}\vert \in \hat{A}$. Choisissons un entier $N(\nu)$ tel que $p^{N(\nu)}a_{\nu}\in \hat{A}$; puisque
$$
\vert p^{N(\nu)}a_{\nu+\mu}-p^{N(\nu)}a_{\nu}\vert=\frac{1}{p^{N(\nu)}}\vert a_{\nu+\mu}-a_{\nu}\vert<\varepsilon
$$
on en d\'eduit que, pour tous entiers naturels $\nu, \mu, \nu\geqslant N_{0}$, on a $p^{N(\nu)}a_{\nu+\mu}\in \hat{A}$ et que, pour  $\nu\geqslant N_{0}$, $(p^{N(\nu)}a_{\nu+\mu})_{\mu}$ est une suite de Cauchy de $\hat{A}$; notons $b\in \hat{A}$ la limite de cette suite. Il est alors clair que 
$$\frac{b}{p^{N(\nu)}}=\lim_{n\to \infty}a_{n}\ .
$$

(ii) D'apr\`es [B-G-R, \S2.1.7, prop 4] et la compl\'etude de $\hat{A}_{K}$ on a les isomorphismes isom\'etriques suivants:
$$
\widehat{\left( A^{\dag}_{K} \right)}\simeq\widehat{A^{\dag}}\hat{\otimes}_{\hat{\mathcal{V}}}\hat{K}=\hat{A}\hat{\otimes}_{\mathcal{V}}K \ ,
$$
$$
\hat{A}_{K}=\widehat{\hat{A}_{K}}\simeq\hat{A}\hat{\otimes}_{\mathcal{V}}K\ ;
$$
  d'o\`u le (ii). \\
  
  (iii) Notons $M'=\hat{A}_{K}\otimes_{A^{\dag}_{k}}M$: c'est un $\hat{A}_{K}$-module de type fini que l'on peut d\'ecrire comme le quotient $  M'= M''/N  $ d'un module libre de type fini $M''$ par un sous-module $N$. D'apr\`es [B-G-R, \S2.1.1 prop 3 et \S2.1.5 prop 6] $M''$ est complet, donc $M'$ aussi [B-G-R, \S2.1.2 prop 3]. Ainsi, compte tenu du (ii), on a $M'\simeq\hat{M'}$, i.e [B-G-R, \S2.1.7 prop 4 et prop 6 (i)]
  $$
  M'\simeq\hat{A}_{K}\hat{\otimes}_{\hat{A}_{K}}\hat{M}\simeq\hat{M}\ .
  $$ 
  D'o\`u le lemme.    $\square$\\

  On va en d\'eduire le crit\`ere suivant de surconvergence pour les sous-objets des $F$-isocristaux, pour la d\'efinition desquels nous renvoyons \`a [B 2]:
  
 \vskip 3mm
\noindent \textbf{Corollaire (5.4.2.7)}.\textit{Soient X comme en (5.4.1), $(E,\phi_{E})\in F\mbox{-}Isoc^{\dag}(X/ K)$ un $F$-isocristal surconvergent  [B 2, (2.3.7)] et $(\mathcal{E},\phi_{\mathcal{E}})\in F\mbox{-}Isoc(X/ K)$ le $F$-isocristal convergent associ\'e: on note $(M,\phi_{M}, \nabla_{M})$ le $A^{\dag}_{K}$-module associ\'e \`a $E$ par l'\'equivalence de cat\'egories de Berthelot [B 2, (2.5.8)] et $(\mathcal{M},\phi_{\mathcal{M}}, \nabla_{\mathcal{M}})$ le $\hat{A}_{K}$-module associ\'e \`a $\mathcal{E}$ [Et 8]. Consid\'erons un sous-objet $$(\mathcal{E'},\phi_{\mathcal{E'}})\hookrightarrow(\mathcal{E},\phi_{\mathcal{E}})$$
de $(\mathcal{E},\phi_{\mathcal{E}})$ et notons 
$$(\mathcal{N},\phi_{\mathcal{N}}, \nabla_{\mathcal{N}})\hookrightarrow(\mathcal{M},\phi_{\mathcal{M}}, \nabla_{\mathcal{M}})$$
 le $\hat{A}_{K}$-module associ\'e \`a $\mathcal{E'}$ [Et 8]. Alors on a \'equivalence entre les propri\'et\'es suivantes:
\begin{enumerate}
\item[(i)] $(\mathcal{E'},\phi_{\mathcal{E'}})$ provient d'un objet $(E',\phi_{E'})\in F\mbox{-}Isoc^{\dag}(X/ K)$
par le foncteur d'oubli [B 2, (2.3.9)]
$$
F\mbox{-}Isoc^{\dag}(X/ K)\rightarrow F\mbox{-}Isoc(X/ K)\ .
$$
\item[(ii)] $\left(\mathcal{N}\cap M\right)\otimes_{A^{\dag}_{K}}\hat{A}_{K}\simeq \mathcal{N}$.
\end{enumerate}
}

\vskip 3mm
\noindent \textit{Preuve du corollaire.} Supposons (i): il existe $E'\in F\mbox{-}Isoc^{\dag}(X/ K)$ avec pour image $\mathcal{E'}$ par le foncteur d'oubli
$$
F\mbox{-}Isoc^{\dag}(X/ K)\rightarrow F\mbox{-}Isoc(X/K)\ .
$$
L'injection $\mathcal{E'}\hookrightarrow \mathcal{E}$ se rel\`eve alors par la pleine fid\'elit\'e du foncteur d'oubli \'etablie par Kedlaya [Ked 2, Theo 5.2.1] en une injection $E'\hookrightarrow E$: par l'\'equivalence de cat\'egories de Berthelot [B 2, (2.5.8)] il en r\'esulte une injection entre $A^{\dag}_{K}$- modules projectifs de type fini $N\hookrightarrow M$ telle que
$$
\mathcal{N}\simeq N\otimes_{A^{\dag}_{K}}\hat{A}_{K}.
$$
Le (ii) r\'esulte alors de la proposition (5.4.2.3)(e; propri\'et\'e (a)).\\

R\'eciproquement supposons (ii) et posons $N=\mathcal{N}\cap M$. D'apr\`es [(5.4.2.3)(a) et (e)] $N$ est muni d'une connexion $\nabla_{N}$ et d'un Frobenius $\phi_{N}$ qui est un isomorphisme horizontal: d'apr\`es l'\'equivalence de cat\'egories de Berthelot [B 2, (2.5.8)] il existe un objet $(E',\phi_{E'})\in F\mbox{-}Isoc^{\dag}(X/ K)$ associ\'e \`a $(N,\phi_{N}, \nabla_{N})$. Compte tenu de l'isomorphisme
$$
(N,\phi_{N}, \nabla_{N})\otimes_{A^{\dag}_{K}}\hat{A}_{K}\simeq(\mathcal{N},\phi_{\mathcal{N}}, \nabla_{\mathcal{N}})\ ,
$$
$(\mathcal{E'},\phi_{\mathcal{E'}})$ provient de $(E',\phi_{E'})\in F\mbox{-}Isoc^{\dag}(X/ K)$
par le foncteur d'oubli. D'o\`u le corollaire. $\square$

\subsection*{5.4.3 Premier contre-exemple: les courbes elliptiques de Legendre}

Soient $k=\mathbb{F}_{p}$ le corps fini \`a $p$ \'el\'ements, $\overline {k}$ une cl\^oture alg\'ebrique de $k$, $\lambda\in \overline {k},\ \lambda\neq 0, 1, k(\lambda)= \mathbb{F}_{q_{\lambda}}$ le corps r\'esiduel de $\lambda$, et $X_{\lambda}$ la courbe elliptique d'\'equation affine $y^{2}=x(x-1)(x-\lambda)\ .$\\

 Si $p\neq 2$, soit $H_{p}$ le polyn\^ome suivant
$$
H_{p}(\lambda)=\sum_{i=0}^{i=\frac{p-1}{2}}(_{i}^{\frac{p-1}{2}})^{2}\ \lambda^{i} \ .
$$ 
Alors par d\'efinition on a
\begin{itemize}
\item[(i)]  $X_{\lambda}$ est  supersinguli\`ere $\iff\ H_{p}(\lambda)=0$
\item[(ii)]  $X_{\lambda}$ est  ordinaire $\iff\ H_{p}(\lambda)\neq0$.
\end{itemize}
Sur cette caract\'erisation on constate que les courbes elliptiques ordinaires sont les plus nombreuses; pour $p>2$ donn\'e il y a au plus [p/12]+2 courbes elliptiques supersinguli\`eres (\`a isomorphisme pr\`es).  \\

La fonction z\^eta de $X_{\lambda}$ s'\'ecrit:
$$
Z(X_{\lambda}/k(\lambda), t)=\frac{(1-\alpha_{\lambda}t)(1-\beta_{\lambda}t)}{(1-t)(1-qt)}=\frac{det(1-tF;H^{1}_{cris}(X_{\lambda}/W_{\lambda}))}{(1-t)(1-qt)}\ ,
$$
\noindent o\`u $W_{\lambda}$ est l'anneau $W(k(\lambda))$ des vecteurs de Witt de $k(\lambda)$ et  $\alpha_{\lambda}$ , $\beta_{\lambda}$ sont des entiers alg\'ebriques \'el\'ements d'une extension finie $K_{\lambda}$ du corps des fractions de $W_{\lambda}$ tels que $\beta_{\lambda}=q_{\lambda}/\alpha_{\lambda}$. Soit $v$ l'extension \`a $K_{\lambda}$ de la valuation $p$-adique telle que $v(q_{\lambda})=1$. Il n'y a alors que deux possibilit\'es:
\begin{itemize}
\item[(i)] ou bien $v(\alpha_{\lambda})=v(\beta_{\lambda})=1/2$, auquel cas $X_{\lambda}$ est  supersinguli\`ere,
\item[(ii)] ou bien $v(\alpha_{\lambda})=0, v(\beta_{\lambda})=1$, auquel cas $X_{\lambda}$ est  ordinaire. \\ 
\end{itemize}

Dans le cas ordinaire (qui est "le plus courant") on s'int\'eresse \`a voir comment $\alpha_{\lambda}$ varie avec $\lambda$: pour \c ca on consid\`ere la famille des courbes elliptiques de Legendre $f: X\rightarrow S$ sur le $\mathbb{F}_{p}$-sch\'ema lisse $S=Spec(\mathbb{F}_{p}[\lambda][\frac{1}{\lambda(1-\lambda)H_{p}(\lambda)}])$ et les fonctions g\'en\'eratrices

$$L(R^1 f_{\textrm{cris}\ast}(\mathcal{O}_{X/W}), t) :=\prod_{s\in\vert S\vert}[(1-\beta_{s}t^{deg\ s})(1-\alpha_{s}t^{deg\ s})]^{-{1}} \ ,
$$

$$L_{0}(R^1 f_{\textrm{cris}\ast}(\mathcal{O}_{X/W}), t) = L(R^1 f_{\textrm{\'et}\ast} (\mathbb		{Z}_{p}) \otimes_{\mathbb{Z}_{p}} \mathcal{O}_{X/W}, t)=\prod_{s\in\vert S\vert}(1-\alpha_{s}t^{deg\ s})^{-{1}}
$$
\noindent dans lesquelles 
$$
deg\ s=[k(s):\mathbb{F}_{p} ] \  .
$$

\noindent Par [B-B-M] on sait que $R^1 f_{\textrm{cris}\ast}(\mathcal{O}_{X/W})$ est un $F$-cristal localement libre de rang 2 et $R^1 f_{\textrm{\'et}\ast} (\mathbb{Z}_{p}) \otimes_{\mathbb{Z}_{p}} \mathcal{O}_{X/W}$ est son sous-$F$-cristal unit\'e. Le $F$-isocristal convergent associ\'e \`a $R^1 f_{\textrm{cris}\ast}(\mathcal{O}_{X/W})$ par la construction de Berthelot [B 2, (2.4.2)] n'est autre que $\mathcal{E}:=R^{1}f_{conv\ast}(\mathcal{O}_{X/\mathbb{Q}_{p}})$ et ce dernier provient en fait du $F$-isocristal surconvergent $E:=R^{1}f_{rig\ast}(\mathcal{O}_{X/\mathbb{Q}_{p}})$ [Et 4, th\'eor\`eme 7]: il r\'esulte alors de [E-LS 1] et de [Ked 1] que la fonction $L(R^1 f_{\textrm{cris}\ast}(\mathcal{O}_{X/W}), t)$ est rationnelle.\\
 D'autre part on sait depuis Dwork [Dw 3] que la fonction $ L_{0}(R^1 f_{\textrm{cris}\ast}(\mathcal{O}_{X/W}), t)$ est m\'eromorphe.\\

Notons $A_{0}=\mathbb{F}_{p}[\lambda][\frac{1}{\lambda(1-\lambda)H_{p}(\lambda)}], A=\mathbb{Z}_{p}[\lambda][\frac{1}{\lambda(1-\lambda)H_{p}(\lambda)}], A^{\dag}$ le compl\'et\'e faible de $A$, $\hat{A}$ le compl\'et\'e de $A$, $A^{\dag}_{\mathbb{Q}_{p}}=A^{\dag}\otimes_{\mathbb{Z}_{p}}\mathbb{Q}_{p}$, $\hat{A}_{\mathbb{Q}_{p}}=\hat{A}\otimes_{\mathbb{Z}_{p}}\mathbb{Q}_{p}$, $F_{A^{\dag}}:A^{\dag}\rightarrow A^{\dag}$ un rel\`evement du Frobenius de $A_{0}$ au-dessus du Frobenius $\sigma$ de $\mathbb{Z}_{p}$, $F_{\hat{A}}= F_{A^{\dag}}\otimes_{A^{\dag}}\hat{A}$ et $F_{A^{\dag}_{\mathbb{Q}_{p}}}=F_{A^{\dag}}\otimes_{\mathbb{Z}_{p}}\mathbb{Q}_{p}$, $F_{\hat{A}_{\mathbb{Q}_{p}}}=F_{\hat{A}}\otimes_{\mathbb{Z}_{p}}\mathbb{Q}_{p}$.\\
Rappelons que $F_{A^{\dag}}$ est automatiquement fini et fid\`element plat [Et 3, th\'eo 17], et que la donn\'ee d'un rel\`evement $F_{\hat{A}}: \hat{A}\rightarrow\hat{A}$ du Frobenius de $A_{0}$ au-dessus du Frobenius $\sigma$ de $\mathbb{Z}_{p}$ est \'equivalente [Et 6, cor. (3.1.4)] \`a la donn\'ee d'un $F_{A^{\dag}}$ tel que $F_{\hat{A}}= F_{A^{\dag}}\otimes_{A^{\dag}}\hat{A}$.\\ 
 Puisque $S$ est un $\mathbb{F}_{p}$-sch\'ema affine et lisse, la donn\'ee du $F$-isocristal surconvergent $E=R^{1}f_{rig\ast}(\mathcal{O}_{X/\mathbb{Q}_{p}})$ est \'equivalente [B 2, (2.5.8)] \`a la donn\'ee d'un $A^{\dag}_{\mathbb{Q}_{p}}$-module projectif de type fini $M_{\mathbb{Q}_{p}}$, muni d'une connexion int\'egrable $\nabla_{M_{\mathbb{Q}_{p}}}$ et d'un isomorphisme horizontal $\phi_{M_{\mathbb{Q}_{p}}}: F^{\ast}_{A^{\dag }_{\mathbb{Q}_{p}}}(M_{\mathbb{Q}_{p}})\rightarrow M_{\mathbb{Q}_{p}}$. Soit $(\mathcal{M}_{\mathbb{Q}_{p}}, \phi_{\mathcal{M}_{\mathbb{Q}_{p}}}, \nabla_{\mathcal{M}_{\mathbb{Q}_{p}}})$ le $\hat{A}_{\mathbb{Q}_{p}}$-module d\'eduit de $(M_{\mathbb{Q}_{p}}, \phi_{M_{\mathbb{Q}_{p}}}, \nabla_{M_{\mathbb{Q}_{p}}})$ par le changement de base $ A^{\dag}_{\mathbb{Q}_{p}}\rightarrow \hat{A}_{\mathbb{Q}_{p}}$: $\mathcal{M}_{\mathbb{Q}_{p}}$ n'est autre que le $\hat{A}_{\mathbb{Q}_{p}}$-module associ\'e dans [Et 8] au $F$-isocristal convergent correspondant \`a $R^1 f_{\textrm{cris}\ast}(\mathcal{O}_{X/W})$.\\

 Il se trouve aussi que la donn\'ee de $R^1 f_{\textrm{cris}\ast}(\mathcal{O}_{X/W})$ fournit d'apr\`es Katz [K 2, \S\ II] un $\hat{A}$-module projectif de type fini $\mathcal{M}$ muni d'une connexion int\'egrable $\nabla_{\mathcal{M}}$ et d'une isog\'enie $\phi_{\mathcal{M}}: F^{\ast}_{\hat{A}}(\mathcal{M})\rightarrow\mathcal{M}$ commutant \`a la connexion, et l'on a un isomorphisme canonique
 $$\mathcal{M}_{\mathbb{Q}_{p}}= \mathcal{M}\otimes_{\hat{A}}\hat{A}_{\mathbb{Q}_{p}}\ .$$
Il est \`a noter que le Frobenius $\phi_{\mathcal{M}}$ d\'epend du rel\`evement $F_{\hat{A}}$ du Frobenius de $A_{0}$: pour deux tels rel\`evements $F_{1\hat{A}}$ et $F_{2\hat{A}}$ les Frobenius $\phi_{1\mathcal{M}}$, $\phi_{2\mathcal{M}}$ correspondants sont reli\'es par un isomorphisme $\chi(F_{1},F_{2})$ provenant de la connexion et qui rend commutatif le diagramme suivant [K 1, \S(1.3)] :
$$
\xymatrix{
F_{1\hat{A}}^{\ast}(\mathcal{M})\ar[rr]^{\phi_{1\mathcal{M}}}\ar[d]^{\simeq}_{\chi(F_{1},F_{2})}&&\mathcal{M}\\
F_{2\hat{A}}^{\ast}(\mathcal{M})\ar[urr]_{\phi_{2\mathcal{M}}}&&\ .
}
$$
\vskip 3mm

 Pr\'ecisons la provenance de $\mathcal{M}$: soit $h: \mathcal{X}\rightarrow \mathcal{S}= Spec\ A$ le rel\`evement \'evident de $f$ au-dessus de $\mathbb{Z}_{p}$ et $h^{\dag}:\mathcal{X}^{\dag}\rightarrow\mathcal{S}^{\dag}:= Spec\ A^{\dag}$ [resp. $\hat{h}:\hat{\mathcal{X}}\rightarrow\hat{\mathcal{S}}:= Spec\ \hat{A}$]  l'image inverse de $h$ par le changement de base de $A$ \`a $A^{\dag} $ [resp. de $A$ \`a $\hat {A}$]; alors $\mathcal{M}$ est le $H^1 de\ Rham$ de $\hat{h}$ [K 1, \S\ 8] et $\mathcal{M}$ est un $\hat{A}$-module libre de rang 2 sur $\omega$ et $\omega'$ o\`u
\begin{enumerate}
\item[] $\omega$\  est la classe de la diff\'erentielle de premi\`ere esp\`ece\  $dx/y$\ ,

  \item[]$\omega'= \nabla(d/d\lambda)(\omega)$\ .
\end{enumerate}
La filtration de Hodge est sp\'ecifi\'ee par
$$
Fil^{1}\mathcal{M}= \hat{A}\omega\subset \mathcal{M}\ .
$$ 
La donn\'ee de $\mathcal{M}$ et $M_{\mathbb{Q}_{p}}$ fournit d'apr\`es [F-R, \S4] un $A^{\dag}$-module projectif de type fini $M=M_{\mathbb{Q}_{p}}\cap\mathcal{M}$ et l'on a des isomorphismes canoniques
$$
\mathcal{M}=M\otimes_{A^{\dag}}\hat{A}\ ,\qquad M_{\mathbb{Q}_{p}}= M \otimes_{A^{\dag}}A^{\dag}_{\mathbb{Q}_{p}} \ .
$$
Notons 
$$
M^{\sigma}= F^{\ast}_{A^{\dag}}(M)\ , \qquad M_{\mathbb{Q}_{p}}^{\sigma}= F^{\ast}_{A_{\mathbb{Q}_{p}}^{\dag}}(M_{\mathbb{Q}_{p}}) \ ,
$$
$$
\mathcal{M}^{\sigma}= F^{\ast}_{\hat{A}}(\mathcal{M})\ , \qquad \mathcal{M}_{\mathbb{Q}_{p}}^{\sigma}= F^{\ast}_{\hat{A}_{\mathbb{Q}_{p}}}(\mathcal{M}_{\mathbb{Q}_{p}}) \ ;
$$
on a alors des isomorphismes
$$
\mathcal{M}^{\sigma}\simeq M^{\sigma}\otimes_{A^{\dag}}\hat{A}\  ,\qquad (M_{\mathbb{Q}_{p}})^{\sigma}\simeq M^{\sigma}\otimes_{A^{\dag}} A^{\dag}_{\mathbb{Q}_{p}}=(M^{\sigma})_{\mathbb{Q}_{p}}\ ,
$$
$$
(\mathcal{M}_{\mathbb{Q}_{p}})^{\sigma}\simeq (\mathcal{M}^{\sigma})_{\mathbb{Q}_{p}}\simeq (M_{\mathbb{Q}_{p}})^{\sigma}\otimes_{A^{\dag}_{\mathbb{Q}_{p}}}\hat{A}_{\mathbb{Q}_{p}}\simeq M^{\sigma}\otimes_{A^{\dag}}\hat{A}_{\mathbb{Q}_{p}}\ .
$$
\vskip 3mm
\noindent \textbf{Lemme (5.4.3.1)}. \textit {Avec les notations pr\'ec\'edentes on a un isomorphisme}
$$
M^{\sigma}\simeq \mathcal{M}^{\sigma}\cap M^{\sigma}_{\mathbb{Q}_{p}}\subset \mathcal{M}^{\sigma}_{\mathbb{Q}_{p}} \ .
$$
\vskip 3mm
\noindent \textit{Preuve du lemme.} Le $A^{\dag}$-module $M^{\sigma}$ est plat pour le $A^{\dag}$-module $\hat{A}_{\mathbb{Q}_{p}}$ [Bour, AC I, \S2, \no2, d\'ef 1], car $M^{\sigma}$ est un $A^{\dag}$-module plat, puisque $F_{A^{\dag}}$ est plat et $M$ est projectif de type fini sur $A^{\dag}$. Donc on a un isomorphisme canonique [Bour, AC I, \S2, \no6, prop. 6 et Rq 1] et [Et 4, cor de prop.2]
$$
\mathcal{M}^{\sigma}\cap (M_{\mathbb{Q}_{p}})^{\sigma}=(M^{\sigma}\otimes_{A^{\dag}}\hat{A})\cap(M^{\sigma}\otimes_{A^{\dag}}A^{\dag}_{\mathbb{Q}_{p}})\simeq M^{\sigma}\otimes_{A^{\dag}}(\hat{A}\cap A^{\dag}_{\mathbb{Q}_{p}})=M^{\sigma}\ .\ \square
$$\\
Les Frobenius $\phi_{\mathcal{M}}: F^{\ast}_{\hat{A}}(\mathcal{M})\rightarrow\mathcal{M}$ et $\phi_{M_{\mathbb{Q}_{p}}}: F^{\ast}_{A^{\dag }_{\mathbb{Q}_{p}}}(M_{\mathbb{Q}_{p}})\rightarrow M_{\mathbb{Q}_{p}}$ permettent alors, gr\^ace au lemme (5.4.3.1), de d\'efinir par intersection, \`a la mani\`ere de [Et 4, cor de prop. 7], un morphisme de Frobenius 
$$
\phi_{M}: M^{\sigma}\rightarrow M
$$
tel que $\phi_{M}=\phi_{\mathcal{M} \vert M^{\sigma}}= \phi_{M_{\mathbb{Q}_{p}}\vert M^{\sigma}}= \phi_{\mathcal{M}_{\mathbb{Q}_{p}}\vert M^{\sigma}}$ .\\

De la m\^eme fa\c con on \'etablit le lemme suivant:
\vskip 3mm
\noindent \textbf{Lemme (5.4.3.2)}. \textit {Avec les notations pr\'ec\'edentes on a un isomorphisme}
$$
M\otimes_{A^{\dag}}\Omega^{1}_{A^{\dag}}\simeq (\mathcal{M}\otimes_{\hat{A}}\Omega^{1}_{\hat{A}})\cap(M_{\mathbb{Q}_{p}}\otimes_{A^{\dag}_{\mathbb{Q}_{p}}}\Omega^{1}_{A^{\dag}_{\mathbb{Q}_{p}}})\ .
$$
Gr\^ace \`a ce lemme (5.4.3.2) les connexions 
$$\nabla_{\mathcal{M}}: \mathcal{M}\rightarrow \mathcal{M}\otimes_{\hat{A}}\Omega^{1}_{\hat{A}}$$
 et
 $$
 \nabla_{M_{\mathbb{Q}_{p}}}: M_{\mathbb{Q}_{p}}\rightarrow M_{\mathbb{Q}_{p}}\otimes_{A^{\dag}_{\mathbb{Q}_{p}}}\Omega^{1}_{A^{\dag}_{\mathbb{Q}_{p}}}
 $$
 permettent de d\'efinir encore par intersection, \`a la mani\`ere de [Et 4, \S4], une connexion 
 $$
 \nabla_{M}: M\rightarrow M\otimes_{A^{\dag}}\Omega^{1}_{A^{\dag}}
 $$
  telle que 
  $$
  \nabla_{M}=\nabla_{\mathcal{M}\vert M}= \nabla_{M_{\mathbb{Q}_{p}}\vert M}=\nabla_{\mathcal{M}_{\mathbb{Q}_{p}}\vert M} \ .
  $$
  On a vu ci-dessus que $\nabla_{M_{\mathbb{Q}_{p}}}$ et $\phi_{M_{\mathbb{Q}_{p}}}$ commutent, il en donc de m\^eme pour $ \nabla_{M}$ et $\phi_{M}$.\\

Le sous-$F$-cristal unit\'e $R^1 f_{\textrm{\'et}\ast} (\mathbb{Z}_{p}) \otimes_{\mathbb{Z}_{p}} \mathcal{O}_{X/W}$ de $R^1 f_{\textrm{cris}\ast}(\mathcal{O}_{X/W})$ fournit de m\^eme par Katz [K 2, \S II] un $\hat{A}$-module projectif de type fini $\mathcal{M}_{0}$ muni d'une connexion int\'egrable $\nabla_{\mathcal{M}_{0}}$ et d'une isog\'enie $\phi_{\mathcal{M}_{0}}: F^{\ast}_{\hat{A}}(\mathcal{M}_{0})\rightarrow\mathcal{M}_{0}$ commutant \`a la connexion; $\mathcal{M}_{0}$ est aussi le sous-module de $\mathcal{M}$ d\'efini dans le lemme (5.4.2.1) ci-dessus . En fait $\mathcal{M}_{0}$ est libre de rang un sur $\hat{A}$ [K 1, \S8], [vdP, (7.12)] avec pour base $u= \beta\omega-\lambda(1-\lambda)\omega' \  ,\beta\in\hat{A}$, telle que $\phi_{\mathcal{M}_{0}}(u)=\phi_{\mathcal{M}}(u)=u$. La filtration
$$
0\subset \mathcal{M}_{0}\subset \mathcal{M}
\leqno(5.4.3.3)
$$
en sous-modules $\phi_{\mathcal{M}}$-stables fait de $(\mathcal{M}, \phi_{\mathcal{M}})$ un $F$-module ordinaire [K 1] et l'on a un isomorphisme de $\hat{A}$-modules
$$
\mathcal{M}\simeq \mathcal{M}_{0}\oplus Fil^{1}\ \mathcal{M} \ .
\leqno(5.4.3.4)
$$
En proc\'edant \`a une \'etude locale de la famille de Legendre [De$\ell$ 2, \S2.1.4], [K 1, \S7] il se trouve que $\mathcal{M}_{0}$ et $Fil^{1}\ \mathcal{M}$ sont en dualit\'e par la dualit\'e de Poincar\'e qui respecte l'action du Frobenius [loc. cit.]: pour les \guillemotleft rel\`evements excellents\guillemotright du Frobenius qui induisent une action sur $Fil^{1}\ \mathcal{M}$ [cf (5.2)], la dualit\'e de Poincar\'e permet d'en d\'eduire l'action du Frobenius sur la partie unit\'e $\mathcal{M}_{0}$.

Nous sommes \`a pr\'esent en mesure d'\'enoncer notre premier contre-exemple promis en 5.1:
\vskip 3mm
\noindent \textbf{Th\'eor\`eme (5.4.3.5)}. \textit {Avec les notations pr\'ec\'edentes la filtration (5.4.3.3) ne se rel\`eve pas en une filtration de $M$, plus pr\'ecis\'ement il n'existe pas de sous-$A^{\dag}$-module $M_{0}$ de $M$ tel que l'on ait un isomorphisme canonique}
$$
\mathcal{M}_{0}=M_{0}\otimes_{A^{\dag}}\hat{A}\ .
$$
\vskip 3mm
\noindent \textit{Preuve du th\'eor\`eme.} Par l'absurde supposons l'existence d'un sous-$A^{\dag}$-module $M_{0}$ de $M$ v\'erifiant un isomorphisme canonique
$$
\mathcal{M}_{0}=M_{0}\otimes_{A^{\dag}}\hat{A}\ .
$$
 Notons 
$$M'_{0}=\mathcal{M}_{0}\cap M \ ;$$
 comme on a les inclusions \'evidentes
$$
M_{0}\subset M'_{0}\subset \mathcal{M}_{0} \ ,
$$
que le compl\'et\'e de $M_{0}$ est $\mathcal{M}_{0}$ et que $M'_{0}$ est un $A^{\dag}$-module de type fini, il en r\'esulte une \'egalit\'e
$$
M_{0}\otimes_{A^{\dag}}\hat{A}=M'_{0}\otimes_{A^{\dag}}\hat{A}\ ,
$$
d'o\`u l'\'egalit\'e
$$
M_{0}=M'_{0}=\mathcal{M}_{0}\cap M \\ .
$$
Comme pour le lemme (5.4.3.1) on \'etablit l'isomorphisme
$$
M_{0}^{\sigma}\simeq \mathcal{M}_{0}^{\sigma}\cap M^{\sigma}\subset \mathcal{M}^{\sigma}\ ;
$$
d'o\`u un morphisme de Frobenius 
$$
\phi_{M_{0}}: M_{0}^{\sigma}\rightarrow M_{0}
$$
d\'efini par $\phi_{M_{0}}=\phi_{\mathcal{M}_{0} \vert M_{0}^{\sigma}}= \phi_{M\vert M_{0}^{\sigma}}= \phi_{\mathcal{M}\vert M_{0}^{\sigma}}$ .\\

\noindent Par le lemme (5.4.2.5) on a un isomorphisme
$$
M_{0}\otimes_{A^{\dag}}\Omega^{1}_{A^{\dag}}\simeq (\mathcal{M}_{0}\otimes_{\hat{A}}\Omega^{1}_{\hat{A}})\cap(M\otimes_{A^{\dag}}\Omega^{1}_{A^{\dag}})\ ;
$$
d'o\`u une connexion
 $$
 \nabla_{M_{0}}: M_{0}\rightarrow M_{0}\otimes_{A^{\dag}}\Omega^{1}_{A^{\dag}}
 $$
  d\'efinie par 
  $$
  \nabla_{M_{0}}=\nabla_{\mathcal{M}_{0}\vert M_{0}}= \nabla_{M\vert M_{0}}=\nabla_{\mathcal{M}\vert M_{0}} \ .
  $$
Comme $\nabla_{M}$ et $\phi_{M}$ commutent, il en est de m\^eme de $\nabla_{M_{0}}$ et $\phi_{M_{0}}$ et de $\nabla_{M_{0\mathbb{Q}_{p}}}$ et $\phi_{M_{0\mathbb{Q}_{p}}}$: par l'\'equivalence de cat\'egories de Berthelot [B 2, (2.5.8)] le sous-$F$-isocristal convergent $\mathcal{E}_{0}$ du $F$-isocristal surconvergent $R^{1}f_{rig\ast}(\mathcal{O}_{X/\mathbb{Q}_{p}})$ correspondant \`a $R^1 f_{\textrm{\'et}\ast} (\mathbb{Z}_{p}) \otimes_{\mathbb{Z}_{p}} \mathcal{O}_{X/W}$ serait surconvergent, i.e. serait dans l'image essentielle du foncteur d'oubli [B 2, (2.3.8)]
$$
F\mbox{-}Isoc^{\dag}(S/ \mathbb{Q}_{p})\rightarrow F\mbox{-}Isoc(S/ \mathbb{Q}_{p})\ .
$$
Nous allons prouver qu'il n'en est rien.\\
Compte tenu de l'isomorphisme $\mathcal{M}_{0}=M_{0}\otimes_{A^{\dag}}\hat{A}$,   il existe \\
\begin{enumerate}
\item[]$\mu\in \hat{A}^{\times}$:= \{ \'el\'ements inversibles de $\hat{A}$\} tel que \\
\item[]$\mu u= \mu\beta\omega-\mu\lambda(1-\lambda)\omega' \in M_{0}$; \\
\end{enumerate}
puisque $\lambda(1-\lambda)\in A^{\dag\times}$:= \{ \'el\'ements inversibles de $A^{\dag}$\}, on en d\'eduit que $\mu\in A^{\dag\times}$, donc aussi $\beta\in A^{\dag\times}$. Au final ceci prouverait que
le g\'en\'erateur $u$ de $\mathcal{M}_{0}$ comme $\hat{A}$-module serait aussi g\'en\'erateur de $M_{0}$ comme $A^{\dag}$-module: on aurait alors (et ceci quel que soit le rel\`evement $F_{\hat{A}}$ du Frobenius de $A_{0}$)
\begin{enumerate}
\item[](5.4.3.6)\quad $\phi_{M_{0}}(u)=\phi_{\mathcal{M}}(u)\ \in M_{0}\subset \mathcal{M}_{0}$, \\

\item[](5.4.3.7)\quad $\nabla_{M_{0}}(d/d\lambda)(u)=\nabla_{\mathcal{M}}(d/\lambda)(u) \ \in M_{0}\subset \mathcal{M}_{0}$\ .
\end{enumerate}
Arriv\'e \`a ce stade du raisonnement nous allons \`a pr\'esent choisir pour rel\`evement $F_{\hat{A}}$ le rel\`evement canonique $\varphi_{can}$ (rel\`evement excellent du Frobenius dans la terminologie de Dwork) et nous noterons
$$
\phi_{\mathcal{M}}(\varphi_{can}): \varphi^{\ast}_{can}(\mathcal{M})\rightarrow \mathcal{M}
$$
le Frobenius de $\mathcal{M}$ pour insister sur sa d\'ependance en $\varphi_{can}$. D'apr\`es [vdP, (7.14), (7.16)], [Dw 4] on a 

$$\phi_{\mathcal{M}}(\varphi_{can})(u)=\xi u \ , \quad
\xi=(-1)^{\frac{p-1}{2}}\frac{\alpha(\lambda)}{\alpha(\varphi_{can}(\lambda))}\in A^{\dag}
\leqno(5.4.3.8)
$$
o\`u $\alpha$ est la fonction hyperg\'eom\'etrique
$$
\alpha(\lambda):= F(\frac{1}{2}, \frac{1}{2}, 1,\lambda):=\sum_{i=0}^{\infty}\left(\frac{(\frac{1}{2})_{_{i}}}{i!}\right)^{2}\lambda^{i}\ , \qquad \left(\frac{1}{2}\right)_{i}:=\prod_{j=0}^{j=i-1}\left(\frac{1}{2}+j\right) \ ,
$$
 $$\nabla_{\mathcal{M}}(d/\lambda)(u)=\eta u \ , \quad
\eta= - \frac{\alpha'(\lambda)}{\alpha(\lambda)} \in \hat{A}\setminus A^{\dag} \ .
\leqno(5.4.3.9)
$$

\noindent La relation (5.4.3.9) contredit (5.4.3.7). Ceci ach\`eve la preuve du th\'eor\`eme (5.4.3.5). $\square$\\

Chemin faisant nous avons en partie prouv\'e le corollaire suivant:

\vskip 3mm
\noindent \textbf{Corollaire (5.4.3.10)}. \textit {Avec les notations de (5.4.3) le $F$-isocristal convergent associ\'e \`a $R^1 f_{\textrm{cris}\ast}(\mathcal{O}_{X/W})$ par la construction de Berthelot [B 2, (2.4.2)] n'est autre que $\mathcal{E}:=R^{1}f_{conv\ast}(\mathcal{O}_{X/\mathbb{Q}_{p}})$ et ce dernier provient du $F$-isocristal surconvergent $E:=R^{1}f_{rig\ast}(\mathcal{O}_{X/\mathbb{Q}_{p}})$ par le foncteur d'oubli [B 2, (2.3.8)]}
$$
F\mbox{-}Isoc^{\dag}(S/ \mathbb{Q}_{p})\rightarrow F\mbox{-}Isoc(S/ \mathbb{Q}_{p})\ .
$$

\textit{Le sous-$F$-isocristal convergent unit\'e $\mathcal{E}_{0}$ de $\mathcal{E}$ associ\'e \`a \\ 
$R^1 f_{\textrm{\'et}\ast} (\mathbb{Z}_{p}) \otimes_{\mathbb{Z}_{p}} \mathcal{O}_{X/W}$ n'est pas surconvergent, i.e. il n'est pas dans l'image essentielle du foncteur d'oubli }
$$
F\mbox{-}Isoc^{\dag}(S/ \mathbb{Q}_{p})\rightarrow F\mbox{-}Isoc(S/ \mathbb{Q}_{p})
$$
\textit{bien que $\mathcal{E}$ lui-m\^eme soit dans l'image de ce foncteur d'oubli.}

\vskip 3mm
\noindent \textit{Preuve du corollaire.} Supposons par l'absurde qu'il existe $E_{0}\in F\mbox{-}Isoc^{\dag}(S/ \mathbb{Q}_{p})$ avec pour image $\mathcal{E}_{0}$ par le foncteur d'oubli
$$
F\mbox{-}Isoc^{\dag}(S/ \mathbb{Q}_{p})\rightarrow F\mbox{-}Isoc(S/ \mathbb{Q}_{p})\ .
$$
L'injection $\mathcal{E}_{0}\hookrightarrow \mathcal{E}$ se rel\`eve alors par la pleine fid\'elit\'e du foncteur d'oubli \'etablie par Kedlaya [Ked 2, Theo 5.2.1] en une injection $E_{0}\hookrightarrow E$: par l'\'equivalence de cat\'egories de Berthelot [B 2, (2.5.8)] il en r\'esulte une injection entre $A^{\dag}_{\mathbb{Q}_{p}}$- modules projectifs de type fini $M_{0\mathbb{Q}_{p}}\hookrightarrow M_{\mathbb{Q}_{p}}$; avec les notations utilis\'ees dans la preuve du th\'eor\`eme (5.4.3.5) posons
$$
M_{0}:= \mathcal{M}_{0}\cap M_{0\mathbb{Q}_{p}}\ .
$$
 Gr\^ace \`a [F-R, \S4] on en d\'eduit une injection $M_{0}\hookrightarrow M$ telle que
 $$
\mathcal{M}_{0}=M_{0}\otimes_{A^{\dag}}\hat{A}\ ,\quad M_{0\mathbb{Q}_{p}}=M_{0}\otimes_{A^{\dag}}A^{\dag}_{\mathbb{Q}_{p}}\ ;
$$
ce qui est impossible d'apr\`es le th\'eor\`eme (5.4.3.5). D'o\`u le corollaire. $\square$

\vskip3mm
\noindent \textbf{Remarques (5.4.3.11)}. Le corollaire (5.4.3.10) et les relations (5.4.3.8)(5.4.3.9) peuvent se paraphraser en disant deux choses :\\
(i) Dans la base $u$ de  $\mathcal{M}_{0}$, $\phi_{\mathcal{M}_{0}}(\varphi_{can})$ est surconvergent et $\nabla_{\mathcal{M}_{0}}$ ne l'est pas: le corollaire (5.4.3.10) nous dit qu'on ne peut pas faire mieux, i.e. qu'il n'existe aucun rel\`evement $F_{\hat{A}}$ du Frobenius de $A_{0}$ ni aucune base de $\mathcal{M}_{0}$ dans laquelle \`a la fois  $\phi_{\mathcal{M}_{0}}(F_{\hat{A}})$ et $\nabla_{\mathcal{M}_{0}}$ soient surconvergents.\\
\textbf{(ii) Bien que $\mathcal{M}_{0}$ soit surconvergent au sens de Dwork} (i.e. existence d'un Frobenius $F_{\hat{A}}$ et d'une base de $\mathcal{M}_{0}$ tels que la matrice de $\phi_{\mathcal{M}_{0}}$ soit surconvergente, c'est-\`a-dire \`a coefficients dans $A^{\dag}$: avec les notations ci-dessus la matrice de $\phi_{1\mathcal{M}}$ peut \^etre surconvergente et pas celle de $\phi_{1\mathcal{M}}$, la diff\'erence provenant de la matrice de $\chi(F_{1},F_{2})$),\\
\textbf{la filtration par les pentes (5.4.3.3)  $\mathcal{M}_{0}\subset \mathcal{M}$ ne se rel\`eve pas en une filtration analogue de $M$, sinon $\mathcal{M}_{0}$ serait surconvergent au sens de Berthelot, ce qui n'est pas le cas}.\\
 Cette surconvergence au sens de Dwork de $\mathcal{M}_{0}$ permet cependant de lui appliquer la formule des traces de Monsky et d'en d\'eduire la m\'eromorphie de la fonction $L(X,\mathcal{M}_{0},t)$.

\newpage
\subsection*{5.4.4 Deuxi\`eme contre-exemple: les courbes modulaires}
Soient $k=\mathbb{F}_{p}, N\geqslant 3$ un entier premier \`a $p$, $Y_{0}(N)$ la courbe modulaire affine lisse sur $\mathbb{Z}_{p}$ espace de modules de courbes elliptiques avec structure de niveau de type $\Gamma_{0}(N)$. On note $A\rightarrow Y_{0}(N)$ la courbe elliptique universelle. Soient $Y_{0}(N)_{k}$ la fibre sp\'eciale de $Y_{0}(N)$, $X\subset Y_{0}(N)_{k}$ l'ouvert qui param\'etrise les courbes elliptiques ordinaires [F-C,V, \S7, p.192] et $U\subset Y_{0}(N)$ un sous-sch\'ema ouvert de fibre sp\'eciale $X$. On d\'esigne par $\mathcal{X}$ le compl\'et\'e formel de $U$ le long de $X$: c'est un $\mathbb{Z}_{p}$-sch\'ema formel ind\'ependant du choix de $U$. Consid\'erons le probl\`eme de module $\mathcal{A}^{ord}_{N}$ qui associe \`a tout $\mathbb{Z}_{p}$-sch\'ema localement noeth\'erien $S$ dans lequel $p$ est localement nilpotent, l'ensemble des classes d'isomorphismes
$$
\mathcal{A}^{ord}_{N}(S)=\{(B,\delta_{N})\}/ \simeq
$$
o\`u $B$ est un $S$-sch\'ema ab\'elien de dimension relative 1 muni d'une polarisation principale tel que toutes les fibres g\'eom\'etriques de $B\rightarrow S$ sont des courbes elliptiques ordinaires et $\delta_{N}$ est une structure de niveau $N$. Alors $\mathcal{A}^{ord}_{N}$ est ind-repr\'esentable par $\mathcal{X}$ [A-M, \S8]. La courbe elliptique formelle universelle $\mathcal{C}\rightarrow \mathcal{X}$ est la compl\'et\'ee formelle de l'image inverse de $A\rightarrow Y_{0}(N)$ par la fl\`eche $U\hookrightarrow Y_{0}(N)$; de plus $\mathcal{C}\rightarrow \mathcal{X}$ rel\`eve la courbe elliptique ordinaire universelle $f: C\rightarrow X$. Ainsi $E=R^{1}f_{rig\ast}(C/\mathcal{X})$ est un $F$-isocristal surconvergent sur $X$ de rang 2 [B 2, (2.3.8)(iii)] [Et 4, th\'eo 7]; notons $\mathcal{E}$ le $F$-isocristal convergent associ\'e \`a $E$ par le foncteur d'oubli
$$
F\mbox{-Isoc}^{\dag}(X/\mathbb{Q}_{p})\rightarrow F\mbox{-Isoc}(X/\mathbb{Q}_{p})\ .
$$
Le sous-$F$-isocristal unit\'e $\mathcal{E}_{0}$ de $\mathcal{E}$ est l'isocristal associ\'e, par la construction de Berthelot [B 2, 2.4], au $F$-cristal unit\'e de rang 1 $R^{1}f_{\textrm{\'et}\ast}(\mathbb{Z}_{p})\otimes_{\mathbb{Z}_{p}}\mathcal{O}_{X/\mathbb{Z}_{p}}$: ce $F$-cristal unit\'e correspond \`a une repr\'esentation [K 4] 
$$
\rho: \Pi_{1}(X, \overline{x})\rightarrow \mathbb{Z}_{p}^{\times}
$$
o\`u $\overline{x}$ est un point g\'eom\'etrique de $X$ fix\'e. D'apr\`es un th\'eor\`eme de Igusa [Ig], [K 4, 4.3] l'image par $\rho$ du groupe d'inertie en chaque point supersingulier est \'egale \`a $\mathbb{Z}_{p}^{\times}$ tout entier: en particulier $\rho$ n'est pas \`a monodromie locale finie [C, \S3], donc d'apr\`es Crew [C, theo 4.12] $\mathcal{E}_{0}$ ne provient pas d'un $F$-isocristal surconvergent par le foncteur d'oubli ci-dessus. Cependant, d'apr\`es un th\'eor\`eme de Brinon-Mokrane [Bri-Mo, th\'eo 1.1], le \guillemotleft rel\`evement canonique\guillemotright \ $\phi: \mathcal{X}\rightarrow \mathcal{X}$  du Frobenius de $X$, issu de la th\'eorie des sous-groupes canoniques de Abbes-Mokrane [A-M], fournit un Frobenius surconvergent
$$
\Phi_{\mathcal{E}_{0}}: \phi^{\ast}\mathcal{E}_{0}\rightarrow\mathcal{E}_{0}\ .
$$
\textbf{Ainsi les analogues du th\'eor\`eme (5.4.3.5), du corollaire (5.4.3.10) et des remarques (5.4.3.11) s'appliquent int\'egralement \`a ce nouvel exemple:} en particulier \\
(i) $\mathcal{E}_{0}$ est surconvergent au sens de Dwork mais pas au sens de Berthelot.\\
(ii) malgr\'e la surconvergence de $\mathcal{E}_{0}$ au sens de Dwork la filtration $\mathcal{E}_{0}\subset\mathcal{E}$ ne se rel\`eve pas en une filtration analogue de $E$.


\newpage

\begin{thebibliography}{1234567}
\markboth{\sc j.-y. etesse}{\sc Bibliographie}
\bibitem[A-M] {A. Abbes, A. Mokrane} A. Abbes, A. Mokrane: \textit{Sous-groupes canoniques et cycles \'evanescents $p$-adiques pour les vari\'et\'es ab\'eliennes}, Publ. Math. IHES 99 (2004), 117-162.
\bibitem[A] {Y. Amice} Y. Amice: \textit{Les nombres p-adiques}, Presses Universitaires de France (1975).
\bibitem[B 1]  {P. Berthelot 1}  P. Berthelot: \textit{Cohomologie rigide et th\'eorie de Dwork: le cas des sommes exponentielles},  Ast\'erisque \no 119-120 (1984), 17-49.
\bibitem[B 2]  {P. Berthelot 2} P. Berthelot: \textit{Cohomologie rigide et cohomologie rigide \`a supports propres}, Pr\'epublication 93-03 de Rennes (1996).
\bibitem[B-B-M]  {P. Berthelot , L. Breen, W. Messing} P. Berthelot , L. Breen, W. Messing: \textit{Th\'eorie de Dieudonn\'e cristalline II}, Lecture Notes in Math. 930, Springer (1982).
\bibitem[B-M] {P. Berthelot, W. Messing } P. Berthelot, W. Messing: \textit{Th\'eorie de Dieudonn\'e cristalline III : th\'eor\`emes d'\'equivalence et de pleine fid\'elit\'e}, The Grothendieck Festschrift, vol. 1, Progress in Math. 86, Birkha¬\"user (1990).
\bibitem[B-G-R]{S. Bosch, U. Guntzer, R. Remmert} S. Bosch, U. G\"untzer, R. Remmert:  \textit{Non-archimedean analysis}, Grundlehren der Math. Wissenschaften 261, Springer Verlag (1984).
\bibitem[Bour]{N. Bourbak}N. Bourbaki: \textit{Alg\`ebre} [A] chap. I \`a VII; \textit{Alg\`ebre commutative} [AC] chap. I \`a X.
\bibitem[Bri-Mo]{O. Brinon, F. Mokrane}O. Brinon, F. Mokrane: \textit{Surconvergence de la monodromie p-adique des familles universelles de vari\'et\'es ab\'eliennes ordinaires}, preprint 13 oct. 2010.
\bibitem[C]{R. Crew}R. Crew: \textit{$F$-isocrystals and $p$-adic representations} in \textit{Algebraic Geometry, Bowdoin 1985}, Proceedings of Symposia in Pure Math., Vol. 46, AMS (1987), 111-138.
\bibitem[De$\ell$ 1]{P. Deligne 1}P. Deligne: \textit{La conjecture de Weil II}, Pub. Math. IHES 52, (1980), 137-252.
\bibitem[De$\ell$ 2]{P. Deligne 2}P. Deligne: \textit{Cristaux ordinaires et coordonn\'ees canoniques} in \textit{Surfaces alg\'ebriques}, Lecture Notes in Math. 868, Springer (1981).
\bibitem[De$\ell$ 3]{P. Deligne 3}P. Deligne: \textit{Rapport sur la formule des traces} in \textit{SGA 4 1/2, Cohomologie Etale}, Lecture Notes in Math. 569, Springer (1977), 76-109.
\bibitem[De$\ell$ 4]{P. Deligne 4}P. Deligne: \textit{Fonctions $L$ modulo $\ell^{n}$ et modulo $p$} in \textit{SGA 4 1/2, Cohomologie Etale}, Lecture Notes in Math. 569, Springer (1977), 110-128.
\bibitem[Dw 1]{B. Dwork 1}B. Dwork: \textit{On the rationality of the zeta function of an algebraic variety}, Amer. J. Math., vol 88, (1960), 631-648.
\bibitem[Dw 2]{B. Dwork 2}B. Dwork: \textit{p-adic Cycles}, Pub. Math. IHES \no 37 (1969), 27-115.
\bibitem[Dw 3]{B. Dwork 3}B. Dwork: \textit{On Hecke polynomials}, Invent. Math. 12, (1971), 249-256.
\bibitem[Dw 4]{B. Dwork 4}B. Dwork: \textit{Normalized period matrices I}, Annals of Math. 94, \no2, (1971), 337-388.
\bibitem[Dw 5]{B. Dwork 5}B. Dwork: \textit{Normalized period matrices II}, Annals of Math. 98, \no1 (1973), 1-57.
\bibitem[Dw 6]{B. Dwork 6}B. Dwork: \textit{Lectures on p-adic differential equations}, Springer Grundlehren 253 (1982).
\bibitem[Dw-S]{B. Dwork, S. Sperber}B. Dwork, S. Sperber: \textit{Logarithmic decay and overconvergence of the unit root and associated zeta functions}, Annales Scient. Ec. Norm. Sup., 4\`eme s\'erie, t. 24, (1991), 575-604.
\bibitem[EGA]{A. Grothendieck, J. Dieudonne}A. Grothendieck, J. Dieudonn\'e: \textit{El\'ements de G\'eom\'etrie Alg\'ebrique}: Chap. I, Springer Grundlehren 166; Chap. II, III, IV, Pub. Math. IHES \no 8, 11, 17, 20, 24, 28, 32.
\bibitem[Eis]{D. Eisenbud}D. Eisenbud: \textit{Commutative Algebra with a View Toward Algebraic Geometry}, Graduate Texts in Math. 150, Springer, (1996).
\bibitem[E$\ell$]{R. Elkik}R. Elkik: \textit{Solutions d'\'equations \`a coefficients dans un anneau hens\'elien}, Annales Scient. Ec. Norm. Sup., 4\`eme s\'erie, t. 6, (1973), 553-604.
\bibitem[Et 1]{J.-Y. Etesse 1}J.-Y. Etesse: \textit{Rationalit\'e et valeurs de fonctions L en cohomologie cristalline}, Annales Inst. Fourier, t. 38, fasc. 4, (1988), 33-92.
\bibitem[Et 2]{J.-Y. Etesse 2}J.-Y. Etesse: \textit{Rel\`evement de sch\'emas ab\'eliens, F-cristaux et fonctions L}, J. reine angew. Math. 535, (2001), 51-63.
\bibitem[Et 3]{J.-Y. Etesse 3}J.-Y. Etesse: \textit{Rel\`evement de sch\'emas et alg\`ebres de Monsky-Washnitzer: th\'eor\`emes d'\'equivalence et de pleine fid\'elit\'e}, Rendiconti Sem. Mat. Univ. Padova, Vol. 107, (2002), 111-138.
\bibitem[Et 4]{J.-Y. Etesse 4}J.-Y. Etesse: \textit{Descente \'etale des F-isocristaux surconvergents et rationalit\'e des fonctions L de sch\'emas ab\'eliens}, Annales Scient. Ec. Norm. Sup., 4\`eme s\'erie, t. 35, (2002), 575-603.
\bibitem[Et 5]{J.-Y. Etesse 5}J.-Y. Etesse: \textit{Introduction to L- functions of F-isocrystals}, in \textit{Geometric Aspects of Dwork Theory}, Vol. II, de Gruyter, (2002), 701-710.
\bibitem[Et 6]{J.-Y. Etesse 6}J.-Y. Etesse: \textit{Rel\`evement de sch\'emas et alg\`ebres de Monsky-Washnitzer: th\'eor\`emes d'\'equivalence et de pleine fid\'elit\'e II}, Rendiconti Sem. Mat. Univ. Padova, Vol. 122, (2009), 205-234.
\bibitem[Et 7]{J.-Y. Etesse 7}J.-Y. Etesse: \textit{Images directes I: Espaces rigides analytiques et images directes}, hal-00425909/arXiv: 0910.4433.
\bibitem[Et 8]{J.-Y. Etesse 8}J.-Y. Etesse: \textit{Images directes II: F-isocristaux convergents}, hal-00425919/arXiv: 0910.4434.
\bibitem[Et 9]{J.-Y. Etesse 9}J.-Y. Etesse: \textit{Images directes III: F-isocristaux surconvergents}, hal-00425922/arXiv: 0910.4435.
\bibitem[Et 10]{J.-Y. Etesse 10}J.-Y. Etesse: \textit{Fonctions L en g\'eom\'etrie rigide II: F-(iso)cristaux et conjecture de Katz}, preprint.
\bibitem[Et 11]{J.-Y. Etesse 11}J.-Y. Etesse: \textit{Fonctions L en g\'eom\'etrie rigide III: Sch\'emas ab\'eliens ordinaires}, preprint.
\bibitem[E-LS 1]{J.-Y. Etesse, B. Le Stum 1}J.-Y. Etesse, B. Le Stum: \textit{Fonctions L associ\'ees aux F-isocristaux surconvergents I: Interpr\'etation cohomologique}, Math. Annalen 296, (1993), 557-576.
\bibitem[E-LS 2]{J.-Y. Etesse, B. Le Stum 2}J.-Y. Etesse, B. Le Stum: \textit{Fonctions L associ\'ees aux F-isocristaux surconvergents II: Z\'eros et p\^oles unit\'es}, Invent. Math. 127, (1997), 1-31.
\bibitem[F-C]{G. Faltings, C.-L. Chai}G. Faltings, C.-L. Chai: \textit{Degeneration of Abelian Varieties}, Ergebnisse der Mathematik und ihrer Grenzgebiete 3. Folge, Band 22, Springer (1990).
\bibitem[F-R]{D. Ferrand, M. Raynaud}D. Ferrand, M. Raynaud: \textit{Fibres formelles d'un anneau local noeth\'erien}, Annales Scient. Ec. Norm. Sup., 4\`eme s\'erie, t. 3, (1970), 295-311.
\bibitem[G-K 1]{E. Grosse-Klonne 1}E. Gro\ss e-Kl\" onne: \textit{de Rham-Kohomologie in der rigiden Analysis}, Preprintreihe der Universit\" at M\" unster SFB 478, Heft 39 (1999).
\bibitem[G-K 2]{E. Grosse-Klonne 2}E. Gro\ss e-Kl\" onne: \textit{Rigid analytic spaces with overconvergent stucture sheaf}, Journal f\" ur die reine und angewandte Math. 519,(2000), 73-95.
\bibitem[G 1]{A. Grothendieck 1}A. Grothendieck: \textit{Fondements de la G\'eom\'etrie Alg\'ebrique}, Extraits du S\'eminaire Bourbaki 1957-1962, Secr\'etariat Math\'ematique (1962).
\bibitem[G 2]{A. Grothendieck 2}A. Grothendieck: \textit{Formule de Lefschetz et rationalit\'e des fonctions L},  S\'eminaire Bourbaki 279, dans \textit{Dix expos\'es sur la cohomologie des sch\'emas}, Masson, North-Holland, (1968).
\bibitem[Ig]{J.-I. Igusa}J.-I. Igusa: \textit{On the algebraic theory of elliptic modular functions}, J. Math. Soc. Japan 20 (1968), 96-106.
\bibitem[K 1]{N. Katz 1}N. Katz: \textit{Travaux de Dwork}, S\'eminaire Bourbaki 409, Lecture Notes in Math. 383, Springer (1972).
\bibitem[K 2]{N. Katz 2}N. Katz: \textit{Slope filtration of F-crystals}, Ast\'erisque 63 (1979), 113-163.
\bibitem[K 3]{N. Katz 3}N. Katz: \textit{Appendice \`a Cristaux ordinaires et coordonn\'ees canoniques} in \textit{Surfaces alg\'ebriques}, Lecture Notes in Math. 868, Springer (1981).
\bibitem[K 4]{N. Katz 4}N. Katz: \textit{p-adic properties of modular schemes and modular forms} in \textit{Modular Functions of One Variable III}, Lecture Notes in Math. 350, Springer (1973).
\bibitem[Ked 1]{K. Kedlaya 1}K. Kedlaya: \textit{Finiteness of rigid cohomology with coefficients}, Preprint, arxiv: math.AG/0208027. Duke Math. J. 134 (2006), 15-97.
\bibitem[Ked 2]{K. Kedlaya 2}K. Kedlaya: \textit{Semistable reduction for overconvergent F-isocrystals, I: Unipotence and logarithmic extensions}, Preprint, arxiv: math.NT/0405069 v3, 24 Jul 2005. Compositio Math. 143 (2007), 1164-1212.
\bibitem[Ma] { H. Matsumura}  H. Matsumura: \textit{Commutative ring theory}, Cambridge Studies in Advanced Mathematics 8 (1997).
\bibitem[Mi] { J.-S. Milne} J.-S. Milne: \textit{Etale cohomology}, Princeton University Press (1980).
\bibitem[Mo 1] { P. Monsky 1} P. Monsky: \textit{p-Adic Analysis and Zeta Functions}, Kinokuniya Book-Store Co,Tokyo (1970). 
\bibitem[Mo 2] { P. Monsky 2} P. Monsky: \textit{Formal Cohomolgy III}, Annals of Math. 93, \no 2 (1971), 315-343. 
\bibitem[M-W] { P. Monsky, G. Washnitzer} P. Monsky, G. Washnitzer: \textit{Formal Cohomolgy I}, Annals of Math. 88, \no 2 (1968), 181-217. 
\bibitem[vdP] { M. van der Put} M. van der Put: \textit{The cohomology of Monsky and Washnitzer}, Bulletin de la SMF, m\'emoire \no23, t. 114/fasc. 2 (1986), 33-60.
\bibitem[R] { M. Raynaud } M. Raynaud: \textit{Anneaux locaux hens\'eliens}, Lecture Notes in Math. 169, Springer (1970).
\bibitem[Ro] {A.M. Robert} A.M. Robert: \textit {A Course in p-adic Analysis}, Graduate Texts in Mathematics 198, Springer (2000).
\bibitem[S 1] { J.-P. Serre 1} J.-P. Serre: \textit{Corps locaux}, Hermann (1968).
\bibitem[S 2] { J.-P. Serre 2} J.-P. Serre: \textit{Endomorphismes compl\`etement continus des espaces de Banach $p$-adiques}, Pub. Math. IHES \no 12 (1962), p. 69-85.
\bibitem[Sp ] { S. Sperber} S. Sperber: \textit{Congruence properties of the Hyperkloostermann sum}, Compositio Math. 40 (1980), 3-33.
\bibitem[SGA 1] { A. Grothendieck} A. Grothendieck: \textit{Rev\^etements \'etales et goupe fondamental}, Lecture Notes in Math. 224, Springer (1971).
\bibitem[SGA 3] { M. Demazure, A. Grothendieck} M. Demazure, A. Grothendieck: \textit{Sch\'emas en groupes}, Lecture Notes in Math. 151, 152, 153, Springer (1970).
\bibitem[SGA 4] {M. Artin, A. Grothendieck, J.-L. Verdier}M. Artin, A. Grothendieck, J.-L. Verdier: \textit{Th\'eorie des topos et cohomologie \'etale des sch\'emas}, Lecture Notes in Math. 269, 270, 305, Springer (1972, 1973).
\bibitem[SGA 7, II] {P. Deligne, N. Katz}P. Deligne, N. Katz: \textit{Groupes de Monodromie en G\'eom\'etrie Alg\'ebrique}, Lecture Notes in Math. 340, Springer (1973).
\bibitem[Shi 1] { A. Shiho 1} A. Shiho: \textit{Crystalline Fundamental Groups II- Log Convergent Cohomology and Rigid Cohomology}, J. Math. Sci. Univ. Tokyo 9 (2002), 1-163.
\bibitem[Shi 2] { A. Shiho 2} A. Shiho: \textit{Relative Log Convergent Cohomology and Relative Rigid Cohomology I}, arXiv: 0707.1742v1 [math.NT] 12 Jul 2007.
\bibitem[Shi 3] { A. Shiho 3} A. Shiho: \textit{Relative Log Convergent Cohomology and Relative Rigid Cohomology II}, arXiv: 0707.1743v1 [math.NT] 12 Jul 2007.
\bibitem[W 1] { D. Wan 1} D. Wan: \textit{Meromorphic continuation of L-functions of p-adic representations}, Annals of Math. 143 (1996), 469-498.
\bibitem[W 2] { D. Wan 2} D. Wan: \textit{Dwork's conjecture on unit-root zeta functions}, Annals of Math. 150 (1999), 867-927.
\bibitem[W 3] {D. Wan 3} D. Wan: \textit{Higher rank case of Dwork's conjecture}, J. Amer. Math. Soc. 13 (2000), 807-852.
\bibitem[W 4] { D. Wan 4} D. Wan: \textit{Rank one case of Dwork's conjecture}, J. Amer. Math. Soc. 13 (2000), 853-908.
 \end{thebibliography}
\end{document}